\documentclass[11pt]{article}
\usepackage{a4wide,
amssymb,latexsym,epsfig,amsmath,amsfonts,mathtools,array,graphicx,verbatim}
\usepackage[english]{babel}
\newtheorem{thm}{Theorem}
\newtheorem{lem}[thm]{Lemma}
\newtheorem{prop}[thm]{Proposition}

\newcommand{\sect}[1]{\section{#1}\setcounter{thm}{0}\setcounter{equation}{0}}
\renewcommand{\thesection}{\arabic{section}.}

\renewcommand{\theequation}{\thesection\arabic{equation}}
\newcommand{\be}{\begin{equation}}
\newcommand{\ee}{\end{equation}}
\newcommand{\ba}{\begin{array}}
\newcommand{\ea}{\end{array}}

\newcommand{\noi}{\noindent}
\newcommand{\pa}{\partial}
\newcommand{\na}{\nabla}
\newcommand{\nt}{{\widetilde \nabla}}
\newcommand{\De}{\Delta}
\newcommand{\ti}{\times}

\newcommand{\cpl}{{\mathbb C}}
\newcommand{\HP}{{\mathbb H}}

\newcommand{\eee}{{\bf e}}

\newcommand{\eps}{\varepsilon}
\newcommand{\ro}{{r\!o}}
\newcommand{\roh}{{\widehat{r\!o}}}

\newcommand{\vi}{\varphi}

\newcommand{\vih}{{\widehat \varphi}}
\newcommand{\vii}{\psi}
\newcommand{\ppsi}{{\psi\kern-0.65em \psi}}
\newcommand{\ze}{\zeta}

\newcommand{\ga}{\gamma}
\newcommand{\al}{{\alpha}}
\newcommand{\la}{\lambda}
\newcommand{\ka}{\kappa}
\newcommand{\si}{\sigma}

\newcommand{\lab}{\overline{ \lambda}}

\newcommand{\lala}{{\lambda\kern-0.51em \lambda}}
\newcommand{\Om}{\Omega}
\newcommand{\Omt}{{\widetilde{\Omega}}}

\newcommand{\FF}{\overline{ F'}}
\newcommand{\EE}{\overline{ E}}

\newcommand{\disp}{\displaystyle}
\newcommand{\ints}{{\mathbb Z}}

\newcommand{\nat}{{\mathbb N}}

\newcommand{\R}{{\rm R\kern-0.1em e\,}}
\newcommand{\Rc}{{\rm R\kern-0.1em e_c}}
\newcommand{\RE}{{\rm R\kern-0.1em e_E}}
\newcommand{\Rx}{{\rm R\kern-0.1em e_E^x}}
\newcommand{\Ry}{{\rm R\kern-0.1em e_E^y}}
\newcommand{\qed}{\hfill\hspace*{2mm}\hfill $\Box$}

\newcommand{\real}{{\mathbb R}}
\newcommand{\ra}{\rightarrow}
\newcommand{\rd}{{\rm d}}

\newcommand{\du}{{\delta\! u}}
\newcommand{\dv}{{\delta\! v}}
\newcommand{\bb}{{\bf b}}
\newcommand{\aaa}{{\bf a}}
\newcommand{\aaat}{{\widetilde{\bf a}}}

\newcommand{\xx}{{\bf x}}

\newcommand{\nn}{{\bf n}}
\newcommand{\ttt}{{\bf t}}

\newcommand{\BC}{{\cal B}}
\newcommand{\BtC}{{\widetilde{{\cal B}}}}
\newcommand{\CC}{{\cal C}}

\newcommand{\FC}{{\cal F}}
\newcommand{\FCt}{\widetilde{\cal F}}
\newcommand{\HC}{{\cal H}}
\newcommand{\IC}{{\cal I}}
\newcommand{\LC}{{\cal L}}
\newcommand{\KC}{{\cal K}}
\newcommand{\KCt}{\widetilde{\cal K}}
\newcommand{\RC}{{\cal R}}
\newcommand{\TC}{{\cal T}}

\newcommand{\VC}{{\cal V}}

\newcommand{\HH}{{\bf H}}
\renewcommand{\FF}{{\bf F}}

\newcommand{\Ht}{\widetilde{H}}
\newcommand{\HHt}{\widetilde{{\bf H}}}

\newcommand{\DD}{{\bf D}}

\newcommand{\Ih}{{\widehat I}}

\newcommand{\Ct}{{\widetilde C}}

\newcommand{\de}{\delta}
\newcommand{\dt}{{\tilde \delta}}
\newcommand{\ati}{{\tilde a}}

\newcommand{\at}{{\widetilde \alpha}}

\newcommand{\ct}{{\tilde c}}

\newcommand{\lt}{{\tilde l}}
\newcommand{\pt}{{\widetilde p}}
\newcommand{\qt}{{\widetilde q}}
\newcommand{\ut}{{\tilde u}}
\newcommand{\vt}{{\tilde v}}

\newcommand{\pst}{\widetilde{\psi}}

\newcommand{\Lt}{{\widetilde L}}

\newcommand{\ph}{{\widehat p}}

\newcommand{\qb}{{\overline{q}}}
\newcommand{\zb}{{\overline{z}}}
\newcommand{\fb}{{\overline{f}}}
\newcommand{\hb}{{\overline{h}}}
\newcommand{\Fb}{{\overline{F}}}

\newcommand{\ARb}{{\overline{A_R}}}
\newcommand{\QRb}{{\overline{Q_R}}}
\newcommand{\Qib}{{\overline{Q_\infty}}}

\begin{document}
\title{Geophysical intensity problems: the axisymmetric case} 
\author{Ralf Kaiser \\[2ex]      
\small 
Fakult\"at f\"ur Mathematik, Physik und Informatik\\[1ex]
 \small Universit\"at Bayreuth, D-95440 Bayreuth, Germany\\[1ex]
          \small ralf.kaiser@uni-bayreuth.de} 
\date{December, 2025}
\maketitle
%%%
%%
%
%
%     Abstract
%
%
%
\begin{abstract}
Considering the earth or any other celestial body the main sources of the gravitational as well as of the magnetic field lie inside the body. Above the surface both fields are in good approximation harmonic vector fields determined by their values at the body's surface or any other surface enclosing the body. The intensity problem seeks to determine harmonic vector fields vanishing at infinity and with prescribed intensity of the field at the surface. This problem constitutes a nonlinear boundary value problem, whose general solvability is not yet established. From a practical point of view uniqueness of solutions is the even more urgent problem: intensity data can often be obtained with greater accuracy and at less cost than field directions. So, when modelling the (especially magnetic) field of the earth the question arises whether intensity data are sufficient to obtain a reliable field model.

In this paper {\em axisymmetric} harmonic fields $\HH$ outside the unit sphere $S^2 \subset \real^3$ are studied and, given an axisymmtric H\"older continuous intensity function $I\neq 0$ on $S^2$, the existence of infinitely many solutions of the intensity problem is proved. These solutions can more precisely be characterized as follows: fix a number $\de \in \nat\setminus \{1 \}$ and a meridional plane $M$ through the symmetry axis  $S\!A$, and in $M$  a unit circle $S^1$ (symmetric with respect to $S\!A$) and, furthermore, $2\, N$, $N \in \nat_0$, points $z_n \in M$ (symmetric with respect to $S\!A$, avoiding $S\!A$, and outside $S^1$), then the existence of an (up to a sign) unique harmonic field $\HH$ is established that vanishes at (the axisymmetric circles piercing $M$ at) $z_n$ and nowhere else,
that has intensity $I$ at $S^2$ and (exact) decay order $\de$ at infinity.

The proof is based on the solution of a nonlinear elliptic equation with discontinuous coefficients, which are, moreover, singular at the symmetry axis. Its combination with fixed boundary conditions was the basis of a recent treatment of the ``geomagnetic direction problem'' \cite{KR22}. Here we have instead natural boundary conditions, which provide less information, and which require, therefore, in part new solution techniques and sharper estimates. 


\noindent MSC-Classification (2020): 35J65, 86A25, 86A30.
\end{abstract}
%
%
%
%
%
%
%
%            Einleitung
%
%
%
%
%
%
\sect{Introduction}
In the theory of harmonic vector fields the standard decomposition of the field vector at the boundary in tangential and normal components, which each lead to well posed boundary value problems, is not always the best choice, especially not in geophysical applications.
Considering measurement devices the distinction between direction and intensity makes in fact a better choice, since these are the quantities that are for historical, technical, or economical reasons quite differently available. Concerning the magnetic field of the earth, most of the time in the history of magnetic field measurements only directions could be detected \cite{MM83}. Meanwhile the complete field vector can be determined; however, measuring only the intensity of the magnetic field can be easier and cheaper than determining the complete field vector. In situations, where the attitude of the magnetometer is difficult to control, intensity data are the only or the more reliable magnetic data. Early satellite surveys, for example, provided only intensity data, and even nowadays directional information is still susceptible to errors in determining the attitude of the satellite in orbit. Similar problems arise at sea when the field is measured by magnetometers towed far behind oceanographic ships to escape the ship's own magnetic field (for more information about the significance of the intensity problem we refer to \cite{B68,PG90,AOP04,GH07} and references therein). Concerning the gravitational field of the earth, which is dominated by its monopole contribution, measuring directions is obviously not the most sensitive way to uncover the field's fine structure. Intensity data from satellite measurements, however, as in the magnetic case are nowadays widely available and the corresponding intensity problem is known as the ``{\em fixed} geodetic boundary value problem" \cite{KP72,H89}. This problem  has to be distinguished from the more general (and more ambitious) geodetic boundary value problem that seeks to determine both the gravity field and the shape of the earth from the complete field vector and the gravitational potential at the earth's surface (for more information we refer to \cite{H76,TM12}).

Here we consider an ideal situation, where the intensity data are known all over a spherical surface $S^2$ of radius $1$ and where the effects of rotation or of additional sources of the harmonic field in the exterior space $E \subset \real^3$ of the unit ball are neglected. Given a positive continuous function $I: S^2 \ra \real_+$ (the ``intensity function'', see Fig.\ 1 for an example) and $\de \in \nat\setminus \{1\}$ (the ``decay order'' of the harmonic field at infinity) the {\bf intensity problem} $P_I$ asks for all vector fields $\HH \in C^1 (E) \cap C (\EE)$ satisfying the conditions
\begin{equation}\label{1.1}
\left.
\begin{array}{cc}
\nabla \times \HH = 0,\ \ \nabla \cdot \HH = 0 & \mbox{ in }E, \\[1ex]
|\HH (\xx)| = O(|\xx|^{-\de}) & \mbox{ for } |\xx|\ra \infty, \\[1ex]
|\HH | = I & \mbox{ on } S^2 .
\end{array} \right\}
\end{equation}
The magnetic and gravitational problems differ in that magnetic fields have no monopole part, a condition that can be taken into account by the choice $\de = 3$, whereas gravitational fields (due to the positive mass density) are only admissible if 
\be\label{1.2}
\HH \cdot \xx \neq 0  \qquad \mbox{ on } S^2 ,
\ee
which requires $\de = 2$ in (\ref{1.1})$_2$ .

\begin{figure}
\begin{center}
\includegraphics[width=0.6\textwidth, angle=-90]{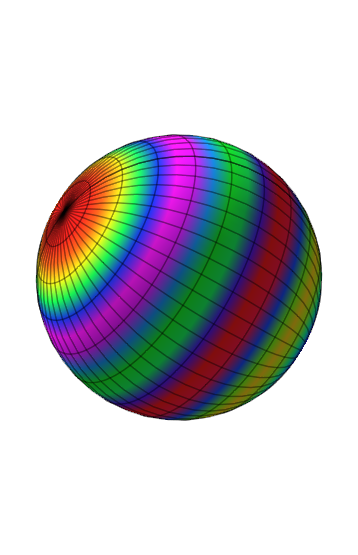}
\caption{Illustration of an axisymmetric intensity function
on a sphere with intensity coded by color.}
\end{center}
\end{figure}

For the intensity problem there is obviously no linear relation between solution and boundary data, which means that the usual solution techniques for the standard boundary value problems such as integral equation or Hilbert space methods are not applicable. In consequence there are so far only a few results concerning existence and uniqueness: in the geodetic case, when condition (\ref{1.2}) holds additionally, a simple argument based on maximum principles proves uniqueness \cite{B68}; in the general (magnetic) case non-uniqueness is known by example \cite{B70}. Adding some directional information, such as the knowledge of the ``dip equators'' (loci on $S^2$, where the normal component of the field vector vanishes), can, however, restore uniqueness \cite{KHL97}. Concerning existence there are in the geodetic case results for small data, i.e.\ here data close in some norm to those of the monopole field \cite {BS83,J87,SS89}; in the general case there are viscosity solutions for a related problem inspired by the intensity problem \cite{DDO06}. These solutions are simultaneously solutions of the intensity problem if they happen to satisfy condition (\ref{1.2}) (which is conjectured by the authors of \cite{DDO06}). In general, however, this is not the case as we will see below.  

In the following we will present an (almost) complete solution of the {\em axisymmetric} intensity problem. The method is inspired by the solution of the two-dimensional version of this problem \cite{B68}, which used methods of complex analysis. Axisymmetry leads - in cylindrical coordinates $\rho$, $\ze$ in a meridional plane $\theta = const$ - likewise to a two-dimensional problem, which, although complicated by the coordinate singularity at $\rho = 0$, is amenable to a complex formulation. With $(H_\ze, H_\rho)$ representing the nonvanishing components of an axisymmetric harmonic field we make the ansatz
\be \label{1.3}
H_\ze (\ze , \rho) + i\, H_\rho (\ze, \rho) = h(\ze + i \rho) \exp \Big(\frac{1}{2} \big(p(\ze, \rho) + i\, q(\ze, \rho)\big)\Big)
\ee
with a given complex analytic function $h$ representing the zeroes and the asymptotic behaviour of the harmonic field, and ``correction functions'' $p$ and $q$. This ansatz is the same as in \cite{KR22} for the solution of the axisymmetric direction problem, where we derived an elliptic boundary value problem for the direction-related variable $q$. An analogous boundary value problem for the intensity related variable $p$ does not exist. Instead, we derive for a modified variable $q'$ the following boundary value problem in $A_\infty := \{(\ze, \rho) \in \real^2 : \ze^2 + \rho^2 > 1 \} \subset \real^2$:
\begin{equation}\label{1.4}
\left.
\begin{array}{rll}
- \Delta\, q' & = \disp - \pa_\ze \Big( \frac{1}{\rho} \cos (q' - \Om)\Big) + \pa_\rho \Big(\frac{1}{\rho} \sin (q' - \Om) \Big)  & \mbox{ in } A_\infty , \\[2ex]
\nn \cdot \na q' & = \disp  \frac{n_\ze}{\rho} \Big(\cos (q' - \Om) - 1\Big) -\frac{n_\rho}{\rho} \sin (q' - \Om) & \mbox{ on } S^1 , 
\end{array} \right\}
\end{equation}
where $\nn$ denotes the exterior normal at the unit circle $S^1$, and 
the ``data function'' $\Om$ depends on $h$ and the intensity function $I$.
Problem (\ref{1.4}) differs from the corresponding direction problem only in the ``natural'' boundary condition (\ref{1.4})$_2$. In a weak formulation this is implemented by test functions, which are free at the boundary in contrast to test functions vanishing at the boundary in the direction problem. It is this loss of information at the boundary that requires additional measures in the solution strategy for problem (\ref{1.4}), which is in principle the same as for the direction problem: the general solution of a linearized version of (\ref{1.4}) by applying a Lax-Milgram-type solution criterion and, subsequently, setting up a suitable iteration procedure in $L^p$-spaces, to which Schauder's fixed point principle can be applied. At the linear level it was a major problem in \cite{KR22} to obtain control by the left-hand side
over the second term on the right-hand side in (\ref{1.4})$_1$, which is not small and which behaves near the symmetry axis $\{\rho =0\}$ like a second order term. By the use of Hardy-type inequalities, weighted by some power of the distance to the symmetry axis, and by introducing further parameters, exploiting the fine structure of the right-hand side in (\ref{1.4})$_1$, a subtle balance of these parameters allowed, finally, to gain this control. 

For the intensity problem the above mentioned loss of information is not complete. Solutions and test functions are antisymmetric with respect to the symmetry axis and hence vanish there. A formulation of the problem in polar coordinates is now more appropriate (section 2); it introduces, however, additional (lower order) terms, which are not small and which change the problem to one of Fredholm-type. Uniqueness of the homogeneous problem is now essential (section 3) and proves together with a careful choice of (an enlarged set of) parameters existence in the linearized problem (section 4).

At the nonlinear level the same iteration, which worked in the direction problem, works here as well (section 5), since the essential Hardy-type estimates continue to hold in the present situation with less boundary information though with different constants (appendix A). Another essential ingredient for the working of the iteration is a global bound on the solution, which in the direction problem was a byproduct of the coercivity estimates for Lax-Milgram. Here, additionally, Fredholm constants are involved, which require spectral estimates of the associated Fredholm operator (appendix C).  

Regularity of the solution follows in principle as in the direction problem; only in two respects additional effort is necessary: at the boundary it needs some care to recover the prescribed intensity from the free boundary condition (section 6) and at the axis the evaluation of some condition on the intensity function requires some non-standard H\"older-type estimates (appendix B).

The solution procedure explained so far works in bounded regions. In the unbounded case we solve the intensity problem in a sequence of bounded regions with artificial intensity function at the exterior boundary. The global bound and a ``diagonal argument'' as in \cite{KR22} are then the crucial ingredients to justify the passage to the unbounded region establishing thus the solution of the intensity problem in exterior space.   

Finally, we have a look on the related problem \cite{DDO06} mentioned above and show by example the difference of both problems (section 7). 
\sect{Mathematical setting and results} 
This section provides the mathematical setting for our treatment of the intensity problem, derives the system (\ref{1.4}), and presents our results in the theorems 2.1 -- 2.4. Where this section overlaps with the corresponding section in \cite{KR22} it is a condensed version stressing the differences and omitting most of the motivations and comments given in \cite{KR22}.

Axisymmetric harmonic vector fields $\HH$, expressed in cylindrical coordinates $(\rho,\theta, \zeta)$, have just two nontrivial components $H_\rho$ and $H_\ze$, depending on $\rho$ and $\ze$, which satisfy the system  
\be \label{2.1}
\left. \ba{c} \pa_\ze H_\rho - \pa_\rho H_\ze = 0 , \\[1ex] 
\disp \pa_\ze H_\ze + \pa _\rho H_\rho + \frac{1}{\rho}\, H_\rho = 0 . 
\ea \right\}
\ee
For the time being $\HH$ is defined on the half-plane $\HP = \{(\ze, \rho) \in \real^2: \rho >0\}$ bounded by the symmetry axis $\{ \rho = 0\}$ but can be extended to $\real^2$ by (anti)symmetric continuation: 
\be \label{2.2}
\left.
\ba{l}
H_\ze (\ze, -\rho) := H_\ze (\ze, \rho) , \\[1ex]
H_\rho (\ze, -\rho) := - H_\rho (\ze, \rho) , 
\ea \right\}
\qquad (\ze, \rho) \in \HP .
\ee
Note that (\ref{2.1})$_2$ implies (if defined) $H_\rho (\ze, 0) = 0$. 
On $\real^2$ also polar coordinates $(r, \vi) \in (0, \infty) \ti (-\pi , \pi]$ with basis vectors $\eee_r$ and $\eee_\vi$ are useful.
They are related to $(\ze, \rho)$ by 
\be \label{2.3}
\ze = r \cos \vi\, , \quad \rho = r \sin \vi
\ee
(see Fig.\ 2). 
\begin{figure}
\begin{center}
\includegraphics[width=0.45\textwidth]{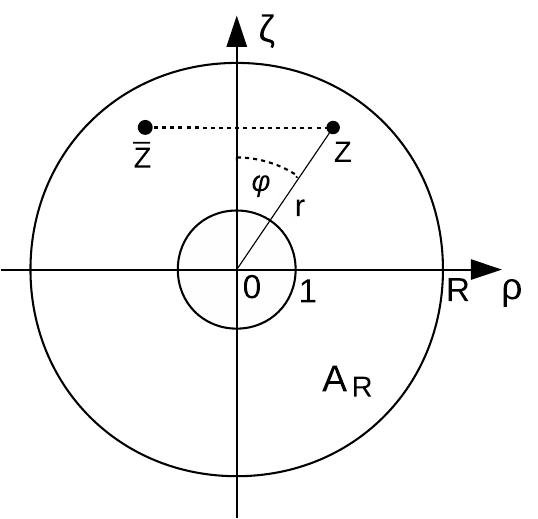}
\caption{Various coordinates in the meridional cross section $A_R$.}
\end{center}
\end{figure}
In these coordinates condition (\ref{2.2}) takes the form
\be \label{2.4}
\left.
\ba{l}
\Ht_\ze (r, -\vi) = \Ht_\ze (r, \vi) , \\[1ex]
\Ht_\rho (r, -\vi) = - \Ht_\rho (r, \vi), 
\ea \right\}
\ee
where $\HH(r \cos \vi, r \sin \vi) =: \HHt (r, \vi)$. 
Finally, to take advantage of the exponential mapping in complex analysis relating Cartesian and polar coordinates, a complex formulation is useful. 
With the identifications
\be \label{2.5}
\Re z := \ze\, , \qquad \Im z := \rho\, , \qquad \Re f := H_\ze \, ,\qquad \Im f := - H_\rho\, ,
\ee
eqs.\ (\ref{2.1}) then take the form 
\be \label{2.6}
\pa_\zb f - \frac{1}{2}\, \frac{f - \fb}{z - \zb} = 0
\ee
and the symmetry condition (\ref{2.2}) amounts to 
\be \label{2.7}
f(z, \zb) = \fb (\zb , z)\, . 
\ee
We made there use of the Wirtinger derivatives $\pa_z := \frac{1}{2} (\pa_\ze - i \pa_\rho)$ and $\pa_\zb := \frac{1}{2} (\pa_\ze + i \pa_\rho)$
acting on functions $f: \cpl \ti \cpl \ra \cpl$, $(z, \zb) \mapsto f(z, \zb)$.

Intensity functions $I : S_1 \ra \real_+$, parametrized by the polar angle $\vi \in (-\pi, \pi]$, are called symmetric, if $I(\vi) = I(-\vi)$, and fields $H$, $\Ht$, and complex functions $f$ are called symmetric, if they satisfy (\ref{2.2}), (\ref{2.4}), and (\ref{2.7}), respectively.\footnote{ We skip henceforth the upper index $1$ at $S$ indicating the dimension of the ``sphere'', and introduce a lower index indicating the radius of the circle/sphere.} The {\em axisymmetric} intensity problem now reads:

\hspace{1ex}

\noi\textbf{Problem $P_I(A_\infty)$}. {\it Let $I \in C(S_1, \real_+)$ be a positive symmetric function  and $\delta \in \nat \setminus \{1\}$. Determine all symmetric solutions $\HH \in C^1 (A_\infty) \cap C(\overline{A_\infty})$ of {\rm (\ref{2.1})} with decay order $\de$ and boundary condition $|\HH|\big|_{S_1} = I$ .
}

\hspace{1ex}

\noi Introducing an additional intensity function $\Ih$ on the circle $S_R$ of radius $R > 1$ the analogous problem in the bounded region $A_R := \{(\ze, \rho) \in \real^2 : 1 < \ze^2 + \rho^2 < R^2\}$ reads:

\hspace{1ex}

\noi\textbf{Problem $P_{I,\Ih}\, (A_R)$}. {\it Let $I \in C(S_1, \real_+)$ and $\Ih \in C(S_R, \real_+)$ be positive symmetric functions. Determine all symmetric solutions $\HH \in C^1 (A_R) \cap C(\overline{A_R})$ of {\rm (\ref{2.1})} satisfying the boundary conditions $|\HH|\big|_{S_1} = I$ and $|\HH|\big|_{S_R} = I_0\Ih$, where the positive constant $I_0$ is part of the solution. 
}

\hspace{1ex}

\noi Note that in the bounded case the second boundary condition is assumed only up to a positive factor $I_0$. This is comparable to the standard boundary value problem with prescribed normal component of the field at the boundary. In the case of constant normal component at both boundaries the monopole field is the unique solution, whose amplitude can be fixed only once for both boundaries. 

The zeroes of harmonic fields play a crucial role in our solution ansatz. The following facts are well known: 
Axisymmetric harmonic fields $\HH$ have only a finite number of isolated zeroes with finite negative indices (``x-points'') in $A_\infty$. Let these zeroes be contained in the annulus $A_R \subset A_\infty$, then its number $\nu (\HH, A_R)$ can be computed by 
\be \label{2.8}
\nu = \nu(\HH, A_R) = \ro - \roh\, ,
\ee
where $\ro$ and $\roh$ mean the rotation numbers of $\HH$ along the circles $S_1$ and $S_R$, respectively. The rotation number $\ro \in \ints$ counts the number of turns, the field vector makes when circling once around $S_1$. The number $\nu$ counts a zero as often as indicated by its index.  Equation (\ref{2.8}) is clearly an analogue of the argument principle in complex analysis and it has likewise some invariance properties with respect to continuous deformations (see section 10 in \cite{KR22}); in particular, $\roh$ is constant in the limit $R\ra \infty$, which yields for $A_\infty$ the relation  
\be \label{2.9}
\nu = \nu(\HH, A_\infty ) = \ro - \dt + 1\, ,
\ee
where $\dt$ denotes the {\em exact} decay order of $\HH$, i.e.\  
$|\HH (\xx)| = O(|\xx|^{-\dt})$ but not  $|\HH (\xx)| = O(|\xx|^{-(\dt + 1)})$ for $|\xx| \ra \infty$.\footnote{ {\em Exact} decay orders are often denoted by 
$\dt$ as opposed to $\de$ for decay orders in general.}

The basic idea to solve the intensity problem is the same as for the direction problem: expressing a two-dimensional vector field via the exponential mapping by argument and modulus, i.e.\ direction and intensity, and to derive from eq.\ (\ref{2.6}) corresponding equations for these quantities. 
The zeroes of the harmonic field, whose number by (\ref{2.8}) is related to the rotation numbers at the boundaries, must be treated separately. This leads to the ansatz 
\be \label{2.10}
f(z,\zb) = h(z)\, e^{g(z,\zb)} 
\ee
with the analytic function
\be \label{2.11}
h(z) := \prod_{n=1}^{\ro - \roh} (z - z_n)\, z^{-\ro} .
\ee
Inserting (\ref{2.10}) into (\ref{2.6}) one obtains
\be \label{2.12}
-\pa_\zb\, g = \frac{1}{2}\, \frac{1}{z-\zb} \Big( \frac{\hb}{h}\, e^{- 2i \Im g} - 1 \Big) .
\ee
Introducing the angle-type variable $\Psi$ by $\hb/h =: \exp (i \Psi)$, writing $2 g =: p + i q$, and using the real variables $(\ze, \rho)$ in $A_R$, eq.\ (\ref{2.12}) takes the real form
\begin{equation}\label{2.13}
\left.
\begin{array}{rl}
\disp \pa_\ze p - \pa_\rho q - \frac{1}{\rho}\sin (q - \Psi)& = 0 , \\[1.5ex]
\disp \pa_\rho p + \pa_\ze q - \frac{1}{\rho}\Big(\cos (q - \Psi) - 1\Big) & = 0 . 
\end{array} \right\}
\end{equation}
The symmetry condition (\ref{2.7}) is satisfied by the choice of a symmetric set $S$ of zeroes in (\ref{2.11}), i.e.\ $z \in S$ implies $\zb \in S$, and by the following symmetries of $p$ and $q$:
\be \label{2.14}
p (\ze , -\rho) = p(\ze, \rho) \, ,\qquad q(\ze, -\rho) = - q (\ze, \rho)\, .
\ee
The boundary conditions on $f \triangleq \HH$ translate into ones for $p$:
\be \label{2.15}
p\big|_{S_1} := 2\, \Big( \ln I - \sum_{n=1}^{\ro- \roh} \ln |z - z_n|\Big|_{S_1} \Big) , 
\ee
\be \label{2.16}
p\big|_{S_R} := 2\, \Big( \ln \Ih - \sum_{n=1}^{\ro- \roh} \ln |z - z_n|\Big|_{S_R} + p_0 \Big) , 
\ee
where $p_0$ combines the constants $\ln I_0$ and $\ro \ln R$.
The system (\ref{2.13}) allows the reduction to a semilinear elliptic equation in the direction-related variable $q$ alone, which was the basis of the solution of the direction problem in \cite{KR22}. An analogous reduction to the intensity-related variable $p$ seems not to be possible. Instead, we transform to zero boundary conditions in $p$, which induce ''natural`` boundary conditions for a new direction-related variable $u$: let $(p_l, q_l)$ be the solution of the linear auxiliary boundary value problem (see appendix B)
\be \label{2.17}
\left.
\begin{array}{rl}
\disp \pa_\ze p_l - \pa_\rho q_l & = 0\,  \\[1ex]
\disp \pa_\rho p_l + \pa_\ze q_l & = 0\,  
\end{array} \right\} \mbox{ in } A_R\, ,
\end{equation}
\be \label{2.18}
p_l\big|_{S_1} := p\big|_{S_1}\; ,\qquad p_l\big|_{S_R} := p\big|_{S_R}\, ,
\ee
where $p|_{S_1}$ and $p|_{S_R}$ are given by eqs.\ (\ref{2.15}) and (\ref{2.16}), respectively. In the variables $\wp := p - p_l$ and $u := q - q_l$ and with 
$\Om := \Psi - q_l$, the system (\ref{2.13}) now takes the form
\begin{equation}\label{2.19}
\left.
\begin{array}{rl}
\disp \pa_\ze \wp - \pa_\rho u - \frac{1}{\rho}\sin (u - \Om)& = 0 , \\[1.5ex]
\disp \pa_\rho \wp + \pa_\ze u - \frac{1}{\rho}\Big(\cos (u - \Om) - 1\Big) & = 0  
\end{array} \right\}
\end{equation}
with boundary condition 
\be \label{2.20}
\wp \big|_{\pa A_R} = 0\, .
\ee
Note that the data of the problem, viz., positions of the zeroes of the solution and intensity values at the boundary, are now represented by the ``data function'' $\Om$.   
Eliminating the variable $\wp
$ from (\ref{2.19}) yields
\be \label{2.21}
\Delta u + \pa_\ze \Big( \frac{1}{\rho} \big(1 - \cos (u - \Om) \big) \Big) + \pa_\rho \Big(\frac{1}{\rho}  \sin (u - \Om) \Big) = 0
\ee
or
\be \label{2.22}
\na \cdot \Big( \na u + \aaa[u - \Om]\, \frac{u}{\rho}\Big) = \na \cdot \Big( \aaa [u - \Om]\, \frac{\Om}{\rho} \Big)
\ee
with 
\be \label{2.23}
a_\ze [x] := \frac{1 - \cos x}{x}\, ,\qquad a_\rho [x] := \frac{\sin x}{x}\, ,
\ee
which is formally the same equation as in \cite{KR22}. To obtain boundary conditions for $u$ we evaluate condition (\ref{2.20}) in the form $\ttt \cdot \na \wp |_{\pa  A_R} = 0$ by (\ref{2.19}) to obtain
\begin{equation}\label{2.24}
\begin{array}{rl}
 0 = - \ttt\cdot \na \wp \big|_{\pa A_R} & = (n_\rho \pa_\ze - n_\ze \pa_\rho )\,  \wp\big|_{\pa A_R} \\[1.5ex]
 & \disp = \Big[n_\ze \Big(\pa_\ze u + \frac{1}{\rho} \big( 1 - \cos(u - \Om)\big) \Big) \\[1.5ex] 
 & \disp \qquad\qquad + n_\rho \Big(\pa_\rho u + \frac{1}{\rho} \sin (u- \Om)\Big) \Big]_{\pa A_R} \\[2ex]
 & \disp = \nn \cdot \Big[ \na u + \aaa [u - \Om]\, \frac{u - \Om}{\rho} \Big]_{\pa A_R}\, ,
\end{array}
\end{equation}
where $\nn = (n_\ze , n_\rho)$ and $\ttt = (-n_\rho , n_\ze)$ denote unit normal and tangential vectors, respectively, at $\pa A_R$. It is this Neumann-type nonlinear boundary condition in contrast to a zero boundary condition, which makes the difference between intensity and direction problem.

The symmetry requirements on $\wp$, $u$, and $\Om$ follow by (\ref{2.14})--(\ref{2.18}) (see also appendix B):
\be \label{2.25}
\wp(\ze, -\rho) = \wp(\ze, \rho)\, ,\quad u(\ze, -\rho) = - u(\ze, \rho)\, ,\quad \Om(\ze, -\rho) = - \Om (\ze, \rho)\, .
\ee

Equation (\ref{2.22}) suggests a weak formulation: let 
$$
C_{as}^1 (\ARb) := \big\{ \psi \in C^1(\ARb) : \psi (\ze, -\rho) = - \psi (\ze, \rho) \big\}
$$
and
$$
H_{as}^1 (A_R) :={\rm clos} \big(C_{as}^1 (\ARb)\, , \| \na \cdot \|_{L^2 (A_R)} \big) 
$$
be antisymmetric versions of the usual spaces of $C^1$-test-functions and $H^1$-functions. A function $u \in H^1_{as} (A_R)$ is then called weak solution of the boundary value problem (\ref{2.22})--(\ref{2.25}), if $u$ satisfies 
\be \label{2.26}
\int_{A_R} \na u \cdot \na \vii\, \rd \ze \rd \rho + \int_{A_R} \frac{u}{\rho}\, \aaa[u - \Om] \cdot \na \vii\, \rd \ze \rd \rho =
\int_{A_R} \frac{\Om}{\rho}\, \aaa[u - \Om] \cdot \na \vii\, \rd \ze \rd \rho
\ee
for all test functions $\psi \in C^1_{as} (\ARb)$. Obviously, by integrating by parts eq.\ (\ref{2.26}) is equivalent to eqs.\ (\ref{2.22}) and (\ref{2.24}) on functions $u \in C^2 (A_R) \cap C^1 (\ARb)$.

Antisymmetry allows us to reduce the problem to the half-annulus $A_R^+ := A_R \cap \{\rho >0\}$, bounded by the half-circles $S_1^+$, $S_R^+$, and (parts of) the symmetry-axis $S\!A$. Solution and test functions vanish by antisymmetry on the symmetry-axis, which makes the essential difference between the boundary components $S_{1,R}^+$ and $S\!A$ of $A_R^+$. Polar coordinates $(r, \vi)$ are here more appropriate to take into account this difference. They map $A_R^+$ on the rectangle $Q_R := (1,R)\times (0,\pi)$ with side walls $\{\vi = 0\}$ and $\{\vi =\pi\}$, at which solution and test functions vanish; this allows the application of Hardy-type inequalities, a crucial ingredient of our solution procedure. 

Equation (\ref{2.22}) now takes the form\footnote{The tilde denotes again dependence on polar coordinates: $\ut (r, \vi) := u(r \cos \vi, r \sin \vi)$, etc.}
\be \label{2.26a}
\begin{array}{c}
\disp \frac{1}{r} \pa_r (r \pa_r \ut) + \frac{1}{r^2} \pa_\vi^2 \ut + \frac{1}{r} \pa_r \Big( r\, a_r [\ut - \Omt]\, \frac{\ut}{r \sin \vi}\Big) + \frac{1}{r} \pa_\vi \Big(a_\vi [\ut - \Omt]\, \frac{\ut}{r \sin \vi}\Big) \\[1.5ex]
\disp \qquad \qquad = \frac{1}{r} \pa_r \Big( r\, a_r [\ut - \Omt]\, \frac{\Omt}{r \sin \vi}\Big) + \frac{1}{r} \pa_\vi \Big(a_\vi [\ut - \Omt]\, \frac{\Omt}{r \sin \vi}\Big) ,
\end{array}
\ee
where we made use of the abbreviations
\be \label{2.27}
\left.
\begin{array}{rl}
a_r & := a_\ze \cos \vi + a_\rho \sin \vi ,\\[1ex]
a_\vi & := - a_\ze \sin \vi + a_\rho \cos \vi .
\end{array}
\right\}
\ee
The weak formulation (\ref{2.26}) reads
\be \label{2.28}
\begin{array}{l}
 \disp \int_{Q_R}\Big(\pa_r \ut\, \pa_r \pst + \frac{1}{r^2}\, \pa_\vi \ut\, \pa_\vi \pst \Big) \rd \vi\, r \rd r \\[1ex]
 \disp \qquad + \int_{Q_R} \frac{\ut}{\sin \vi} \Big( a_r [\ut - \Omt]\, \pa_r \pst + a_\vi [\ut - \Omt]\, \frac{1}{r} \pa_\vi \pst \Big) \rd \vi \rd r \\[1ex]
 \disp \qquad\qquad = \int_{Q_R} \frac{\Omt}{\sin \vi} \Big( a_r [\ut - \Omt]\, \pa_r \pst + a_\vi [\ut - \Omt]\, \frac{1}{r} \pa_\vi \pst \Big) \rd \vi \rd r  
\end{array}
\ee
with $\pst \in \CC^1 (\QRb)$, $\ut \in \HC_0 (Q_R)$, and appropriate spaces
$$
\CC^1 (\QRb) := \big\{\pst \in C^1(\QRb) : \pst (\cdot\, , 0) = \pst(\cdot\, , \pi) = 0 \big\} ,
$$
$$\HC_0 (Q_R) := {\rm clos} \bigg(\CC^1 (\QRb)\, , \int_{Q_R} \Big( (\pa_r \,\cdot )^2 + \frac{1}{r^2} (\pa_\vi\, \cdot)^2 \Big) \rd \vi\, r \rd r\bigg) .
$$
Equations (\ref{2.26a}), (\ref{2.27}) with given data function $\Omt$ constitute {\bf problem} $P_\Omt (Q_R)$, for which the following theorem holds.
\begin{thm}
 Let $R > 1$ and $\Omt \in L^\infty (Q_R)$ with bound
 \be \label{2.29}
 |{\widetilde \Om} (\cdot\, , \vi)| \leq K\, \sin \vi \qquad \mbox{in }\, Q_R
 \ee
 for some constant $K > 0$, then problem $P_\Omt (Q_R)$ has a unique weak (in the sense of eq.\ (\ref{2.28})) solution $u \in \HC_0 (Q_R)$ satisfying the bound
 \be \label{2.30}
 \int_{Q_R} |\nt \ut|^2\,  \vih^{-\beta } \rd \vi\, r^{1-\eta} \rd r \leq C
 \ee
 with $\beta = 1/5$, $\eta = 10^{-3}$, and $C$ denoting some constant that depends on $K$, but does not depend on $R$. 
\end{thm}
In (\ref{2.30}), $\nt$ denotes the gradient in polar coordinates  and $\vih$ the ``tent function'' (see appendix A).

In the unbounded case a function $u \in \HC_0^{loc} (Q_\infty)$ with $Q_\infty := (1, \infty)\ti (0,\pi)$ is called weak solution of eqs.\ (\ref{2.26a}), (\ref{2.27}), if $\ut$ satisfies eq.\ (\ref{2.28}) (with $Q_R$ replaced by $Q_\infty$) for all functions $\psi \in \CC^1 (\Qib)$, where 
$$
\ba{l}
\CC^1 (\Qib) := \big\{\pst \in C^1(\Qib) : \pst (\cdot\, , 0) = \pst(\cdot\, , \pi) = 0 \mbox{ and } \\[1ex]
\disp \qquad \qquad \qquad \pst(r, \cdot) = 0 \mbox{ for all } r > R_0 \mbox{ and some } R_0 > 1 \big\} .
\ea
$$
For the corresponding {\bf problem} $P_\Omt (Q_\infty)$ holds the theorem:
\begin{thm}
 Let $\Omt \in L^\infty (Q_\infty)$ with bound
 \be \label{2.31}
 |{\widetilde \Om} (\cdot\, , \vi)| \leq K\, \sin \vi \qquad \mbox{in }\, Q_\infty
 \ee
 for some constant $K > 0$, then problem $P_\Omt (Q_\infty)$ has a unique weak solution $\ut \in \HC_0^{loc} (Q_\infty)$
 satisfying the bound
 \be \label{2.32}
 \int_{Q_\infty} |\nt \ut|^2\,  \vih^{-\beta } \rd \vi\, r^{1-\eta} \rd r \leq C
 \ee
 with $\beta = 1/5$, $\eta = 10^{-3}$, and some constant $C > 0$.
\end{thm}
Based on these results the following both theorems answer the original problems $P_{I, \Ih} (A_R)$ and $P_I (A_\infty)$. H\"older continuity of the data $I$ and $\Ih$ is here sufficient to ensure continuity of the solution up to the boundary; only at the symmetry axis we need some more regularity, viz.,
\be \label{2.33}
\left.
\begin{array}{cl}
|I(\vi) - I(0)|  = O \big(|\vi|^{1 + \at}\big) & \quad\mbox{ for } |\vi|\ra 0\, ,\\[1.5ex]
|I(\vi) - I(\pi)|  = O \big(|\vi - \pi|^{1 + \at}\big) & \quad\mbox{ for } |\vi - \pi|\ra 0
\end{array}
\right\}
\ee
with some $0< \at < 1$, to ensure conditions (\ref{2.29}), (\ref{2.31}). 
\begin{thm}
 Let $I$ and $\Ih$ be H\"older continuous, symmetric intensity functions 
 satisfying conditions (\ref{2.33}), let $\roh \in \nat$, and let $\{z_1, \ldots,z_{2 N} \}$, $N \in \nat_0$, be a symmetric set of points (avoiding the symmetry axis) in $A_R$.\footnote{In the case $N=0$ the set is empty.}
 Then, problem $P_{I, \Ih} (A_R)$ has a (up to a sign) unique solution $\HH$ with rotation number $\roh$ at $S_R$ and vanishing at  $z_1, \ldots,z_{2N}$ and nowhere else.
\end{thm}

\begin{thm}
 Let $I$ be a H\"older continuous, symmetric intensity function 
 satisfying conditions (\ref{2.33}), let $\dt \in \nat\setminus \{1\}$, and let $\{z_1, \ldots,z_{2 N} \}$, $N \in \nat_0$, be a symmetric set of points (avoiding the symmetry axis) in $A_\infty$. Then, problem $P_{I} (A_\infty)$ has a (up to a sign) unique solution $\HH$ with exact decay order  $\dt$ and vanishing at  $z_1, \ldots,z_{2N}$ and nowhere else.
\end{thm}
We close this section with two comments:

\noi 1. The restriction to an even number of zeroes in the theorems 2.3 and 2.4 is clearly not optimal, since axisymmetric harmonic fields with an odd number of zeroes (and necessarily at least one of them on the symmetry axis) are well known. The same restriction appeared provisionally in the direction problem \cite{KR22}, where, however, it could be removed in the unbounded case by taking suitable linear combinations of ``even solutions'' to generate an odd solution; a method, which obviously does not work in the intensity problem. So, for the present it is a conjecture that theorem 2.4 holds for an odd number of zeroes as well.\footnote{It is this deficit why we qualified in section 1 our solution of the axisymmetric intensity problem as only   ``almost'' complete.} 

\noi 2. Without further information theorem 2.4 predicts for a given intensity function $I$ an infinite set of solutions characterized by the exact decay order and positions of an (even) number of zeroes. Applying the theorem to the gravitational field of the earth we have additionally condition (\ref{1.2}), which clearly implies a rotation number $\ro = 1$ and by (\ref{2.9}) $\dt = 2$ and no zeroes, hence uniqueness in agreement with \cite{B68}. Concerning the geomagnetic field the present situation is characterized by a dominating dipole-field, which means $\ro =2 $, $\dt = 3$ and hence again no zeroes and uniqueness. In earlier epochs, however, especially during pole reversals, $\ro > 2$ is possible, while a (however small) dipole part always implies $\dt = 3$, which implies by (\ref{2.9}) zeroes of the magnetic field and hence nonuniqueness. In that case the reconstruction of the magnetic field from intensity data alone is in principle not possible.    
\sect{Uniqueness in the problems $P_\Om (Q_R)$ and $P_\Om (Q_\infty)$}
The function $\Om$ is a functional of the intensity function $I$ (and $\Ih$ in the bounded problem) and the set of zeroes parametrizing the nonuniqueness of the original problems $P_I (A_\infty)$ and $P_{I, \Ih} (A_R)$. Starting with the bounded case, $\Om : Q_R \ra \real$ is in this section a fixed function with property (\ref{2.29}) and $u_1$, $u_2$ are two weak solutions of eq.\ (\ref{2.21}) in polar coordinates, i.e.\ $u_1$ satisfies\footnote{We use in this and the subsequent sections exclusively polar coordinates; the tilde indicating so far these coordinates can thus safely be omitted.} 
$$
\ba{l}
 \disp \int_{Q_R} \Big( \pa_r u_1\, \pa_r \psi + \frac{1}{r^2}\, \pa_\vi u_1 \,\pa_\vi \psi \Big) \rd \vi\, r \rd r \\[1.5ex]
 \disp \qquad \qquad + \int_{Q_R} \frac{1}{r \sin \vi}\Big(\big( \cos \vi - \cos (u_1 - \Om + \vi) \big) \pa_r \psi \\[1ex]
 \disp \qquad \qquad \qquad \qquad \qquad + \big( \sin ( u_1 - \Om + \vi) - \sin \vi \big) \frac{1}{r} \pa_\vi \psi \Big) \rd \vi \, r \rd r = 0
\ea
$$
with test function $\psi \in \CC^1 (\QRb)$, and analogously for $u_2$.
The ``perturbation'' $\du := u_1 - u_2$ then satisfies the equation\be \label{3.1}
\int_{Q_r} \na \du \cdot \na \psi \, \rd \mu + \int_{Q_R} \aaa' \cdot \na \psi \,\frac{\du}{r \sin \vi}\, \rd \mu  = 0
\ee
with $\na  = \eee_r \pa_r + \eee_\vi \frac{1}{r} \pa_\vi$, $\rd \mu = \rd \vi\, r \rd r$, $\aaa' = a'_r \eee_r + a'_\vi \eee_\vi$, and 
\be \label{3.2}
\left. 
\ba{l}
\disp a'_r := \sin \Big( \frac{1}{2}\, (u_1 + u_2 ) - \Om + \vi\Big)\, \frac{\sin \big( (u_1 - u_2)/2\big)}{(u_1 - u_2)/2} \, ,\\[2ex]
\disp a'_\vi := \cos \Big( \frac{1}{2}\, (u_1 + u_2 ) - \Om + \vi\Big)\, \frac{\sin \big( (u_1 - u_2)/2\big)}{(u_1 - u_2)/2} \, .
\ea
\right\}
\ee
In order to prove $\du = 0$ in (\ref{3.1}) we make use of the same sequence of test functions as in \cite{KR22}:
$$
(\psi_n) := \Big(\frac{\du^+}{\du^+\! + 1 /n}\Big) ,
$$
where $\du^+$ denotes the positive part of $\du$. With $\du \in \HC_0 (Q_R)$ we also have $\du^+$ and $\psi_n \in \HC_0 (Q_R)$ with
\be \label{3.3}
\na \psi_n = \left\{ 
\ba{cc}
\disp \frac{1}{n}\, \frac{\na \du^+}{(\du^+\! + 1/n)^2} & \quad \mbox{on supp}\,\du^+ , \\[3ex]
\disp 0 & \quad \mbox{ else}\, .
\ea\right.
\ee
Note that any $\na \psi_n$ can be approximated in $L^2 (Q_R)$ by elements from $\CC^1 (\QRb)$, thus $\psi_n$ itself is an admissible (test) function. Inserting (\ref{3.3}) into (\ref{3.1}) one obtains suitably splitted:
\be \label{3.4}
\ba{l}
\disp 0= \frac{1}{n} \int_{Q_R} \frac{|\na \du^+|^2}{(\du^+\! + 1/n)^2} \, \rd \mu + \frac{1}{n} \int_{Q_R} \frac{a'_r}{r \sin \vi} \, \frac{\du^+ \pa_r \du^+}{(\du^+\! + 1/n)^2}\, \rd \mu \\[3ex]
\disp \qquad \qquad + \frac{1}{n} \int_1^R \int_0^{\pi/2} \Big(\frac{a'_\vi - 1}{r \sin \vi} + \frac{1}{r \sin \vi}\Big) \, \frac{\du^+ \frac{1}{r} \pa_\vi \du^+}{(\du^+\! + 1/n)^2}\, \rd \vi\, r \rd r \\[3ex]
\disp \qquad \qquad + \frac{1}{n} \int_1^R \int_{\pi/2}^{\pi} \Big(\frac{a'_\vi + 1}{r \sin \vi} - \frac{1}{r \sin \vi}\Big) \, \frac{\du^+ \frac{1}{r} \pa_\vi \du^+}{(\du^+\! + 1/n)^2}\, \rd \vi\, r \rd r\, .
\ea
\ee
Estimating the $\aaa'$-related terms by Cauchy-Schwarz, eq.\ (\ref{3.4}) can be arranged as follows:
\be \label{3.5}
\ba{l}
\ba{ll}
\disp \int_{Q_R} \frac{|\na \du^+|^2}{(\du^+\! + 1/n)^2} \, \rd \mu + \int_1^R \frac{1}{r^2}\bigg[ & 
\disp \int_0^{\pi/2} \frac{1}{\sin\vi} \, \frac{\du^+ \pa_\vi \du^+}{(\du^+\! + 1/n)^2}\, \rd \vi \\[2ex]
& \disp - \int_{\pi/2}^\pi \frac{1}{\sin\vi} \, \frac{\du^+ \pa_\vi \du^+}{(\du^+\! + 1/n)^2}\, \rd \vi \bigg] r \rd r 
\ea \\[7ex]
\disp \qquad \leq \Bigg[\bigg(\int_{Q_R} \Big(\frac{a_r'}{r \sin \vi}\Big)^2 \rd \mu \bigg)^{\frac{1}{2}} + 
\bigg( \int_1^R \int_0^{\pi/2} \Big(\frac{a_\vi' - 1}{r \sin \vi} \Big)^2 \rd \vi\, r \rd r \bigg)^{\frac{1}{2}} \\[3ex]
 \disp \qquad \qquad + \bigg(\int_1^R \int_{\pi/2}^\pi \Big(\frac{a_\vi' + 1}{r \sin \vi} \Big)^2 \rd \vi\, r \rd r \bigg)^{\frac{1}{2}}\Bigg] \ti
\bigg( \int_{Q_R} \frac{|\na \du^+|^2}{(\du^+\! + 1/n)^2}\,  \rd \mu \bigg)^{\frac{1}{2}} . 
\ea
\ee
The three terms in the bracket on the right-hand side of (\ref{3.5}) can further be estimated by (\ref{3.2}), (\ref{2.29}), and (\ref{A.2}): 
$$
\ba{l}
\disp \int_{Q_R} \Big(\frac{a_r'}{r \sin\vi}\Big)^2 \rd \mu \leq \int_{Q_R} \bigg(\frac{ \sin \big( (u_1 + u_2)/2 - \Om + \vi\big)}{r \sin \vi} \bigg)^2 \rd \mu \\[3ex]
\disp \qquad \leq \int_{Q_R} \bigg(1 + \frac{ \big|\sin \big( (u_1 + u_2)/2 - \Om \big)\big|}{r \sin \vi} \bigg)^2 \rd \mu \\[3ex]
\disp \qquad \leq 2 \int_{Q_R} \rd \mu + 
2 \int_{Q_R} \bigg(\Big(\frac{u_1}{r \sin \vi}\Big)^2 + \Big(\frac{u_2}{r \sin \vi}\Big)^2 \bigg) \rd \mu + 4 \int_{Q_R} \Big(\frac{\Om}{r \sin\vi}\Big)^2 \rd \mu \\[3ex]
\disp \qquad \leq \pi (1 + 2 K^2) R^2 + 8 \max_{[0, \pi]}
\Big( \frac{\vih}{\sin \vi}\Big)^2 \int_{Q_R} \Big[ \Big(\frac{1}{r} \pa_\vi u_1\Big)^2 + \Big(\frac{1}{r} \pa_\vi u_2\Big)^2 \Big] \rd \mu \\[3ex]
\disp \qquad  \qquad \leq \pi (1 + 2 K^2) R^2 + 2 \pi^2 \int_{Q_R} \big(| \na u_1|^2 + |\na u_2|^2 \big)\, \rd \mu =: C_{\aaa'}\, .
\ea
$$
Elementary estimates yield
$$
\ba{l}
\disp \bigg( \cos \Big(\frac{1}{2}\, (u_1 + u_2) - \Om + \vi \Big) \frac{\sin\big((u_1 - u_2)/2\big)}{(u_1 - u_2)/2} -1 \bigg)^2 \\[3ex]
\disp \quad \leq \bigg(\Big|1 - \cos \Big(\frac{1}{2}\, (u_1 + u_2) - \Om + \vi \Big) \Big| + \Big|1 - \frac{\sin\big((u_1 - u_2)/2\big)}{(u_1 - u_2)/2}\Big|\bigg)^2 \\[3ex]
\disp \quad \quad \leq \Big( \Big|\frac{1}{2}\, (u_1 + u_2) - \Om + \vi \Big| + \Big|\frac{1}{2} (u_1 - u_2) \Big|\Big)^2
\leq 3 (u_1^2+ u_2^2) + 8 ( \Om^2 + \vi^2)\, .
\ea
$$
Thus the second term can be estimated as follows:
$$
\ba{l}
\disp \int_1 ^R \int_0^{\pi/2} \Big(\frac{a_\vi' - 1}{r \sin\vi}\Big)^2 \rd \vi\, r \rd r \\[3ex]
\disp \leq 
 \int_1^R \int_0^{\pi/2} \bigg[ 3 \bigg(\Big(\frac{u_1}{r \sin \vi}\Big)^2 + \Big(\frac{u_2}{r \sin \vi}\Big)^2 \bigg) + 8 \bigg(\Big(\frac{\Om}{ \sin\vi}\Big)^2 + \Big(\frac{\vi}{\sin \vi }\Big)^2\bigg) \bigg]\rd \vi\, r \rd r \\[3ex]
 \disp \qquad  \qquad \leq 3 \pi^2 \int_{Q_R} \big(| \na u_1|^2 + |\na u_2|^2 \big)\, \rd \mu + 2 \pi \Big( K^2 + \frac{\pi^2}{4} \Big) R^2 =: \Ct_{\aaa'}\, .
\ea
$$
Similarly, by
$$
\ba{l}
\disp \bigg( \cos \Big(\frac{1}{2}\, (u_1 + u_2) - \Om + \vi \Big) \frac{\sin\big((u_1 - u_2)/2\big)}{(u_1 - u_2)/2} + 1 \bigg)^2 \\[3ex]
\disp \qquad = \bigg( \cos \Big(\frac{1}{2}\, (u_1 + u_2) - \Om + \vi - \pi\Big) \frac{\sin\big((u_1 - u_2)/2\big)}{(u_1 - u_2)/2} - 1 \bigg)^2 
\\[3ex]
\disp \qquad \qquad \leq 3 (u_1^2+ u_2^2) + 8 \big( \Om^2 + (\pi - \vi)^2\big)\, ,
\ea
$$
one obtains the same estimate as before for the third term:
$$
\int_1 ^R \int_{\pi/2}^\pi \Big(\frac{a_\vi' + 1}{r \sin\vi}\Big)^2 \rd \vi\, r \rd r \leq \Ct_{\aaa'}\, .
$$
The integrand in the second term on the left-hand side of inequality (\ref{3.5}) turns out to be a $\vi$-derivative:
\be \label{3.5a}
\frac{\du^+ \pa_\vi \du^+}{(\du^+\! + 1/n)^2} = \pa_\vi f(n\, \du^+)\, ,
\ee
where $f$ denotes the nonnegative function
\be \label{3.5b}
f: [0,\infty) \ra \real_+ \, ,\quad x\mapsto \ln (1+x) - \frac{x}{1 + x}\, .
\ee
Thus by integrating by parts one obtains
$$
\ba{l}
\disp \int_0^{\pi/2} \frac{1}{\sin \vi} \, \pa_\vi f(n\, \du^+)\, \rd \vi -
\int_{\pi/2}^{\pi} \frac{1}{\sin \vi} \, \pa_\vi f(n\, \du^+)\, \rd \vi \\[2ex]
\disp \qquad = \int_0^{\pi/2} \frac{\cos \vi}{(\sin \vi)^2} \, f(n\, \du^+)\, \rd \vi -
\int_{\pi/2}^{\pi} \frac{\cos \vi}{(\sin \vi)^2} \, f(n\, \du^+)\, \rd \vi \\[3ex]
\disp \qquad \qquad \qquad + 2 \, f(n\, \du^+)\big|_{\vi =\pi/2} \geq 0\, ,
\ea
$$
which implies the second term on the left-hand side of (\ref{3.5}) to be nonnegative.

Applying all these estimates, inequality (\ref{3.5}) boils down to
$$ \int_{Q_R} \frac{|\na \du^+|^2}{(\du^+\! + 1/n)^2} \, \rd \mu \leq  \big(C_{\aaa'}^{1/2} + 2\, \Ct_{\aaa'}^{1/2}\big)
\bigg( \int_{Q_R} \frac{|\na \du^+|^2}{(\du^+\! + 1/n)^2}\,  \rd \mu \bigg)^{\frac{1}{2}} 
$$
or 
\be \label{3.6}
\ba{l}
\disp 9\, \Ct_{\aaa'} \geq \int_{Q_R} \frac{|\na \du^+|^2}{(\du^+\! + 1/n)^2} \, \rd \mu = \int_{Q_R} \big| \na \ln ( 1 + n\, \du^+)\big|^2 \rd \mu \\[3ex]
\disp \qquad\quad \geq \frac{1}{\pi^2 R^2} \int_{Q_R} \big( \ln (1 + n\, \du^+)\big)^2 \rd\mu\, ,
\ea
\ee
where for the last estimate we made use of inequality (\ref{A.3}). By Levi's theorem we conclude convergence in the sequence $\big( \ln (1 + n\, \du^+) \big)_{n \in \nat}$, which in turn implies $\du^+ = 0$ a.e. The argument applies verbatim to $- \du^- = - \min\{ \du, 0\}$, hence $\du = 0$ a.e. We summarize this result in the following proposition.
\begin{prop}
 Let $u_1$, $u_2 \in \HC_0 (Q_R)$ be weak solutions of problem $P_\Om (Q_R)$ with $1 < R < \infty$ and $\Om \in L^\infty (Q_R)$ satisfying the bound (\ref{2.29}). Then $u_1$ and $u_2$ coincide a.e.\ in $Q_R$.
\end{prop}

In the unbounded case we deal with solutions $u \in \HC^{loc}_\beta (Q_\infty)$ of problem $P_\Om (Q_\infty)$ that satisfy a bound of the form 
$$
\|\na u \|^2_{(\beta, \eta)}  = \int_{Q_\infty} |\na u|^2 \, \vih^{-\beta} \rd\vi\, r^{1-\eta} \rd r < C\, , \quad 0< \beta < \frac{1}{2}\, ,\quad 0 < \eta < 1\, . 
$$
In $Q_\infty$ the weights, especially $r^{-\eta}$, cannot be neglected. In order to incorporate the weights into the governing equation (\ref{3.1})  we replace $u$ by the variable $v:= \vih^\beta r^\eta  u$ to obtain
\be \label{3.7}
\ba{l}
\disp \int_{Q_\infty} \na \dv \cdot \na \psi\, \rd \mu_\beta^\eta + \int_{Q_\infty}\Big(\frac{a_r'}{r\sin \vi} - \frac{\eta}{r} \Big) \dv\, \pa_r \psi\, \rd \mu_\beta ^\eta \\[3ex]
\disp  \qquad + \int_1^\infty \frac{1}{r^2}\bigg( \int_0^{\pi/2} \Big(\frac{a_\vi'}{\sin \vi} - \frac{\beta}{\vih}\Big)\dv \, \pa_\vi \psi\, \vih^{-\beta} \rd \vi \\[2ex]
\disp \qquad \qquad \qquad + \int_{\pi/2}^\pi \Big(\frac{a_\vi'}{\sin \vi} + \frac{\beta}{\vih}\Big) \dv \, \pa_\vi \psi\, \vih^{-\beta} \rd \vi \bigg) r^{1-\eta} \rd r = 0
\ea
\ee
with $\rd \mu_\beta^\eta = \vih^{-\beta} \rd \vi \, r^{1-\eta} \rd r$ and $\aaa'$ given as before by (\ref{3.2}). Using now as test functions the sequence 
$$
(\psi_n) := \Big(\frac{\dv^+}{\dv^+\! + 1 /n}\Big) ,
$$
the arguments follow closely those in the bounded case; we skip thus the details. Inserting $\psi_n$ into (\ref{3.7}), estimating and rearranging the various terms as in the bounded case, one obtains instead of (\ref{3.5}):
\be \label{3.8}
\ba{l}
\ba{ll}
\disp \int_{Q_\infty} \frac{|\na \dv^+|^2}{(\dv^+\! + 1/n)^2} \, \rd \mu_\beta^\eta + & \disp \int_1^\infty \frac{1}{r^2}\bigg(  
\int_0^{\pi/2} \Big(\frac{1}{\sin\vi} -\frac{\beta }{\vih}\Big)\, \frac{\dv^+ \pa_\vi \dv^+}{(\dv^+\! + 1/n)^2}\, \vih^{-\beta} \rd \vi \\[3ex]
& \disp \quad - \int_{\pi/2}^\pi \Big(\frac{1}{\sin\vi} - \frac{\beta}{\vih}\Big) \, \frac{\dv^+ \pa_\vi \dv^+}{(\dv^+\! + 1/n)^2}\, \vih^{-\beta} \rd \vi \bigg)
r^{1-\eta} \rd r 
\ea \\[7ex]
\disp \leq \Bigg[\bigg(\int_{Q_\infty} \Big(\frac{a_r'}{r \sin \vi} - \frac{\eta}{r}\Big)^2 \rd \mu_\beta^\eta \bigg)^{\frac{1}{2}} + 
\bigg( \int_1^\infty \int_0^{\pi/2} \Big(\frac{a_\vi' - 1}{r \sin \vi} \Big)^2 \vih^{-\beta} \rd \vi\, r^{1-\eta} \rd r \bigg)^{\frac{1}{2}} \\[3ex]
 \disp \quad + \bigg(\int_1^\infty \int_{\pi/2}^\pi \Big(\frac{a_\vi' + 1}{r \sin \vi} \Big)^2 \vih^{-\beta} \rd \vi\, r^{1-\eta} \rd r \bigg)^{\frac{1}{2}}\Bigg] \ti
\bigg( \int_{Q_\infty} \frac{|\na \dv^+|^2}{(\dv^+\! + 1/n)^2}\,  \rd \mu_\beta^\eta \bigg)^{\frac{1}{2}} . 
\ea
\ee
The three terms in the bracket on the right-hand side of (\ref{3.8}) are estimated by (\ref{2.31}), (\ref{2.32}), and (\ref{A.2}) along the lines of the bounded case. We just present the estimate of the first term: 
$$
\ba{l}
\disp \int_{Q_\infty} \Big(\frac{a_r'}{r \sin\vi} - \frac{\eta}{r} \Big)^2 \rd \mu_\beta^\eta  \leq \int_{Q_\infty} \bigg(\frac{1 + \eta}{r} + \frac{ \big|\sin \big( (u_1 + u_2)/2 - \Om \big)\big|}{r \sin \vi} \bigg)^2 \rd \mu_\beta^\eta \\[3ex]
\disp \quad \leq 2 (1 + \eta)^2 \int_1^\infty \frac{1}{r^2}\, r^{1-\eta} \rd r \int_0^\pi \vih^{-\beta} \rd \vi + \frac{\pi^2}{2} \int_{Q_\infty} \frac{1}{r^2}\bigg( \Big(\frac{u_1}{\vih} \Big)^2 + \Big(\frac{u_2}{\vih}\Big)^2 \bigg) \rd \mu_\beta^\eta \\[3ex]
\disp \qquad \qquad \quad + 4 \int_{Q_\infty} \Big(\frac{\Om}{r \sin\vi}\Big)^2 \rd \mu_\beta^\eta \\[3ex]
\disp \qquad \leq 4 \big((1 + \eta)^2 + 2 K^2\big) \frac{(\pi/2)^{1-\beta}}{\eta (1 - \beta)} + 2 \Big(\frac{\pi}{1 + \beta}\Big)^2 \Big(\| \na u_1\|^2_{(\beta, \eta)}  + \|\na u_2\|^2_{(\beta, \eta)} \Big) \\[2ex]
\disp \qquad \qquad =: C_{\aaa'}^{(\beta, \eta)} \, .
\ea
$$
The second term on the left-hand side of (\ref{3.8}) is again nonnegative. In fact, by (\ref{3.5a}), (\ref{3.5b}), and integration by parts one obtains
$$
\ba{l}
\disp \int_0^{\pi/2} \Big(\frac{1}{\sin \vi} - \frac{\beta}{\vih}\Big)\, \pa_\vi f(n\, \dv^+)\, \vih^{-\beta} \rd \vi -
\int_{\pi/2}^{\pi} \Big(\frac{1}{\sin \vi} - \frac{\beta}{\vih}\Big) \pa_\vi f(n\, \dv^+)\, \vih^{-\beta} \rd \vi \\[3ex]
\disp \qquad = - \int_0^{\pi/2} g' (\vih) \,  f(n\, \dv^+)\, \rd \vi -
\int_{\pi/2}^{\pi} g' (\vih) \,  f(n\, \dv^+)\, \rd \vi \\[3ex]
\disp \qquad \qquad \qquad + 2 g(\pi/2)\, f(n \dv^+)\big|_{\vi =\pi/2} \geq 0
\ea
$$
with
$$
g(y) := \Big( \frac{1}{\sin y} - \frac{\beta}{y}\Big) \, y^{-\beta} 
$$
and 
$$
- g'(y) = \Big( \frac{\cos y}{\sin^2 y} + \frac{\beta}{y \sin y} - \frac{\beta (1 + \beta)}{y^2} \Big) \, y^{-\beta} ,
$$
which is nonnegative for $0 < y \leq \pi/2$ and $0 < \beta < 1/2$ .

The rest of the argument is as in the bounded case with the result $\du = 0$ a.e.\ in $Q_R$ for any $R >1$.
\begin{prop}
 Let $u_1$, $u_2 \in \HC_\beta^{loc} (Q_\infty)$ be weak solutions of problem $P_\Om (Q_\infty)$ with bound (\ref{2.32}) and $\Om \in L^\infty (Q_\infty)$ satisfying the bound (\ref{2.31}). Then $u_1$ and $u_2$ coincide a.e.\ in $Q_\infty$.
\end{prop}

Finally, we prove uniqueness in a linearized, homogeneous version of eq.\ (\ref{2.22}), relevant for the solution of the linearized problem in section 4. Starting point is the weak formulation (\ref{2.28}) in polar coordinates with vanishing right-hand side and the variable $u$ in the vector field $\aaa$ replaced by a fixed function $w$. Using the Cartesian components $a_\ze$ and $a_\rho$ of the vector field $\aaa$, one obtains (cf.\ eq.\ (\ref{4.2}) below):
\be \label{3.9}
\ba{l}
\disp \int_{Q_R} \na u \cdot \na \psi\, \rd \mu
+ \int_{Q_R} \frac{a_\ze [w - \Om]}{r \tan \vi}\, u  \, \pa_r \psi \, \rd \mu 
+ \int_{Q_R} \frac{u}{r \tan \vi} \, \frac{1}{r} \pa_\vi \psi \, \rd \mu 
\\[3ex]
\disp \quad + \int_{Q_R} \frac{a_\rho [w - \Om] - 1}{r \tan \vi}\,  u \, \frac{1}{r} \pa_\vi \psi \, \rd \mu  
+ \int_{Q_R} \Big(\frac{a_\rho}{r} u\, \pa_r \psi  - \frac{a_\zeta}{r} u\, \frac{1}{r} \pa_\vi \psi \Big) \rd \mu = 0\, .
\ea
\ee
The sequence of test functions is now
$$
(\psi_n) := \Big(\frac{u^+}{u^+\! + 1 /n}\Big) ,
$$
and estimating and rearranging as in the bounded case leads now to:
\be \label{3.10}
\ba{l}
\disp \int_{Q_R} \frac{|\na u^+|^2}{(u^+\! + 1/n)^2} \, \rd \mu +  \int_1^R \frac{1}{r^2} \int_0^{\pi} \frac{1}{\tan\vi} \, \frac{u^+ \pa_\vi u^+}{(u^+\! + 1/n)^2}\, \rd \vi\, r \rd r\\[3ex]
\disp \qquad \leq \Bigg[\bigg(\int_{Q_R} \Big(\frac{a_\ze}{r \tan \vi} \Big)^2 \rd \mu \bigg)^{\frac{1}{2}} + 
\bigg(\int_{Q_R} \Big(\frac{a_\rho - 1}{r \tan \vi} \Big)^2 \rd \mu \bigg)^{\frac{1}{2}} \\[3ex]
\disp \qquad \qquad + 2\, \bigg(\int_{Q_R} \frac{1}{r^2} \, \rd \mu \bigg)^{\frac{1}{2}}  
\Bigg] \ti
\bigg( \int_{Q_R} \frac{|\na u^+|^2}{(u^+\! + 1/n)^2}\,  \rd \mu \bigg)^{\frac{1}{2}} . 
\ea
\ee
The estimate (\ref{3.10}) resembles (\ref{3.5}), but differs in the vector field $\aaa$ instead of $\aaa'$, in the trigonometric function $\tan \vi$ instead of $\sin \vi$ and in an additional (unproblematic) term on the right-hand side. Nevertheless, the second term on the left-hand side in (\ref{3.10}) is again nonnegative since by (\ref{3.5a}) and integration by parts we have
$$
\ba{ll}
\disp \int_0^{\pi} \frac{1}{\tan\vi} \, \frac{u^+ \pa_\vi u^+}{(u^+\! + 1/n)^2}\, \rd \vi & \disp = \int_0^{\pi} \frac{1}{\tan\vi} \, \pa_\vi f(n\, u^+)\, \rd \vi \\[3ex]
& \disp = \int_0^{\pi} \frac{1}{(\sin \vi)^2} \, f(n\, u^+)\, \rd \vi \geq 0\, , 
\ea
$$ 
and the three terms in the bracket on the right-hand side are all bounded. We consider again only the first term:
$$
\ba{l}
\disp \int_{Q_R} \Big(\frac{a_\ze}{r \tan\vi} \Big)^2 \rd \mu  \leq \int_{Q_R} \Big(\frac{1}{r \tan \vi}\,  \frac{ 1 - \cos (w - \Om)}{w - \Om} \Big)^2 \rd \mu \leq  \int_{Q_R} \Big(\frac{w -\Om}{r \tan\vi}\Big)^2 \rd \mu \\[3ex]
\disp \qquad \leq 2 \int_{Q_R} \Big(\frac{w}{r \sin \vi}\Big)^2 \rd \mu + 2 \int_{Q_R} \Big(\frac{\Om}{\sin \vi}\Big)^2 \rd \mu \leq  2 \pi^2 \int_{Q_R} | \na w |^2\, \rd \mu + \pi\, K^2 R^2 . 
\ea
$$
We conclude as above that $u= 0$ a.e.\ in $Q_R$, which yields the following proposition:
\begin{prop}
 Let $u \in \HC_0 (Q_R)$, $1 < R < \infty $ be a solution of eq. (\ref{3.9}) with given $w \in \HC_0 (Q_R)$ and $\Om \in L^\infty (Q_R)$ satisfying the bound (\ref{2.29}). Then $u = 0$ a.e.\ in $Q_R$.
\end{prop}
\sect{The linearized problem}
The linearized problem consists in the solution of eq.\ (\ref{2.28}) neglecting the $u$-dependence of the vector field $\aaa$. Thus, by (\ref{2.23}), $a_\rho$ and $a_\ze$ denote here arbitrary measurable functions with bounds 
\be \label{4.1}
-0.22< a_\rho \leq 1\quad  \mbox{ and }\quad 
|a_\ze | < 0.73\quad  \mbox{ on } \real\, .
\ee
Inserting the abbreviations (\ref{2.27}), eq.\ (\ref{2.28}) may be arranged as follows:
\be \label{4.2}
\begin{array}{l}
 \disp \int_{Q_R} \na u \cdot \na \psi\, \rd \mu + \int_{Q_R} \frac{u}{r \tan \vi} \Big( a_\ze \, \pa_r \psi + a_\rho  \frac{1}{r} \pa_\vi \psi \Big) \rd \mu \\[2ex]
 \disp \qquad \qquad + \int_{Q_R} u \Big( a_\rho \, \pa_r \psi - a_\ze \frac{1}{r} \pa_\vi \psi \Big) \frac{1}{r}\, \rd \mu \\[2ex] 
 \disp \qquad\qquad \qquad = \int_{Q_R} \frac{\Om}{r \sin \vi} \Big( a_r \, \pa_r \psi + a_\vi  \frac{1}{r} \pa_\vi \psi \Big) \rd \mu\, ,  
\end{array}
\ee
where $\na = \eee_r \pa_r + \eee_\vi \dfrac{1}{r} \pa_\vi$ and $\rd \mu = \rd \vi \, r \rd r$. The method of choice to prove the existence of a solution of eq.\ (\ref{4.2}) under quite general conditions is the application of the Lax-Milgram criterion, which requires, however, some coercivity property of the left-hand side in (\ref{4.2}). By Hardy's inequality the second term on the left-hand side in (\ref{4.2}) turns out as a (pseudo-) second-order term, which has to be controlled by the first term. A lot of the work in
\cite{KR22} was devoted to establish this control by the combined effect of three measures:   
a change of variables leads to an effective quantity $a_\rho$ with bound better than $1$, the introduction of a $\rho$-dependent weight leads to an improved constant in Hardy's inequality, and a weighted gradient allows better to exploit the asymmetry in the bounds of $a_\rho$ and $a_\ze$. All these measures -- adapted to polar coordinates -- will again be necessary to gain control over the second term on the left-hand side in (\ref{4.2}). Moreover, the third (lower order) term, not present in \cite{KR22}, requires an additional measure, viz., the consideration of an extended bilinear form, which allows the reformulation of (\ref{4.2}) as a Fredholm-type problem, to which Fredholm's alternative can be applied.

Let us start with rewriting eq.\ (\ref{4.2}) in the variable $v:= \vih^\al r^\eta  u$, $\al  \in \real$, $\eta \in \real$, where the factor $\vih^\al$ serves to shift $a_\rho$ by $\al$ and the factor $r^\eta$ to assure convergence in $Q_\infty$. $\vih$ denotes the tent function (see appendix A). We obtain 
$$
\ba{l}
\disp \int_1^R\bigg[ \int_0^{\pi/2}\bigg(\Big(\pa_r v - \eta\, \frac{v}{r} \Big) \pa_r \psi +\frac{1}{r^2} \Big(\pa_\vi v - \al\, \frac{v}{\vih} \Big) \pa_\vi \psi \bigg) \vih^{-\al} \rd \vi \\[2ex]
\disp \qquad \qquad + \int_{\pi/2}^\pi \bigg(\Big(\pa_r v - \eta\, \frac{v}{r} \Big) \pa_r \psi +\frac{1}{r^2} \Big(\pa_\vi v + \al\, \frac{v}{\vih} \Big) \pa_\vi \psi \bigg) \vih^{-\al} \rd \vi \bigg] r^{1 - \eta} \rd r\\[3ex]
\disp \qquad + \int_1^R\bigg[ \int_0^{\pi/2}\frac{v}{\tan \vih}\Big(a_\zeta \pa_r \psi  + a_\rho \frac{1}{r} \pa_\vi \psi \Big) \vih^{-\al} \rd \vi \\[2ex]
\disp \qquad \qquad \qquad - \int_{\pi/2}^\pi \frac{v}{\tan \vih}\Big(a_\zeta \pa_r \psi  + a_\rho \frac{1}{r} \pa_\vi \psi \Big) \vih^{-\al} \rd \vi \bigg] r^{ - \eta} \rd r\\[3ex]
\disp  \qquad \qquad + \int_1^R \int_0^\pi v \Big(a_\rho \pa_r \psi  - a_\zeta \frac{1}{r} \pa_\vi \psi \Big) \vih^{-\al} \rd \vi \, r^{ - \eta} \rd r\\[3ex]
\disp \qquad \qquad \qquad \qquad = \int_1^R \int_0^\pi \frac{\Om}{\sin \vi} \Big(a_r \pa_r \psi  + a_\vi \frac{1}{r} \pa_\vi \psi \Big)  \rd \vi \, \rd r\, .
\ea
$$
Rearranging the first two terms and using the notation introduced in appendix A, eq.\ (\ref{4.2}) takes the form
\be \label{4.3}
\ba{l}
\disp \int_{Q_R} \na v \cdot \na \psi\, \rd \mu_\al^\eta \\[2ex]
\disp \qquad + \int_1^R \int_0^{\pi/2} \frac{v}{r\, \vih} \bigg(\Big(\frac{\vih}{\tan \vih}\, a_\ze - \eta \,\vih\Big) \pa_r \psi + \Big(\frac{\vih}{\tan \vih}\, a_\rho - \al \Big) \frac{1}{r}\pa_\vi \psi \bigg) \rd \mu_\al ^\eta \\[2ex]
\disp \qquad - \int_1^R \int_{\pi/2}^\pi \frac{v}{r\, \vih} \bigg(\Big(\frac{\vih}{\tan \vih}\, a_\ze + \eta \,\vih\Big) \pa_r \psi + \Big(\frac{\vih}{\tan \vih}\, a_\rho - \al \Big) \frac{1}{r}\pa_\vi \psi \bigg)  \rd \mu_\al ^\eta \\[3ex]
\disp  \qquad \qquad + \int_{Q_R} \frac{v}{r} \Big(a_\rho \pa_r \psi  - a_\zeta \frac{1}{r} \pa_\vi \psi \Big) \rd \mu_\al^\eta \\[2ex]
\disp \qquad \qquad \qquad \qquad = \int_{Q_R} \frac{\Om}{r \sin \vi} \Big(a_r \pa_r \psi  + a_\vi \frac{1}{r} \pa_\vi \psi \Big)  \rd \mu\, .
\ea
\ee
The effect of the ($\vih$-part of the) transformation shows up in an improved $a_\rho$-estimate. By (\ref{4.1}) one finds the optimal estimate
\be \label{4.3a}
\Big\|\frac{ \vih}{\tan \vih}\, a_\rho - \al\Big\|_\infty \leq \|a_\rho - \al \|_\infty = \|a_\rho^\al\|_\infty \leq 0.61 \qquad \mbox{for} \quad  \al = 0.39\, ,
\ee
where we have defined $a_\rho^\al := a_\rho - \al$.

The bilinear form on the left-hand side of (\ref{4.3}) is abbreviated by $B[v, \psi]$ with its three components (in this order) $B_0$, $B_1$, and $B_2$, and the linear form on the right-hand side by $L[\psi]$; thus (\ref{4.3}) reads for short:
\be \label{4.4}
B[v, \psi] = B_0[v, \psi] + B_1 [v, \psi] + B_2 [v, \psi] = L[\psi] \, .
\ee
Morevover, we define an extension of $B$ by
\be \label{4.5}
B_{ex} [v, \psi] := B[v, \psi] + M B_3[v, \psi]
\ee
with 
\be \label{4.6}
B_3 [v, \psi] := \int_{Q_R} v\, \psi \,\frac{1}{r^2}\, \rd \mu_\al^\eta
\ee
and some constant $M > 0$ to be determined later.

Solving eq.\ (\ref{4.3}) in weighted spaces, especially when the variables $v$ and $\psi$ lie in different spaces, require a generalized version of the Lax-Milgram criterion. The following proposition presents  a suitable version, whose proof can be found in \cite{KR22}. 

\vspace{.3em}
\noi{\bf Proposition 4.1 (generalized Lax-Milgram criterion).}
{\em
Let $\BC$ and $\BtC$ be reflexive Banach spaces and $B : \BC \ti \BtC \ra \real$ be a continuous bilinear form, i.e.\ for some $K> 0$ holds
\be \label{4.7}
B(u, \ut) \leq K\, \|u\|_\BC \|\ut\|_{\BtC} \qquad \mbox{for all }\, u \in \BC\, ,\; \ut \in \BtC\, .
\ee
Let, furthermore, $\BtC'$ be the dual space of $\BtC$ with norm
\be \label{4.8}
\|\lt \|_{\BtC'} := \sup_{0\neq \ut \in \BtC} \frac{\lt(\ut)}{\|\ut\|_\BtC}\, ,
\qquad \lt \in \BtC' .
\ee
Let, finally, $c> 0$ and $\ct >0$ be constants such that 
 \be \label{4.9}
 \sup_{0\neq \ut \in \BtC} \frac{B(u , \ut)}{\|\ut\|_\BtC} \geq c \, \|u\|_\BC \qquad \mbox{for all }\, u \in \BC\, ,
 \ee
 \be \label{4.10}
 \sup_{0\neq u \in \BC} \frac{B(u , \ut)}{\|u\|_\BC} \geq \ct \, \|\ut\|_\BtC \qquad \mbox{for all }\, \ut \in \BtC\, .
 \ee
 Then, the equation 
 \be \label{4.11}
 B(u, \ut) = \lt( \ut) \qquad \mbox{for all }\, \ut \in \BtC
 \ee
 with $\lt \in \BtC'$ has a unique solution $u \in \BC$ with bound 
 \be \label{4.12}
 \|u \|_\BC \leq \frac{1}{c}\, \|\lt \|_{\BtC'}\, .
 \ee
 } 
 \vspace{.3em}

To make use of the criterion we set 
(see appendix A)\footnote{In this section $\HC_{\al + \ga}$, etc.\ always refer to $Q_R$ without special mention of this fact.}
$$
\ba{c}
\BC := \HC_{\al + \ga} \; \mbox{ with norm } \; 
\| \na_{\!d} \cdot \|_{(\al + \ga , 3\eta)} \\[2ex]
\BtC := \HC_{\al - \ga} \;\mbox{ with norm } \; \| \na_{\!e} \cdot \|_{(\al - \ga , -\eta)}
\ea
$$
and apply it to the bilinear form $B_{ex} [v, \psi]$ given by (\ref{4.3})--(\ref{4.6}). In the subsequent estimates we assume some parameters fixed, viz.,
\be \label{4.13}
\al = 0.39\, , \qquad \eta = 10^{-3} ,
\ee
whereas others are not yet fixed but restricted:
\be \label{4.14}
\ga \geq \al\, ,\quad 0 < d\leq 1\, ,\quad e\geq 1\, , \quad M> 0\, .
\ee
Condition (\ref{4.7}) is the easy part under the requirements of proposition 4.1 and may be checked term by term. Just by Cauchy-Schwarz's inequality one obtains for $B_0$:
$$
\ba{l}
\disp \big|B_0 [v, \psi] \big| = \bigg| \int_1^R \int_0^\pi \na v \cdot \na \psi \, \vih^{-\al} \rd \vi\, r^{1 - \eta} \rd r\bigg| \\[3ex] 
\disp \qquad \qquad\;\, \leq \int_1^R \int_0^\pi |\na v|\, \vih^{-(\al + \ga)/2}\, r^{-\eta}\, | \na \psi|\, \vih^{-(\al - \ga)/2}\, r^\eta \, \rd \vi\, r^{1 - \eta} \rd r \\[3ex]
\disp \qquad \qquad \qquad \;\,\leq \bigg(\int_1^R \int_0^\pi |\na v|^2\, \vih^{-(\al + \ga)} \rd \vi\, r^{1-3 \eta} \rd r \bigg)^{1/2}\\[2ex]
\disp \qquad \qquad \qquad \qquad \qquad \qquad \ti \bigg(\int_1^R \int_0^\pi |\na \psi|^2 \,\vih^{-(\al - \ga)} \rd \vi\, r^{1+ \eta} \rd r \bigg)^{1/2} \\[3ex]
\disp \qquad \qquad \qquad \qquad \;\, \leq \frac{1}{d}\, \|\na_{\!d}\, v\|_{(\al + \ga, 3\eta)}\, \|\na_{\!e} \psi\|_{(\al - \ga, -\eta)}\, . 
\ea
$$
A rough estimate by (\ref{4.1}), (\ref{4.3a}), and Hardy's inequality (\ref{A.2}) yields for $B_1$:
$$
\ba{l}
\disp \big|B_1 [v, \psi]\big| \leq \int_1^R \bigg[\int_0^{\pi/2} \Big|\frac{v}{r\,\vih}\Big| \Big(|\pa_r \psi| + \Big| \frac{1}{r}\, \pa_\vi \psi\Big|\Big) \vih^{-\al} \rd \vi \\[2ex]
\disp \qquad \qquad \qquad \qquad  +
\int_{\pi/2}^\pi \Big|\frac{v}{r\,\vih}\Big| \Big(|\pa_r \psi| + \Big| \frac{1}{r}\, \pa_\vi \psi\Big|\Big) \vih^{-\al} \rd \vi \bigg] r^{1-\eta} \rd r \\[3ex]
\disp \qquad \qquad \quad \leq \sqrt{2} \int_1^R \int_0^{\pi} \Big|\frac{v}{r\,\vih}\Big| |\na_{\!e} \psi| \,\vih^{-\al} \rd \vi\, r^{1 - \eta} \rd r \\[3ex]
\disp \qquad \qquad \qquad \quad \leq \sqrt{2}\, \Big\|\frac{1}{r}\, \frac{v}{\vih} \Big\|_{(\al + \ga, 3\eta)}\, \|\na_{\!e} \psi\|_{(\al - \ga, -\eta)} \\[2ex] 
\disp \qquad \qquad \qquad \qquad \quad \leq \frac{2 \sqrt{2}}{1 + \al + \ga}\, \|\na_{\!d}\, v\|_{(\al + \ga, 3\eta)}\, \|\na_{\!e} \psi\|_{(\al - \ga, -\eta)}\, .
\ea
$$
Similarly one obtains by Poincar\' e's inequality (\ref{A.3}) for $B_2$:
$$
\ba{l}
\disp \big|B_2 [v, \psi]\big| \leq \sqrt{2} \int_1^R \int_0^{\pi} \Big|\frac{v}{r}\Big| |\na_{\!e} \psi| \,\vih^{-\al} \rd \vi\, r^{1 - \eta} \rd r \\[3ex]
\disp \qquad \qquad \;\, \leq \sqrt{2}\, \Big\|\frac{v}{r} \Big\|_{(\al + \ga, 3\eta)}\, \|\na_{\!e} \psi\|_{(\al - \ga, -\eta)} \\[2ex] 
\disp \qquad \qquad \qquad \;\, \leq \frac{\pi}{2}\, \frac{2 \sqrt{2}}{1 + \al + \ga}\, \|\na_{\!d}\, v\|_{(\al + \ga, 3\eta)}\, \|\na_{\!e} \psi\|_{(\al - \ga, -\eta)}\, ,
\ea
$$
and for $B_3$:
\be \label{4.14a}
\ba{l}
\disp \big|B_3 [v, \psi]\big| \leq \Big\|\frac{v}{r}\Big\|_{(\al + \ga, 3\eta)}\, \Big\|\frac{\psi}{r}\Big\|_{(\al - \ga, -\eta)} \\[2ex] 
\disp \qquad \qquad \;\, \leq \frac{\pi^2}{(1 + \al + \ga)(1+ \al - \ga)}\, \|\na_{\!d}\, v\|_{(\al + \ga, 3\eta)}\, \|\na_{\!e} \psi\|_{(\al - \ga, -\eta)}\, .
\ea
\ee

The harder part are the coercivity conditions (\ref{4.9}) and (\ref{4.10}), which require constants $c>0$ and $\ct> 0$ such that
\be \label{4.15}
\inf_{0 \neq v \in \HC_{\al + \ga}}\, \sup_{0 \neq \psi \in \HC_{\al - \ga}} \frac{ B_{ex} [v, \psi]}{\| \na_{\!d} \, v\|_{(\al + \ga, 3\eta)}\, \|\na_{\!e} \psi \|_{(\al - \ga,-\eta)}} \geq c 
\ee
\be \label{4.16}
\inf_{0 \neq \psi \in \HC_{\al - \ga}}\, \sup_{0 \neq v \in \HC_{\al + \ga}} \frac{ B_{ex} [v, \psi]}{\| \na_{\!d} \, v\|_{(\al + \ga, 3\eta)}\, \|\na_{\!e} \psi \|_{(\al - \ga,-\eta)}} \geq \ct\, , 
\ee
respectively.

Let us introduce new variables, viz.,
\be 
\label{4.17}
w:= \vih^{-\beta_+} r^{-3 \eta/2} v\, ,\quad \chi := \vih^{-\beta_-} r^{ \eta/2} \psi\, ,\qquad \beta_\pm := (\al \pm \ga)/2\, ,
\ee
which have the advantage that by (\ref{A.5}) both variables are in $\HC_0$. This allows us to replace the min-max problems by simpler (but rougher) minimization problems. In fact, expressing the variational expression in (\ref{4.15}) by $w$ and $\chi$ and replacing the maximizing $\chi$ by $\chi := w$ one obtains 
\be \label{4.18}
\inf_{0 \neq w \in \HC_0} \frac{ B_{ex} [\vih^{\beta_+} r^{3 \eta/2} w, \vih^{\beta_-} r^{-\eta/2} w]}{\| \na_{\!d} ( \vih^{\beta_+} r^{3 \eta/2} w) \|_{(2 \beta_+, 3\eta)}\, \|\na_{\!e} (\vih^{\beta_-} r^{-\eta/2} w) \|_{(2 \beta_-,-\eta)}} \geq c\, . 
\ee
A solution of problem (\ref{4.18}), i.e.\ a constant $c>0$ satisfying inequality (\ref{4.18}), is clearly sufficient for problem (\ref{4.15}) and by an analogous argument also for problem (\ref{4.16}). 

In $B_{ex}$ the indefinite term $B_1$ must be dominated by the definite term $B_0$, which requires an as sharp as possible upper bound on $B_1$. By (\ref{4.3a}) and using the notation (\ref{4.17}) and
\be \label{4.19}
\vih^{\beta_-} r^{-\eta/2} w = \vih^{-\ga} r^{-2 \eta} v =: \vt\, , 
\ee
we obtain:
\be \label{4.20}
\ba{l}
\disp \big|B_1 [v, \vt]\big| \leq \int_1^R \bigg[ \int_0^{\pi/2} \Big|\frac{v}{r\,\vih}\Big| \bigg(\Big(\|a_\ze\|_\infty + \eta\, \frac{\pi}{2}\Big) |\pa_r \vt| + \|a_\rho^\al\|_\infty \Big| \frac{1}{r}\, \pa_\vi \vt \Big|\bigg) \vih^{-\al} \rd \vi \\[2ex]
\disp \qquad \qquad + \int_{\pi/2}^\pi \Big|\frac{v}{r\,\vih}\Big| \bigg(\Big(\|a_\ze\|_\infty + \eta\, \frac{\pi}{2}\Big) |\pa_r \vt| + \|a_\rho^\al\|_\infty \Big| \frac{1}{r}\, \pa_\vi \vt \Big|\bigg) \vih^{-\al} \rd \vi 
 \bigg] r^{1-\eta} \rd r \\[3ex]
\disp \qquad \qquad \quad \leq \|a_\rho^\al\|_\infty \int_1^R \int_0^{\pi} \Big|\frac{v}{r\,\vih}\Big| \vih^{-\beta_+} r^{- 3\eta/2} \\[3ex]
\disp \hspace{24ex} \ti \bigg( \frac{\|a_\ze \|_\infty + \eta\, \pi/2}{\|a_\rho^\al \|_\infty }\, |\pa_r \vt| + \Big|\frac{1}{r} \pa_\vi \vt\Big|\bigg) \vih^{-\beta_-} r^{\eta/2}\, \rd \vi\, r \rd r \\[3ex]
\disp \qquad \qquad \qquad \quad \leq \sqrt{2}\, \|a_\rho^\al\|_\infty \Big\|\frac{1}{r}\, \frac{v}{\vih} \Big\|_{(2 \beta_+, 3\eta)}\, \|\na_{\!e} \vt \|_{(2 \beta_-, -\eta)} \\[3ex] 
\disp \qquad \qquad \qquad \qquad \quad \leq \frac{2 \sqrt{2}\, \|a_\rho^\al\|_\infty}{1 + 2\beta_+}\, \|\na_{\!d}\, v\|_{(2 \beta_+, 3\eta)}\, \|\na_{\!e} \vt\|_{(2 \beta_-, -\eta)}\, ,
\ea
\ee
were we made use of (\ref{A.2}) and where we have assumed
\be \label{4.21}
e \geq \frac{\|a_\ze \|_\infty + \eta\, \pi/2}{\|a_\rho^\al \|_\infty }\, .
\ee
The lower order term $B_2$ is dominated by $B_0$ and $B_3$, which allows a rougher estimate. By H\"older's inequality we get:
\be \label{4.22}
\ba{l}
\disp \big|B_2 [v, \vt]\big| \leq \sqrt{2} \int_{Q_R} \Big|\frac{v}{r}\Big| |\na \vt| \, \rd \mu_\al^\eta \\[3ex]
\disp \quad = \sqrt{2} \int_1^R \int_0^{\pi} \Big|\frac{v}{r}\Big|\, \vih^{-\beta_+} r^{- \eta}\, \big(|\na \vt|\, \vih^{-\beta_-} r^\eta\big)^{1/2}\big(|\na \vt|\, \vih^{-\beta_-} r^\eta\big)^{1/2} \,\rd \vi\, r^{1- \eta} \rd r \\[3ex]
\disp \qquad \quad \leq \sqrt{2}\, \Big\|\frac{v}{r}\Big\|_{(2 \beta_+, 3\eta)}\, \|\na_{\!e} \vt \|^{1/2}_{(2 \beta_-, -\eta)} 
\big(C(\al, \ga, \eta, d) \|\na_{\!d} v\|_{(2 \beta_+, 3\eta)} \big)^{1/2}\\[3ex] 
\disp \qquad \qquad \quad \leq \frac{C(\al, \ga, \eta, d)}{\eps}\, \Big\|\frac{v}{r}\Big\|^2_{(2 \beta_+ , 3\eta)} + \frac{\eps}{2}\, \|\na_{\!d}\, v\|_{(2 \beta_+, 3\eta)}\, \|\na_{\!e} \vt\|_{(2 \beta_-, -\eta)}\, ,
\ea
\ee
where the third factor has been rewritten by (\ref{A.2}) as follows:
\be \label{4.22a}
\ba{l}
\disp \int_1^R \int_0^\pi \big|\na (\vih^{-\ga} r^{-2\eta} v)\big|^2\, \vih^{-(\al - \ga)} r^{2 \eta} \, \rd \vi \, r^{1- \eta} \rd r \\[2ex]
\disp \qquad \leq \int_1^R \int_0^\pi \bigg((\pa_r v)^2 + \Big(\frac{1}{r} \pa_\vi v\Big)^2 + \Big( 2 \, \eta \frac{v}{r}\Big)^2 + \Big(\frac{\ga}{r}\, \frac{v}{\vih}\Big)^2\bigg) \vih^{-2 \beta_+} \rd \vi\, r^{1-3 \eta}\rd r \\[3ex]
\disp \qquad \qquad \leq 2 \max \bigg\{ \frac{1}{d^2}\, , \bigg( 1 +(\pi^2\eta^2 + \ga^2) \Big( \frac{2}{1 + \al + \ga}\Big)^2 \bigg)\bigg\}\, \| \na_{\!d}\, v\|^2_{(2 \beta_+, 3 \eta)} \\[3ex]
\disp \qquad \qquad \qquad =: \big(C(\al, \ga, \eta, d) \|\na_{\!d} v\|_{(2 \beta_+, 3\eta)} \big)^2 .
\ea
\ee
Note that the first term in (\ref{4.22}) is just $B_3[v, \vt]$:
\be \label{4.23}
B_3 [v, \vt] = \int_{Q_R} v\, \vih^{- \ga} r^{-2 \eta}\, v\, \frac{1}{r^2}\, \rd \mu_\al^\eta = \Big\|\frac{v}{r}\Big\|^2_{(2 \beta_+, 3 \eta)}\, .
\ee
Concerning $B_0$ we aim at an as sharp as possible lower bound expressed in the variational variable $w$; only terms multiplied by the small parameter $\eta$ allow unsparing estimates. We obtain
\be \label{4.24}
\ba{l}
\disp B_0 [\vih^{\beta_+} v^{3 \eta/2} w, \vih^{\beta_-} v^{-\eta/2} w ] = \int_{Q_R} \na (\vih^{\beta_+} v^{3 \eta/2} w) \cdot \na(\vih^{\beta_-} v^{-\eta/2} w )\, \rd \mu_\al^\eta \\[3ex]
\disp \qquad = \int_{Q_R} \Big(\pa_r w + \frac{3}{2}\, \eta \frac{w}{r}\Big) \Big(\pa_r w - \frac{1}{2}\, \eta \frac{w}{r}\Big) \rd \mu \\[3ex]
\disp \qquad \qquad + \int_1^R \frac{1}{r^2}\bigg[\int_0^{\pi/2} \Big(\pa_\vi w + \beta_+ \frac{w}{\vih}\Big) \Big(\pa_\vi w + \beta_- \frac{w}{\vih}\Big) \rd \vi \\[2ex]
\disp \qquad \qquad \qquad \qquad + \int_{\pi/2}^\pi \Big(\pa_\vi w - \beta_+ \frac{w}{\vih}\Big) \Big(\pa_\vi w - \beta_- \frac{w}{\vih}\Big) \rd \vi \bigg] r \rd r\, .
\ea
\ee
The first integral can be estimated by Cauchy-Schwarz and (\ref{A.3}):
\be \label{4.25}
\ba{l}
\disp \int_{Q_R} \bigg( (\pa_r w)^2 + \eta \, \frac{w}{r} \pa_r w - \frac{3}{4} \, \eta^2\Big(\frac{w}{r}\Big)^2 \bigg) \rd \mu \\[3ex]
\disp \qquad \geq \Big(1 - \frac{\eta}{2}\Big) \int_{Q_R} (\pa_r w)^2\, \rd \mu - \Big( \frac{\eta}{2} + \frac{3}{4}\, \eta^2 \Big) \int_{Q_R} \Big( \frac{w}{r}\Big)^2 \rd \mu \\[3ex]
\disp \qquad \quad\; \geq \Big(1 - \frac{\eta}{2}\Big) \int_{Q_R} (\pa_r w)^2\, \rd \mu - \frac{\eta}{2} \Big(1 +  \frac{3}{2}\, \eta \Big) \pi^2 \int_{Q_R} \Big( \frac{1}{r}\pa_\vi w\Big)^2 \rd \mu\, .
\ea
\ee
As to the second integral note the identity due to integration by parts:
$$
\ba{l}
\disp \int_0^{\pi/2} \Big(\pa_\vi w + \beta_+ \frac{w}{\vih}\Big) \Big(\pa_\vi w + \beta_- \frac{w}{\vih}\Big) \rd \vi \\[2ex]
\disp \qquad = \int_0^{\pi/2} \bigg( (\pa_\vi w)^2 + (\beta_+ + \beta_-) \frac{w}{\vih}\, \pa_\vi w + \beta_+ \beta_- \Big(\frac{w}{\vih}\Big)^2 \bigg) \rd \vi \\[3ex]
\disp \qquad \qquad = \int_0^{\pi/2} (\pa_\vi w)^2 \, \rd \vi + \Big( \beta_+ \beta_- + \frac{1}{2} (\beta_+ + \beta_-) \Big) \int_0^{\pi/2} \Big(\frac{w}{\vih}\Big)^2 \rd \vi \\[2ex]
\disp \qquad \qquad \qquad \quad + \frac{1}{2} (\beta_+ + \beta_-) \frac{2}{\pi}\,  w(\, \cdot\, , \pi/2)^2
\ea
$$
and an analogous expression for the complementary $\vi$-integral. This yields for the second integral in (\ref{4.24}):
\be \label{4.26}
\ba{ll}
\disp  \int_1^R \frac{1}{r^2} \Big[ \ldots \Big]\, r \rd r & \disp = \int_{Q_R} \Big(\frac{1}{r} \pa_\vi w\Big)^2 \rd \mu \\[2ex]
& \disp  \qquad + \Big(\beta_+ \beta_- + \frac{1}{2} (\beta_+ + \beta_-)\Big) \int_{Q_R} \Big(\frac{w}{r \vih}\Big)^2 \rd \mu \\[2ex]  
& \disp  \qquad \qquad +  \frac{2}{\pi} (\beta_+ + \beta_-) \int_1^R \frac{1}{r}\, w(r, \pi/2)^2 \,\rd r\, .
\ea
\ee
Similar (however upper) bounds can be obtained for the denominators in (\ref{4.18}), for instance:
$$
\ba{l}
\disp \big\|\na_{\!d} (\vih^{\beta_+} r^{3 \eta/2} w )\big\|^2_{(2 \beta_+, 3\eta)} = \int_{Q_R} d^2 \Big( \pa_r w + \frac{3}{2} \eta \frac{w}{r} \Big)^2
 \rd \mu \\[2ex]
 \disp \qquad \qquad  + \int_1^R \frac{1}{r^2}\bigg[\int_0^{\pi/2} \Big(\pa_\vi w + \beta_+ \frac{w}{\vih}\Big)^2 \rd \vi + \int_{\pi/2}^\pi \Big(\pa_\vi w - \beta_+ \frac{w}{\vih}\Big)^2  \rd \vi \bigg] r \rd r \\[3ex]
 \disp \qquad \leq \Big(1 + \frac{3}{2} \eta \Big) d^2 \int_{Q_R} (\pa_r w)^2\, \rd \mu + \Big( \frac{3}{2} \eta  + \frac{9}{4} \eta^2 \Big) d^2\int_{Q_R} \Big( \frac{w}{r}\Big)^2 \rd \mu \\[2ex]
\disp \qquad \qquad + \int_{Q_R} \Big(\Big(\frac{1}{r} \pa_\vi w\Big)^2 + \beta_+^2 \Big(\frac{w}{r\, \vih}\Big) ^2\Big) \rd \mu \\[2ex]
\disp \qquad \qquad \qquad + \int_1^R \frac{2 \beta_+}{r^2} \bigg[ \int_0^{\pi/2} \frac{w }{\vih}\, \pa_\vi w \, \rd \vi - \int_{\pi/2}^\pi  \frac{w }{\vih}\, \pa_\vi w \, \rd \vi \bigg] r \rd r\\[3ex]
\disp \qquad  \leq \Big(1 + \frac{3}{2} \eta \Big) d^2 \int_{Q_R} (\pa_r w)^2\, \rd \mu + \Big(1 + \frac{3}{2} \eta\Big(1 + \frac{3}{2} \eta \Big) \pi^2 d^2 \Big) \int_{Q_R} \Big( \frac{1}{r} \pa_\vi w\Big)^2 \rd \mu \\[2ex]
\disp \qquad \qquad \qquad + \beta_+ (1 + \beta_+ ) \int_{Q_R} \Big( \frac{w}{r\, \vih}\Big)^2\rd \mu + \frac{4}{\pi}\, \beta_+ \int_1^R\frac{1}{r} \, w(r,\pi/2)^2\,\rd r\, ,
\ea
$$
and analogously
$$
\ba{l}
\disp \big\|\na_{\!e} (\vih^{\beta_-} r^{- \eta/2} w )\big\|^2_{(2 \beta_-, -\eta)}  \\[2ex]
\disp \qquad  \leq \Big(1 + \frac{\eta}{2}\Big) e^2 \int_{Q_R} (\pa_r w)^2\, \rd \mu + \Big(1 + \frac{\eta}{2} \Big(1 + \frac{\eta}{2} \Big) \pi^2 e^2 \Big) \int_{Q_R} \Big( \frac{1}{r} \pa_\vi w\Big)^2 \rd \mu \\[2ex]
\disp \qquad \qquad \qquad + \beta_- (1 + \beta_- ) \int_{Q_R} \Big( \frac{w}{r\, \vih}\Big)^2\rd \mu + \frac{4}{\pi}\, \beta_- \int_1^R\frac{1}{r} \, w(r,\pi/2)^2\,\rd r\, .
\ea
$$
Introducing the abbreviations 
$$
\ba{c}
\disp \de_\pm := \beta_\pm (1 + \beta_\pm)\, ,\qquad \de_0 := \beta_+ \beta_- + \frac{1}{2}(\beta_+ + \beta_-)\, ,\\[2ex]
\disp c_+ := c_+ (\eta, d) := \frac{3}{2} \eta \Big(1 + \frac{3}{2} \eta\Big) \pi^2 d^2 , \quad 
c_- := c_- (\eta, e) := \frac{1}{2} \eta \Big(1 + \frac{1}{2} \eta\Big) \pi^2 e^2 ,\\[2ex]
\disp c_0 := c_0 (\eta) := \frac{1}{2} \eta \Big(1 + \frac{3}{2} \eta\Big) \pi^2  ,
\ea
$$
and the ratios
$$
s:= \frac{\disp \int_{Q_R} (\pa_r w)^2\, \rd \mu }{\disp  \int_{Q_R} \Big(\frac{1}{r} \pa_\vi w \Big)^2 \rd \mu }\, , 
\quad 
t:= \frac{1}{4} \frac{\disp \int_{Q_R} \Big(\frac{w}{r\, \vih}\Big)^2 \rd \mu }{\disp \int_{Q_R} \Big(\frac{1}{r} \pa_\vi w \Big)^2 \rd \mu }\, , \quad 
\tau := \frac{4}{\pi} \frac{\disp \int_1^R \frac{1}{r}\, w (r, \pi/2)^2\, \rd r }{\disp \int_{Q_R} \Big(\frac{1}{r} \pa_\vi w \Big)^2 \rd \mu }\, ,
$$
a lower bound on $B_0$ is thus given by 
\be \label{4.27}
\frac{ B_0 [\vih^{\beta_+} r^{3 \eta/2} w, \vih^{\beta_-} r^{-\eta/2} w]}{\| \na_{\!d} ( \vih^{\beta_+} r^{3 \eta/2} w) \|_{(2 \beta_+, 3\eta)}\, \|\na_{\!e} (\vih^{\beta_-} r^{-\eta/2} w) \|_{(2 \beta_-,-\eta)}} \geq \inf_{\mbox{\scriptsize{$\ba{c}
0\leq s < \infty \\ 0\leq t \leq 1 \\ 0 \leq \tau \leq 1
\ea$ }}}  F(s,t, \tau)
\ee
with
$$
\ba{l}
\disp F(s,t, \tau) := \bigg(\frac{1 - c_0 + \Big(\disp 1 -\frac{\eta}{2}\Big) s + 4 \, \de_0\, t + \frac{1}{2}(\beta_+ + \beta_-)\, \tau}{\disp 1 + c_+ + \Big(1 + \frac{3}{2}\eta\Big) d^2 s + 4\, \de_+\, t + \beta_+  \tau}\bigg)^{1/2} \\[4ex]
\disp \qquad \qquad \qquad \qquad \ti  \bigg(\frac{\disp 1 - c_0 + \Big(\disp 1 -\frac{\eta}{2}\Big) s + 4 \, \de_0\, t + \frac{1}{2}(\beta_+ + \beta_-)\, \tau}{\disp 1 + c_- + \Big(1 + \frac{1}{2}\eta\Big) e^2 s + 4\, \de_-\, t + \beta_- \tau}\bigg)^{1/2} .
\ea
$$
The range of the variable $t$ is restricted by inequality (\ref{A.2}), whereas the elementary inequalities 
$$
w(\,\cdot\, ,\pi/2)^2 \leq \frac{\pi}{2} \int_0^{\pi/2} (\pa_\vi w)^2\, \rd \vi \, ,\quad  w(\, \cdot\, ,\pi/2)^2 \leq \frac{\pi}{2} \int_{\pi/2}^\pi (\pa_\vi w)^2\, \rd \vi\, ,
$$
and hence
$$
w(\,\cdot\, ,\pi/2)^2 \leq \frac{\pi}{4} \int_0^{\pi} (\pa_\vi w)^2 \rd \vi
$$
restrict the variable $\tau$.

Denoting the lower bound on the right-hand side in (\ref{4.27}) by $L\! B = L\!B (\al, \ga,$ $ \eta, d,e)$ we can finally insert the results of (\ref{4.20}), (\ref{4.22}), (\ref{4.23}) by (\ref{4.4})--(\ref{4.6}) into (\ref{4.18}) to obtain
\be \label{4.28}
\ba{l}
\disp \inf_{w \in \HC_0} \frac{ B_{ex} [\vih^{\beta_+} r^{3 \eta/2} w, \vih^{\beta_-} r^{-\eta/2} w]}{\| \na_{\!d} ( \vih^{\beta_+} r^{3 \eta/2} w) \|_{(2 \beta_+, 3\eta)}\, \|\na_{\!e} (\vih^{\beta_-} r^{-\eta/2} w) \|_{(2 \beta_-,-\eta)}} \\[3ex]
\disp \qquad \qquad \geq L\! B (\al, \ga, \eta, d,e) - C\! F (\al, \ga) - \frac{\eps}{2}\, ,
\ea
\ee
where we have set $M:= C(\al, \ga, \eta, d)/\eps$ and introduced the abbreviation 
$$
C\!F(\al, \ga) := \frac{2 \sqrt{2}\, \| a_\rho - \al\|_\infty}{1 + \al + \ga}\, .
$$
The variational problem on the right-hand side of (\ref{4.27}) and the bound (\ref{4.28}) closely resemble those in \cite{KR22}. They differ in the additional parameters $\eta$ and $\eps$, which, however, both can be chosen to be small, and the additional variational variable $\tau$. The minimization of $F$ is again elementary and the result in terms of graphs of $\ga \mapsto L\!B (\al, \ga, \eta, d, e)$ and of the ``comparison function'' $\ga \mapsto C\!F (\al, \ga)$ is shown in Fig.\ 3. 
\begin{figure}
\begin{center}
\resizebox{!}{.69\textwidth}{\input{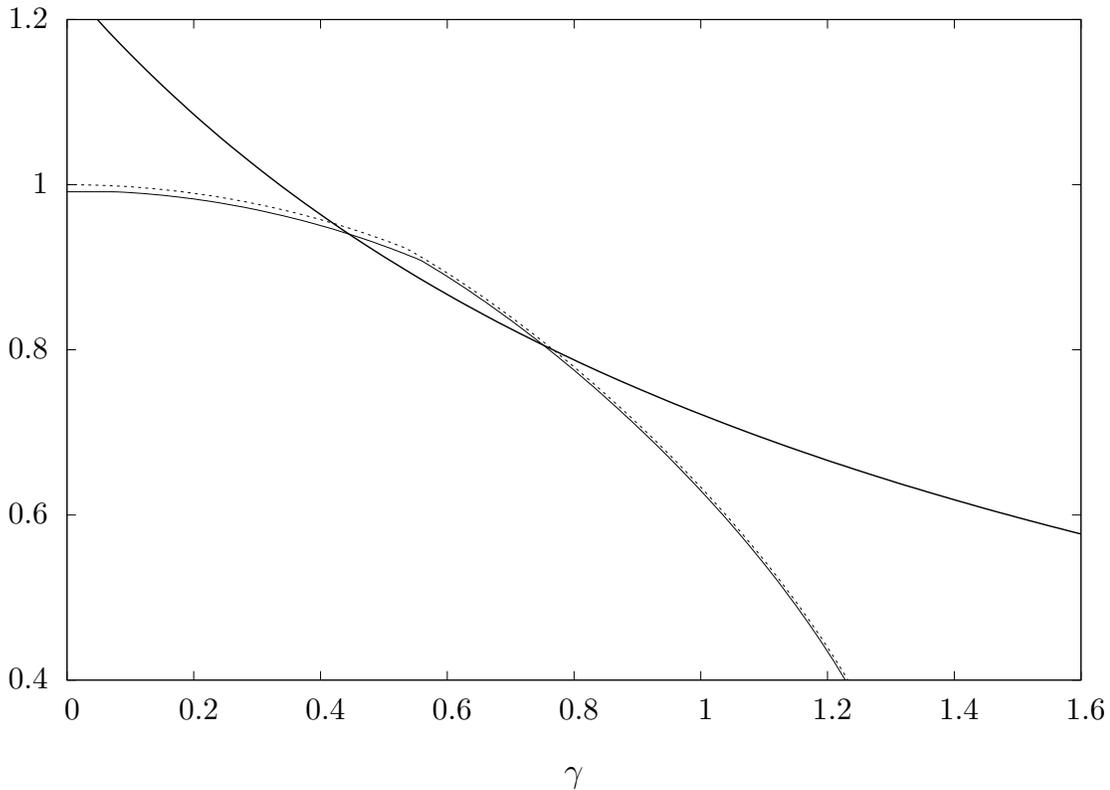}}
\caption{Graphs of the lower bounds $\ga \mapsto L\!B (0.39, \ga, 10^{-3}, 0.1, 1.2)$ (solid thin line),  $\ga \mapsto L\!B (0.39, \ga,0 , 0.1, 1.2)$ (dashed line), and the comparison function $C\!F(0.39, \ga)$ (solid thick line).}
\end{center}
\end{figure}
The parameters $\al$ and $\eta$ have already been fixed in (\ref{4.13}) and $e = 1.2$ is the minimal value compatible with (\ref{4.1}), (\ref{4.3a}), and (\ref{4.21}). $d = 0.1$ turns out to be suitable as in \cite{KR22}. Fig.\ 3 shows again a window of $\ga$-values around $\ga \approx 0.6$, where $L\!B$ exceeds $C\!F$ thus providing (for suitable $\eps >0$) positive coercivity constants in proposition 4.1. To be definite let us fix $\ga =0.6$; we then have
\be \label{4.29}
\De:= L\!B (0.39,\, 0.6,\, 10^{-3},\, 0.1,\, 1.2) - C\!F (0.39,\, 0.6) \geq 0.023 > 0\, ,
\ee
and by taking $\eps := \De$, $\De/2$ is an allowed constant in (\ref{4.18}). 
For comparison Fig.\ 3 also shows the slightly better lower bound from \cite{KR22}, which is just 
$ L\!B (0.39,\, \ga,\, 0,\, 0.1,\, 1.2)$. A closer inspection shows that minimum $F$ is always obtained at $\tau =0$, so that the $\tau$-minimization does in fact not contribute to the minimum; the minimizing values of $s$ and $t$, however, depend on $\ga$ as in \cite{KR22}.

So far proposition 4.1 guarantees the existence of a unique solution $v_{ex} \in \HC_{\al + \ga}$ of the extended equation 
\be \label{4.30}
B_{ex} [v_{ex}, \psi] = L[\psi] \qquad \mbox{ for all }\,\; \psi \in \HC_{\al - \ga}
\ee
with bound (according to (\ref{4.8}), (\ref{4.11}), and (\ref{4.12}))
\be \label{4.31}
\ba{l}
\disp \|\na_{\!d}\, v_{ex}\|_{(\al + \ga, 3 \eta)} 
\leq \frac{2}{\De} \sup_{0 \neq \psi \in \HC_{\al - \ga} } \frac{L[\psi]}{\|\na_{\!e} \psi\|_{(\al - \ga, - \eta)}} \\[3ex]
\disp \qquad \qquad = \frac{2}{\De} \sup_{0 \neq \psi \in \HC_{\al - \ga} } \frac{\disp \int_{Q_R} \frac{\Om}{r \sin \vi} \Big(a_r\, \pa_r \psi + a_\vi \frac{1}{r} \pa_\vi\psi \Big) \rd \mu}{\|\na_{\!e} \psi\|_{(\al - \ga, - \eta)}}\\[3ex]
\disp \qquad \qquad \qquad \qquad \leq \frac{2\, \sqrt{2}}{\De} \, \Big\|\frac{\Om}{r \sin \vi} \Big\|_{(\ga - \al, \eta)}\, .
\ea
\ee
%Note that this bound does not depend on $R$.
Equation (\ref{4.30}) defines a linear and by (\ref{4.31}) bounded operator
\be \label{4.31aa}
\KC_M : \HC_{\al - \ga}' \ra \HC_{\al + \ga}\, , \quad L\mapsto v_{ex}\, ,
\ee
where $\HC_{\al -\ga}'$ denotes the dual space of $\HC_{\al - \ga}$. Defining the weighted $L^2$-space 
$$
\LC^2_{\al + \ga} := {\rm clos} \big( C (Q_R)\, , \, \| \cdot \|_{(\al + \ga, 3 \eta)} \big)\, ,
$$
the operator $\KC_M : \LC^2_{\al + \ga} \ra \LC^2_{\al + \ga}$ is, moreover, compact due to the well-known compact embedding $\HC_{\al + \ga} \Subset \LC^2_{\al + \ga}$ on bounded domains. 

The equation to be solved, eq.\ (\ref{4.3}), reads in the notation (\ref{4.4})--(\ref{4.6}):
\be \label{4.31a}
B_{ex} [v, \psi] = L[\psi] + M B_3[v, \psi] \qquad \mbox{ for all }\,\; \psi \in \HC_{\al - \ga}\, ,
\ee
which corresponds to the operator equation
\be \label{4.32}
v = \KC_M \Big(L + M \frac{v}{r^2}\Big) \quad \Longleftrightarrow \quad \Big(\IC - M\KC_M\, \frac{1}{r^2}\Big) v = \KC_M L\, .
\ee
The multiplication by $1/r^2$ is a bounded operator on $Q_R$ and $\KC_M\, r^{-2}$ is thus compact, too. 
The operator $\IC - M \KC_M\, r^{-2} =: \FC_M$ on the Hilbert space $\LC^2_{\al + \ga}$ is then of Fredholm type, and Fredholm's alternative yields a unique solution of eq.\ (\ref{4.32})$_b$ for any $\LC^2_{\al + \ga}$-right-hand side, provided that the homogeneous equation admits only the trivial solution. This is guaranteed by proposition 3.3, which proves triviality of the solution in the homogeneous version of eq.\ (\ref{4.2}).

An $\LC^2_{\al + \ga}$-bound on this solution $v$ is easily obtained by recalling that $\FC_M$ has a bounded inverse (see, e.g., \cite{GT98}, p.\ 83f.):
\be \label{4.32a}
\sup_{0 \neq u \in \LC^2_{\al + \ga}} \frac{\|\FC_M^{-1} u \|_{(\al + \ga, 3 \eta)}}{\| u\|_{(\al + \ga, 3\eta)}} = : C_F < \infty\, .
\ee
More appropriate, however, is the operator 
$$
\widetilde{\FC}_M := \frac{1 }{r}\,\FC_M \, r :\, \LC^2_{\al + \ga} \ra \LC^2_{\al + \ga}\ ,
$$
which is also Fredholm with inverse constant $\widetilde{C}_F$. By (\ref{A.3}) and (\ref{4.31}) we then obtain for $v/r$:
\be \label{4.33}
\ba{l}
\disp \Big\|\frac{v}{r}\Big\|_{(\al + \ga, 3 \eta)} = \| \widetilde{\FC}_M^{-1} \frac{1}{r}\, \KC_M L\|_{(\al + \ga, 3\eta)} \leq \widetilde{C}_F \Big\|\frac{v_{ex}}{r} \Big\|_{(\al + \ga, 3 \eta)} \\[3ex]
\disp \qquad \leq \widetilde{C}_F\, \frac{\pi}{1 + \al + \ga} \, \| \na_{\!d}\, v_{ex} \|_{(\al + \ga, 3 \eta)} \\[3ex]
\disp \qquad \qquad 
\leq \widetilde{C}_F \,\frac{\pi}{1 + \al + \ga} \, \frac{2}{\De}\,\sup_{0 \neq \psi \in \HC_{\al - \ga} } \frac{L[\psi] }{\|\na_{\!e} \psi\|_{(\al - \ga, - \eta)}} \, .
\ea
\ee
To obtain a gradient-bound on $v$ we replace in the estimate (\ref{4.31}) eq.\ (\ref{4.30}) by eq.\ (\ref{4.31a}) to obtain 
$$
\ba{l}
\disp \|\na_{\!d}\, v\|_{(\al + \ga, 3 \eta)} 
\leq \frac{2}{\De} \sup_{0 \neq \psi \in \HC_{\al - \ga} } \frac{L[\psi] + M B_3 [v, \psi]}{\|\na_{\!e} \psi\|_{(\al - \ga, - \eta)}} \\[3ex]
\disp \qquad \qquad \leq \frac{2}{\De} \bigg[ \sup_{0 \neq \psi \in \HC_{\al - \ga} } \frac{L[\psi] }{\|\na_{\!e} \psi\|_{(\al - \ga, - \eta)}} + M \sup_{0 \neq \psi \in \HC_{\al - \ga} } \frac{ B_3 [v, \psi]}{\|\na_{\!e} \psi\|_{(\al - \ga, - \eta)}} \bigg]\, .
\ea
$$
Splitting $B_3$ similarly to (\ref{4.14a}) and using (\ref{4.33}) we have finally
\be \label{4.34}
\ba{l}
\disp \|\na_{\!d} \, v\|_{(\al + \ga, 3 \eta)} \leq C_1 \sup_{0 \neq \psi \in \HC_{\al - \ga} } \frac{L[\psi] }{\|\na_{\!e} \psi\|_{(\al - \ga, - \eta)}} 
\\[3ex]
\disp \qquad \qquad \qquad \quad \leq C_1\, \sqrt{2}\, \Big\| \frac{\Om}{r \sin \vi} \Big\|_{(\ga - \al ,  \eta)} 
\ea
\ee
with 
\be \label{4.34a}
C_1 := \frac{2}{\De} \Big( 1 + M\, \widetilde{C}_F\, \frac{\pi^2}{(1 + \al + \ga) ( 1 + \al - \ga)} \, \frac{2 }{\De} \Big)\, .
\ee
A weak solution $u\in \HC_0 (Q_R)$ of problem $P_\Om (Q_R)$ corresponds by (\ref{A.4}) to a solution $v \in \HC_{2\al} (Q_R)$ of eq.\ (\ref{4.3}). As $\ga \geq \al$ the above solution is of this kind. 

Note that for the solution of eq.\ (\ref{4.30}) no other properties of the vector field $\aaa$ have been used than the bounds (\ref{4.1}). The solution of eq.\ (\ref{4.31a}), however, relied on the uniqueness result of proposition 3.3, which made use of the particular form (\ref{2.23}) of the vector field.

Concerning the various parameters the choice
\be \label{4.35}
(\al, \ga, \eta, d, e, \De) = (0.39,\, 0.6,\, 10^{-3},\, 0.1,\, 1.2,\, 0.023)
\ee
satisfies all prerequisites of the above described solution procedure, in particular, conditions (\ref{4.3a}), (\ref{4.14}), (\ref{4.21}), (\ref{4.28}), and (\ref{4.29}). $M$ is given by $M= C(\al, \ga, \eta, d)/\De$ with $C(\al, \ga, \eta, d)$ defined in (\ref{4.22a}). So, proposition 4.1 may be applied to obtain 
a solution $v_{ex}$ of eq.\ (\ref{4.30}) and subsequent application of Fredholm's theory yields the desired solution $v$ of eq.\ (\ref{4.31a}).
We summarize this result in the following proposition.

\vspace{.3em}
\noi{\bf Proposition 4.2.}
{\em
Let the vector field $\aaa = a_\ze \eee_\ze + a_\rho \eee_\rho$ be as in (\ref{2.23}) with argument $w - \Om$, where $w \in \HC_0 (Q_R)$ is fixed  and the data function $\Om $ satisfies the bound $\| \Om /(r \sin \vi) \|_{(\ga - \al, \eta)} < \infty$. Let furthermore, the parameters $\al$, $\ga$, $\eta$, $d$, $e$, and $\De$ be fixed according to (\ref{4.35}). Then eq.\ (\ref{4.3}) has a unique solution $v \in \HC_{\al + \ga} (Q_R) \subset \HC_{2 \al} (Q_R)$ satisfying the bound (\ref{4.34}). 
}

\vspace{.3em}
\noi{\em Remark:} The bounding constant $C_1$ as defined by (\ref{4.34a}) does not explicitly depend on the particular form of the vector field $\aaa$ (and hence on its arguments) nor on the domain $Q_R$; implicitly, however, it can depend on these quantities via the Fredholm constant $\widetilde{C}_F$ commonly introduced by a non-constructive argument. In appendix C it is shown that this is not the case: based on spectral estimates of the Fredholm operator $\FCt$ we derive a bound on $\widetilde{C}_F$ that only depends on upper and lower bounds on $\aaa$ and that does not depend on $R$.

\sect{Proofs of theorems 2.1 and 2.2}
The proof of theorem 2.1 is based on Schauder's fixed point principle in a version due to Schaefer (cf.\ \cite{E98}, p.\ 504, Theorem 4), viz.,

\vspace{.3em}
\noi{\bf Proposition 5.1 (fixed point principle).}
{\em
 Let $\TC : \BC \ra \BC$ be a continuous, compact mapping of the Banach space $\BC$ into itself and let the set 
 \be \label{5.1}
 \big\{ v \in \BC : v = \la\, \TC v \mbox{ for some } \la \in [0,1]\big\}
 \ee
 be bounded. Then, $\TC$ has a fixed point.}
\vspace{.3em}

\noi The basic ingredients of the proof are the the same as those employed in \cite{KR22}: The mapping is obtained by solution of the linearized problem with fixed argument in the vector field $\aaa [\ldots]$ according to proposition 4.2. The Banach space is a weighted $L^p$-space, where weight and integration exponent have to be carefully chosen to ensure the mapping to be a continuous endomorphism. Here the lemmata A.2 and A.3 play a crucial role, which are shown to remain valid in the present situation with less boundary control. 
Boundedness of the set (\ref{5.1}) depends on the global bound (\ref{4.34}) and compactness is in fact a standard result. The details of the proof differ, however, from those in \cite{KR22} due to the use of polar coordinates and the necessity of an additional weight. Moreover, the governing equation (\ref{4.3}) is now more complicated and bounds on its solution involve, additionally, Fredholm-type estimates.

Let $\BC := \LC_\de^p (Q_R)$ with
$$
\LC_\de^p (Q_R) = {\rm clos} \big( C (Q_R) \, ,\, \|\cdot\,  \vih^{-\de} \|_{p, (0,3 \eta)} \big)
$$
and 
\be \label{5.2}
\de := (\al + \ga)/2\, ,\quad \eta = 10^{-3},\; \mbox{ and }\; p = 29\, .
\ee
Let, furthermore, $v \in \HC_{\al + \ga} (Q_R)$ be the unique solution, according to proposition 4.2, of eq.\ (\ref{4.3}), where we have set $w:= \vih^{-\al} r^{- \eta} v_0$ for some $v_0 \in \LC^p_\de (Q_R)$ and where $\Om$ satisfies condition (\ref{2.29}). By (\ref{4.34}), (\ref{2.29}) and appendix C this solution satisfies the global bound 
\be \label{5.3}
\|\na_{\!d}\, v\|_{(\al + \ga, 3 \eta)} \leq C_1 \sqrt{2}\, \Big\|\frac{\Om}{\sin \vi} \Big\|_{(\ga - \al, \eta)} \leq K\, C_1 \bigg(\frac{4\, (\pi/2)^{1 + \al -\ga}}{\eta\,(1 + \al - \ga)}\bigg)^{\!\frac{1}{2}} .
\ee
By lemma A.2 with (\ref{5.2})$_a$ we have (in fact for any $p \geq 2$)
\be \label{5.4}
\|v\, \vih^{-\de} \|_{p, (0,3 \eta)} \leq C_>\, \|\na_{\!d}\, v\|_{(\al + \ga, 3 \eta)}\, ,
\ee
hence 
\be \label{5.5}
\HC_{(\al +\ga, 3 \eta)} (Q_R) \subset \LC^p_\de (Q_R)\, ,
\ee
and the mapping
$$
\TC: \LC_\de^p (Q_R) \ra \LC_\de ^p (Q_R)\, ,\quad v_0 \mapsto v 
$$
is well defined .

To prove continuity of $\TC$ let $v$ and $\vt$ be solutions of eq.\ (\ref{4.3}) with different arguments $v_0$ and $\vt_0$ but the same data function $\Om  = \Omt$ in the vector field $\aaa$. Using the abbreviations
$$
\aaa = \aaa[\vih^{-\al} r^{-\eta} v_0  - \Om]\, , \quad 
\aaat = \aaa[\vih^{-\al} r^{-\eta} \vt_0  - \Om]\, ,
$$
we write the difference of (\ref{4.3}) for $v$ and $\vt$ in the form:
\be \label{5.6}
\ba{l}
\disp \int_{Q_R} \na (v - \vt)\cdot \na \psi\, \rd \mu_\al^\eta \\[2ex]
\disp \qquad + \int_1^R \int_0^{\pi/2} \frac{v - \vt}{r\, \vih} \bigg(\Big(\frac{\vih}{\tan \vih}\, a_\ze - \eta \,\vih\Big) \pa_r \psi + \Big(\frac{\vih}{\tan \vih}\, a_\rho - \al \Big) \frac{1}{r}\pa_\vi \psi \bigg) \rd \mu_\al ^\eta \\[2ex]
\disp \qquad - \int_1^R \int_{\pi/2}^\pi \frac{v - \vt}{r\, \vih} \bigg(\Big(\frac{\vih}{\tan \vih}\, a_\ze + \eta \,\vih\Big) \pa_r \psi + \Big(\frac{\vih}{\tan \vih}\, a_\rho - \al \Big) \frac{1}{r}\pa_\vi \psi \bigg)  \rd \mu_\al ^\eta \\[3ex]
\disp  \qquad \qquad + \int_{Q_R} \frac{v - \vt}{r} \Big(a_\rho \pa_r \psi  - a_\zeta \frac{1}{r} \pa_\vi \psi \Big) \rd \mu_\al^\eta \\[3ex]

\disp \;\;  = - \int_1^R \int_0^{\pi/2} \frac{\vt}{r\, \vih} \bigg(\Big(\frac{\vih}{\tan \vih}\, (a_\ze - \ati_\ze) \Big) \pa_r \psi %\\[2ex]
%\disp \qquad \qquad \qquad \qquad \qquad \qquad
+ \Big(\frac{\vih}{\tan \vih}\, (a_\rho - \ati_\rho) \Big) \frac{1}{r}\pa_\vi \psi \bigg) \rd \mu_\al ^\eta \\[2ex]
\disp \qquad + \int_1^R \int_{\pi/2}^\pi \frac{\vt}{r\, \vih} \bigg(\Big(\frac{\vih}{\tan \vih}\, (a_\ze - \ati_\ze) \Big) \pa_r \psi %\[2ex]
%\disp \qquad \qquad \qquad \qquad \qquad \qquad 
+ \Big(\frac{\vih}{\tan \vih}\, (a_\rho - \ati_\rho) \Big) \frac{1}{r}\pa_\vi \psi \bigg)  \rd \mu_\al ^\eta \\[3ex]
\disp  \qquad \qquad - \int_{Q_R} \frac{\vt}{r} \Big((a_\rho - \ati_\rho) \pa_r \psi  - (a_\zeta - \ati_\ze) \frac{1}{r} \pa_\vi \psi \Big) \rd \mu_\al^\eta \\[3ex]
\disp \qquad \qquad \qquad + \int_{Q_R} \frac{\Om}{r \sin \vi} \Big((a_r -\ati_r) \pa_r \psi  + (a_\vi - \ati_\vi) \frac{1}{r} \pa_\vi \psi \Big) \vih^\al r^\eta\, \rd \mu^\eta_\al\, .
\ea
\ee
Denoting the right-hand side of eq.\ (\ref{5.6}) by $\Lt [\psi]$ and using the notation (\ref{4.4}) -- (\ref{4.6}) this equation reads for short
$$
B[v - \vt, \psi] = \Lt[\psi]
$$
or 
$$
B_{ex} [v - \vt, \psi] = \Lt[\psi] + M B_3 [ v - \vt, \psi]\, .
$$
The latter equation is obviously of type (\ref{4.31a}), to which the bound (\ref{4.34}) applies. 
We thus obtain
\be \label{5.7}
\|\na_{\!d} \, (v- \vt) \|_{(\al + \ga, 3 \eta)} \leq C_1 \sup_{0 \neq \psi \in \HC_{\al - \ga} } \frac{\Lt[\psi] }{\|\na_{\!e} \psi\|_{(\al - \ga, - \eta)}}\, .
\ee
The estimate of $\Lt$ proceeds similarly as in ref.\!
\cite{KR22}. Using Lipschitz continuity of $\aaa[\ldots]$, (\ref{2.29}), and H\"older's inequality one obtains
$$
\ba{l}
\disp \big|\Lt [\psi]\big| \leq \int_{Q_R} |v_0 - \vt_0| \bigg(\Big(\frac{|\vt|}{r \tan \vih} + \frac{|\vt|}{r} \Big) \vih^{- \al} r^{- \eta} + \Big|\frac{\Om}{r \sin \vi}\Big| \bigg) |\na \psi|\, \rd \mu_\al^\eta \\[3ex]
\disp \qquad \leq \bigg\| |v_0 - \vt_0| \bigg(\Big( 1 + \frac{\pi}{2}\Big)
\frac{|\vt|}{\vih^{1+ \al}}\, r^{-(1+ \eta)} + \frac{K}{r}\bigg) \bigg\|_{(\al + \ga, 3\eta)}\, \big\|\na \psi\big\|_{(\al - \ga, -\eta)}\\[3ex]
\disp \qquad \qquad \leq \big\| (v_0 - \vt_0)\, \vih^{-(\al + \ga)/2}\big\|_{p, (0, 3 \eta)} \\[2ex]
\disp \qquad \qquad \qquad \ti \bigg(\Big( 1 + \frac{\pi}{2}\Big) \Big\|
\frac{|\vt|}{\vih^{1+ \al}} \Big\|_{q, (0, 3 \eta)} + K \, \Big\|\frac{1}{r}\Big\|_{q, (0, 3\eta)} \bigg) \big\|\na_{\! e} \psi\big\|_{(\al - \ga, -\eta)}
\ea
$$
with $1/p + 1/q = 1/2$. The crucial estimate is again that of $\vt$, which is done by lemma A.3 and which requires $q \lesssim 2.15$ and hence $p \geq 29$:
$$ \Big\|\frac{\vt}{\vih^{1 + \al}}\Big\|_{q, (0, 3 \eta)} \leq C_< \,\big\| \na_{\! d}\, \vt \Big\|_{( \al + \ga, 3 \eta)}\, .
$$
Together with (\ref{5.3}) we thus obtain for $\Lt$:
\be \label{5.8}
\big|\Lt [\psi] \big| \leq C_2\, K\, \big\| (v_0 - \vt_0)\, \vih^{-(\al + \ga)/2}\big\|_{p, (0, 3 \eta)}
\big\|\na_{\! e} \psi\big\|_{(\al - \ga, - \eta)}
\ee
with 
$$
C_2 := \Big(1 + \frac{\pi}{2}\Big) C_<\, C_1 \bigg(\frac{4\, (\pi/2)^{1 + \al - \ga}}{\eta \, (1 + \al - \ga)}\bigg)^{1/2} + \Big(\frac{\pi}{3 \,\eta}\Big)^{1/2} \, .
$$
Finally, applying (\ref{5.4}) to the left-hand side of (\ref{5.7}) and inserting (\ref{5.8}) into the right-hand side we have by (\ref{5.2}):
$$
 \big\| (v - \vt)\, \vih^{-\de}\big\|_{p, (0, 3 \eta)} \leq C_>\, C_1\, C_2\, K\, \big\| (v_0 - \vt_0)\, \vih^{-\de}\big\|_{p, (0, 3 \eta)}\, ,
 $$
 which proves continuity of $\TC$ in $\LC_\de^p (Q_R)$.

 Compactness of the mapping $\TC : \LC_\de^p (Q_R) \ra \LC_\de^p (Q_R)$ is implied by the compactness of the embedding (\ref{5.5}), which is in fact a standard result for any bounded region $\subset \real^2$ and any $1 \leq p < \infty$.
 
 Boundedness of the set (\ref{5.1}) is an immediate consequence of (\ref{5.3}) and (\ref{5.4}):
$$
\ba{l}
\big\|v_0\, \vih^{-\de}\big\|_{p,(0, 3\eta)} \leq C_>\, \big\|\na_{\!d}\, v_0\big\|_{(\al + \ga, 3 \eta)} = C_>\, \big\|\na_{\!d} (\la\, \TC v_0)\big\|_{(\al + \ga, 3\eta)} \\[2ex]
\disp \qquad \qquad \qquad \quad \leq C_>\, \big\|\na_{\!d}\, v\big\|_{(\al + \ga, 3 \eta)} \leq C_>\, K\, C_1 \bigg(\frac{4\, (\pi/2)^{1+ \al - \ga}}{\eta\, (1 + \al - \ga)}\bigg)^{1/2}
\ea
$$
for any $v_0 \in \LC^p_\de (Q_R)$ and any $\la \in [0,1]$ satisfying $v_0 = \la\, \TC v_0$.

With all prerequisites fulfilled proposition 5.1 provides a fixed point $v= v_0$ of eq.\ (\ref{4.3}) with $\aaa = \aaa[\vih^{-\al} r^{-\eta} v_0 - \Om]$. Substituting $u = \vih^{-\al} r^{-\eta} v$ we end up with a solution of eq.\ (\ref{4.2}), equivalently eq.\ (\ref{2.28}), as claimed in theorem 2.1.

The estimate (\ref{2.30}) follows by (\ref{5.3}), (\ref{4.35}), and 
the norm equivalence 
$$
\|\na v \|_{(\al + \ga, 3 \eta)} \sim \|\na u \|_{(\ga - \al, \eta)}\, , 
$$
which holds by (\ref{A.4}). This together with proposition 3.1 concludes the proof of theorem 2.1.

The proof of theorem 2.2 establishing the solution of the unbounded problem $P_\Om (Q_\infty)$ is based on solutions of a suitable sequence of bounded problems $(P_\Om (Q_n))_{n \in \nat\setminus \{1\}}$: Let $\Om  \in L^\infty (Q_\infty)$ with bound (\ref{2.31}) and $\Om_n := \Om|_{Q_n}$. By theorem 2.1 we have a sequence of solutions $(u_n) \subset \HC_0 (Q_n)$ with bound
$$
\|\na u_n \|_{(\beta, \eta)}  \leq \sqrt{C}\, , \qquad n \in \nat\setminus\{1\}\, ,
$$
where the constant $C$ does not depend on $n$. Considering $u_n$ on $Q_m$, $1< m \leq n$ we construct as in \cite{KR22} a sequence of nested subsequences and choose finally the diagonal sequence $(u^{(k)} )_{k \in \nat}$. This sequence has the favourable property that $u^{(k+1)}$ extends $u^{(k)}$ defined on $Q_k$ onto $Q_{k+1}$. Thus $(u^{(k)})_{k \in \nat}$ allows the unique definition of a function $u$ on $Q_\infty$ that satisfies the bound (\ref{2.32}) and that is in fact a weak solution of $P_\Om (Q_\infty)$ (cf.\ \cite{KR22}). Together with proposition 3.2 this proves theorem 2.2.

\sect{Proofs of theorems 2.3 and 2.4}
Theorems 2.3 and 2.4 answer the axisymmetric intensity problems as formulated in the problems $P_{I, \Ih} (A_R)$ and $P_I (A_\infty)$, corresponding to spherical shells and exterior space in $\real^3$, respectively. The existence of weak solutions to both problems follows essentially by theorems 2.1/2.2 and the results gathered in appendix B. Higher interior regularity then follows by standard arguments, whereas boundary regularity, especially continuous assumption of the boundary values, requires more subtle arguments, which differ in parts from those applied to the direction problem.
We present the bounded case in some detail; the unbounded case then follows by obvious modifications. 

The essential prerequisite of theorem 2.1 is the bound (\ref{2.29}) on the data function $\Om$ comprising the zero-positions angle $\Psi$ and the boundary value function $q_l$. Fixing $R > 1$, $\roh \in \nat$, and a (possibly empty) set of points $\{z_1, \ldots , z_{2 N}\} \subset A_R \setminus S\!A$, $N \in \nat_0$, this bound is satisfied by $\Psi$ according to lemma B.1. Fixing, furthermore, H\"older continuous, symmetric intensity functions $I$ and $\Ih$, satisfying the axis condition (\ref{2.33}), the boundary functions $p|_{S_1}$ and $p|_{S_R}$, given by (\ref{2.15}) and (\ref{2.16}), respectively, then satisfy H\"older conditions of type (\ref{B.10}) and (\ref{B.11}). Here
the constant $p_0 =: p_{l,0}$ is fixed by condition (\ref{B.22b}). According to lemma B.3 the function $q_l$, which is part of the solution $(p_l, q_l)$ of the boundary value problem (\ref{2.17}), (\ref{2.18}), then satisfies the above bound as well. Theorem 2.1 thus provides a weak solution $u \in \HC_0 (Q_R)$ of eq.\ (\ref{2.26a}) corresponding -- in cylindrical coordinates -- to a solution $u \in H^1_{as} (A_R)$ of eq.\ (\ref{2.22}). The conjugate variable $\wp \in H^1 (A_R)$ then is obtained from system (\ref{2.19}) as in \cite{KR22} with eq.\ (\ref{2.22}) representing the associated integrability condition. Note that a constant $\wp_0$ is here still free. With $p= \wp + p_l$, $q=u + q_l$ system (\ref{2.19}) is equivalent to  (\ref{2.13}), which in turn takes by $g = (p + i q)/2$ the complex form (\ref{2.12}). Finally, cancelling the substitution (\ref{2.10}), we obtain a weak solution $f$ of 
eq.\ (\ref{2.6}), equivalent -- in cylindrical coordinates -- to a solution $\HH$ of eq.\ (\ref{2.1}). As in \cite{KR22} higher interior regularity is most easily seen by changing once more coordinates. Using Cartesian coordinates in $\real^3$ each Cartesian component of $\HH$ turns out to be weakly harmonic, hence harmonic, and thus a smooth function in the spherical shell $B_R\setminus B_1$.

Different boundary conditions make the essential difference between the direction and the intensity problem. The crucial boundary condition is here (\ref{2.20}), which we have exploited in the form (\ref{2.24}), resulting in the weak formulation (\ref{2.26}), or (\ref{2.28}) when using polar coordinates. The following proposition recovers a weak version of the $\wp$-boundary condition from (\ref{2.28}).
\begin{prop}
Let $u \in \HC_0 (Q_R)$ be a solution of eq.\ (\ref{2.28}) for all $\psi \in \CC^1 (\overline{Q_R})$ and let $\wp \in H^1 (A_R)$ be determined by system (\ref{2.19}), then 
$$
{\rm trace }\; \wp\big|_{\pa A_R} = const.
$$
\end{prop}
\noi\textsc{Proof:} Let us start with writing eq.\ (\ref{2.28}) in the form
\be \label{6.1}
\ba{l}
\disp \int_1^R \int_{-\pi}^\pi \bigg[\bigg( \pa_r u + a_r [u - \Om]\ \frac{u - \Om}{ r \sin \vi}\bigg) \pa_r  \psi \\[2ex]
\disp \qquad\qquad \quad 
+ \bigg( \frac{1}{r} \pa_\vi u + a_\vi [u - \Om]\ \frac{u - \Om}{ r \sin \vi}\bigg) \frac{1}{r} \pa_\vi  \psi \bigg] \rd \vi \, r \rd r\, ,
\ea
\ee
where we made use of the antisymmetric continuations of the variables $u$, $\Om$, and $\psi$.\footnote{ Recall that only in section 2 variables written in polar coordinates are marked by a tilde. In this section we make free use of polar and cylindrical coordinates without further indication.} On the other hand system (\ref{2.19}) takes in polar coordinates and by (\ref{2.23}) the form 
\begin{equation}\label{6.2}
\left.
\begin{array}{rl}
\disp \pa_r \wp &= \disp \frac{1}{r} \pa_\vi u + a_\vi [u -\Om]\,
\frac{u - \Om}{r \sin \vi}\, , \\[2ex]
\disp -\frac{1}{r} \pa_\vi \wp &= \disp \pa_r u + a_r [u - \Om]\, \frac{u -\Om}{r \sin \vi}\, .  
\end{array} \right\}
\end{equation}
Considering the inner boundary component $S_1$ first, for the test function $\psi$ we make the ansatz:
\be \label{6.3}
\psi (r, \vi) =\psi_\eps (r)\, \psi_0 (\vi)\; , \quad 0< \eps < R\, ,
\ee
where $\psi_0$ is a $C^1$-function on $S_1$ and $\psi_\eps$ is given by 
\be \label{6.4}
\psi_\eps := 
\left\{ 
\ba{cl}
1 - (r -1)/\eps &\mbox{ for } \; 1 < r < 1+\eps\, , \\[1ex]
0 & \mbox{ for }\; 1 + \eps  \leq r < R\, .
\ea
\right.
\ee
We then have $\psi \in H^1 (A_R)$, which is admissible in (\ref{6.1}). Inserting (\ref{6.2})$_2$ into the first parantheses in (\ref{6.1}), abbreviating the second parantheses by $A[u, \Om]$, and using (\ref{6.3}) we obtain after integration by parts:
$$ 
\int_1^R \int_{-\pi}^\pi \big(\wp\,  \pa_r \psi_\eps \, \pa_\vi \psi_0 + A[u, \Om] \, \psi_\eps\, \pa_\vi \psi_0\big) \rd \vi \, \rd r = 0\, .
$$
Using (\ref{6.4}) in the first term and estimating the second by Cauchy-Schwarz we further obtain:
\be \label{6.5}
\ba{l}
\disp \bigg|\frac{1}{\eps} \int_1^{1 + \eps} \bigg(\int_{-\pi}^{\pi} \wp\, \pa_\vi \psi_0 \, \rd \vi \bigg) \rd r \bigg| \leq 
\int_1^{1 + \eps} \int_{-\pi}^{\pi} \Big| A\,\frac{1}{r} \pa_\vi \psi_0 \Big| \rd \vi \, r\rd r \\[3ex]
\disp \qquad \leq \bigg(\int_1^{1 + \eps} \int_{-\pi}^\pi A^2\, \rd \vi \, r \rd r \bigg)^{1/2} \bigg(\int_1^{1 + \eps} \int_{-\pi}^\pi \Big(\frac{1}{r} \pa_\vi \psi_0 \Big)^2 \rd \vi \, r \rd r \bigg)^{1/2} \\[2ex]
\disp \qquad \qquad \leq \sqrt{2}\, \bigg( \int_{Q_R} A^2\, \rd \mu \bigg)^{1/2} \sqrt{\eps}\, \bigg( \int_{-\pi}^\pi \big( \pa_\vi \psi_0\big)^2\, \rd \vi \bigg)^{1/2} .
\ea
\ee
The first factor on the right-hand side is bounded by (\ref{A.2}), (\ref{B.5}), and (\ref{B.18}):
$$
\ba{l}
\disp \bigg(\int_{Q_R} A^2\, \rd \mu \bigg)^{1/2} \leq \Big\|\frac{1}{r}\pa_\vi u + a_\vi [u - \Om]\, \frac{u - \Om}{r \sin \vi} \Big\|_{(0,0)} \\[2ex]
\disp \qquad \qquad \qquad \qquad \leq \| \na u\|_{(0,0)} + \sqrt{2}\, \bigg(\Big\| \frac{u}{r \sin \vi}\Big\|_{(0,0)} + \Big\| \frac{\Om}{r \sin \vi}\Big\|_{(0,0)} \bigg) \\[3ex]
\disp \qquad \qquad \qquad \qquad \leq \big( 1 + \sqrt{2}\, \pi\big)\, \|\na u\|_{(0,0)} + \sqrt{\pi}\, K\, R\, ,
\ea
$$
and the second factor represents a normalization of $\psi_0$ assumed to be fixed:
\be \label{6.6}
\int_{-\pi}^\pi \big(\pa_\vi \psi_0\big)^2 \, \rd \vi = 1\, .
\ee
Letting $\eps \ra 0$ in (\ref{6.5}) then yields
$$
\int_{-\pi}^{\pi} {\rm trace }\, \wp\big|_{r = 1} \, \pa_\vi \psi_0 \, \rd \vi = 0
$$
for all antisymmetric $\psi_0 \in C^1 (S_1)$ with normalization (\ref{6.6}), which in turn means 
$$
{\rm trace }\, \wp \big|_{S_1} = const .
$$
With an obviously modified sequence of test functions we obtain 
$ {\rm trace }\, \wp |_{S_R} = const$ as well.
\qed 

Note that in general these constants may be different from zero and different from one another; however, by adding a suitable constant in (\ref{2.19}) we can achieve that
\be \label{6.7}
{\rm trace }\, \wp \big|_{S_1} = 0\; , \quad {\rm trace }\, \wp \big|_{S_R} =: \wp_0\, . 
\ee
\vspace{.2ex}
\begin{figure}
\begin{center}
\includegraphics[width=0.5\textwidth]{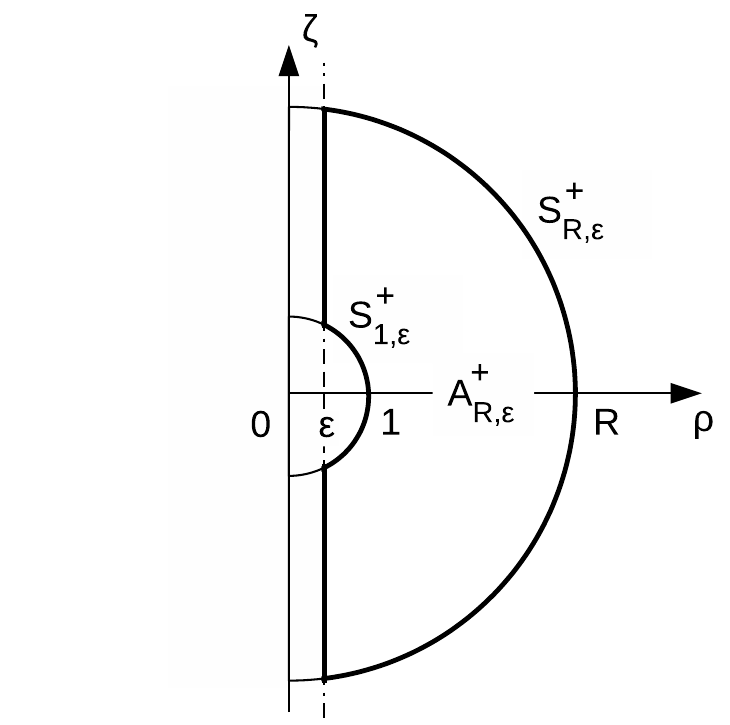}
\caption{$A_{R, \eps}^+$ and the curved boundary components $S_{1, \eps}^+$ and $S_{R, \eps}^+$.}
\end{center}
\end{figure}

Higher boundary regularity then follows as in \cite{KR22}, however, with interchanged variables $\wp$ and $u$:
the second-order equation for the variable $\wp$,
$$
\De \wp = \pa_\ze \Big(\frac{1}{\rho} \sin (u - \Om) \Big) + \pa_\rho \Big(\frac{1}{\rho} \big(\cos (u - \Om ) \big) - 1\Big) =: \pa_\ze F_\ze + \pa_\rho F_\rho \, ,
$$
follows from (\ref{2.19}). On the regularized domain $A^+_{R, \eps}$ (see Fig.\ 4) the vector field $\FF$ is bounded by
$$
\big| F_\ze \big| < \frac{1}{\eps}\; ,\quad \big|F_\rho \big| < \frac{2}{\eps}\, ,
$$
and theorem 14.1 in \cite{LU68} applied to $\wp \in H^1 (A_R)$ then yields $\wp \in C^{0,\al} (\overline{A^+_{R, \eps}} )$ for some $0< \al \leq 1$ (depending on $\eps$, $R$, and $\|\wp\|_{(0,0)}$). This implies, in particular, continuous assumption of the boundary values (\ref{6.7}) on $S^+_{1, \eps} \cup S^+_{R,\eps}$ for any $\eps > 0$, hence on $S^+_{1} \cup S^+_{R}$, and by symmetric continuation on $\pa A_R = S_1 \cup S_R$.

Finally, assume $I \in C^{0,\al} (S_1)$, $\Ih \in C^{0,\al} (S_R)$, then $p|_{\pa A_R} \in C^{0, \al} (\pa A_R)$ by (\ref{2.15}), (\ref{2.16}), and for the solution $p_l$ (and $q_l$) of the boundary value problem (\ref{2.17}), (\ref{2.18}) holds $p_l \in C^{0,\al} (\overline{A_R})$. Therefore, $p =\wp + p_l$ assumes continuously its boundary values, satisfies (\ref{2.15}) at the inner boundary, and (\ref{2.16}) in the form 
$$
p\big|_{S_R} = p_l\big|_{S_R} + \wp\big|_{S_R} 
= 2\, \Big( \ln \Ih - \sum_{n=1}^{\ro- \roh} \ln |z - z_n|\Big|_{S_R} + p_{l,0} \Big) + \wp_{0} 
$$
at the outer boundary. The constant $I_0$ in problem $P_{I, \Ih} (A_R)$ is thus determined by the relation 
$$
\ln I_0 + \ro \ln R = p_0 = p_{l,0} + \wp_0 /2\, .
$$

This argument does not apply to $q$ since in the weak formulation this variable is free at the boundary. Instead, we make use of an argument based on Privaloff's theorem \cite{BJS64} and a representation theorem due to Bers and Nirenberg \cite{BN54}. The former theorem guarantees uniform H\"older continuity of both components of a complex analytic function in a disc if only one component takes continuously boundary values that are in some H\"older class at the boundary; the other theorem represents solutions of system (\ref{2.13}) as sum of a complex analytic function and a uniformly H\"older continuous function. Details of the argument can be found in \cite{KR22} (where the variables $p$ and $q$ have to be interchanged). By (\ref{2.5}), (\ref{2.10}) this establishes $\HH \in C(\overline{A_R})$.

Concerning uniqueness recall that fixed outer rotation number $\roh$, fixed set of zeroes, and fixed intensity functions $I$ and $\Ih$ mean a fixed data function $\Om$. Hence uniqueness of the weak solution as stated in theorem 2.1 implies uniqueness as claimed in theorem 2.3. Of course, a sign always remains free in problem $P_{I, \Ih} (A_R)$. This concludes the proof of theorem 2.3.

With obvious modifications all these arguments work in the unbounded case as well: a unique weak solution $u \in \HC_0^{loc} (A_\infty)$ is provided by theorem 2.2, lemma B.1 (where $\roh$ is replaced by $\dt -1$), and lemma B.2; the conjugate variable $\wp$ is again determined by system (\ref{2.19}); the inner boundary condition (\ref{2.15}) follows as in the bounded case, whereas the outer boundary condition (\ref{2.16}) is no longer present; the regularity arguments, local in character, can be taken over without any change. This proves theorem 2.4.
\sect{A related intensity problem}
In \cite{DDO06} the authors introduce ``an oblique boundary value problem related to the Backus problem in geodesy''. Starting point is the intensity problem (\ref{1.1}) reformulated in terms of a potential $u$ for the gravitational field $\HH$ in the exterior space $E$: 
\begin{equation}\label{7.1}
\left.
\begin{array}{cl}
\De u = 0 & \quad \mbox{ in }E , \\[1ex]
|\na u | = I & \quad \mbox{ on } S^2 , \\[1ex]
u(\xx) \ra 0 & \quad \mbox{ for } |\xx|\ra \infty\, . 
\end{array} \right\}
\end{equation}
By the Kelvin transformation 
$$\xx \ra \xx^* := \frac{\xx}{|\xx|^2}\; , \qquad u(\xx) \ra v(\xx^*) := \frac{1}{|\xx^*|}\, u\Big(\frac{\xx^*}{|\xx^*|^2}\Big)
$$
the exterior problem (\ref{7.1}) can equivalently be formulated as interior problem in the unit ball $B_1$:
\begin{equation}\label{7.2}
\left.
\begin{array}{cl}
\De v = 0 & \quad \mbox{ in }B_1 , \\[1ex]
(v + \nn\cdot \na v)^2 + |\na_t v |^2 = I^2 & \quad\mbox{ on } \pa B_1 \, ,  
\end{array} \right\}
\end{equation}
where $\nn := \xx/|\xx|$ denotes the exterior unit normal at $\pa B_1$ and $\na_t := \na  - \nn\, \nn \cdot \na$ the tangential gradient at $\pa B_1$. 

Under the condition 
\be \label{7.3}
v + \nn\cdot \na v \geq 0  \quad \mbox{ on } \; \pa B_1\, ,
\ee
which is equivalent to
\be \label{7.4}
\nn\cdot \na u \leq 0  \quad \mbox{ on } \; S^2
\ee
in the exterior problem, the boundary condition (\ref{7.2})$_2$ can equivalently be expressed by  
\be \label{7.5}
v + \nn\cdot \na v = \big( I^2 - |\na_t v|^2\big)^{1/2} ,
\ee
which implies, in particular,
\be \label{7.5a}
|\na_t v| \leq I \quad \mbox{ on }\; \pa B_1\, .
\ee
Without condition (\ref{7.5a}) the modified boundary condition\footnote{$(f)_+$ denotes the nonnegative part $\max_{\pa B_1} \{f, 0\}$ of a function $f: \pa B_1 \ra \real$.} 
$$
v + \nn\cdot \na v = \Big(\big( I^2 - |\na_t v|^2\big)_+ \Big)^{1/2} 
$$
still makes sense, and the boundary value problem 
\begin{equation}\label{7.6}
\left.
\begin{array}{cl}
\De v = 0 & \quad \mbox{ in }B_1 , \\%[1ex]
v + \nn\cdot \na v = \Big(\big(I^2 - |\na_t v|\big)_+\Big)^{1/2} & \quad\mbox{ on } \pa B_1  
\end{array} \right\}
\end{equation}
is of oblique type with corresponding solution theory (see, e.g., \cite{L13}). The authors first solve problem (\ref{7.6}) with a regularized boundary condition, viz.,
$$
v + \nn\cdot \na v = \Big(\big( I^2 - |\na_t v|^2 \big)_+ + \eps^2 \Big)^{1/2} , \quad \eps > 0\, ,
$$
and find for $I\geq 0$, satisfying some Lipschitz condition on $\pa B_1$, a unique solution $v_\eps$ of class $C^2 (B_1) \cap C^{1, \beta} (\overline{B_1})$ for some $\beta > 0$. 
Subsequently, letting $\eps \ra 0$, they prove a subsequence to converge to some function $v \in C^\al (\overline{B_1})$ for some $\al  \in (0\, , 1]$, which is, moreover, a viscosity solution of problem (\ref{7.6}) (see \cite{B93}). Concerning regularity and connection to the intensity problem they formulate the following conjecture:

\vspace{.3em}\noi
{\bf Conjecture (3.1 in \cite{DDO06}).} {\em Problem (\ref{7.6}) has a unique solution $v \in C^2 (B_1) \cap C^{1} (\overline{B_1})$, which satisfies, moreover, condition (\ref{7.5a}), i.e., $v$ is also solution of the intensity problem (\ref{7.2}) and satisfies condition (\ref{7.3}).}
\vspace{.3em}

\noi We refute this conjecture by an axisymmetric counterexample that is the {\em unique} solution of the intensity problem but violates condition (\ref{7.3}). This counterexample, formulated in the exterior space $E$, is based on (the uniqueness part of) Theorem 2.4 and (the existence part of) Theorem 2.4 in ref.\ \cite{KR22}. The latter theorem answers the (signed) direction problem $P_\DD (A_\infty)$, which is here reproduced for the reader's convenience:

\vspace{.3em}

\noi\textbf{Problem $P_{\DD}(A_\infty)$.} {\it Let $\DD \in C(S_1, \real^2)$ be a symmetric direction field and $\delta
\in \nat \setminus \{1\}$. Determine all symmetric solutions $\HH \in C^1 (A_\infty) \cap C(\overline{A_\infty})$ of system (\ref{2.1}) with decay order $\de$ and boundary condition
\be \label{7.7}
\exists\, a \in C (S_1, \real_+) : \HH \big|_{S_1} = a\, \DD\, . 
\ee
}

\noi{\bf Theorem (2.4 in \cite{KR22}).} {\em
Let $\DD$ be a H\"older continuous, symmetric direction field with rotation number $\ro \in \nat $ and satisfying condition 
\be \label{7.8}
D_\rho (\vi) = O (\vi) \quad \mbox{for }\, \vi \ra 0\, ,\qquad D_\rho (\vi) = O (\pi - \vi) \quad  \mbox{for }\, \vi \nearrow \pi
\ee
at the symmetry axis. Let, furthermore, $\dt \in \nat \setminus \{1\}$, $\dt \leq \ro +1$, and $\{z_1, \ldots,$\linebreak
$z_{\ro - \dt +1} \}$ be a symmetric set of points in $A_\infty$. Then, problem $P_{\DD} (A_\infty)$ has a (up to a positive constant factor)
unique solution $\HH$ with exact decay order $\dt$ vanishing at  $z_1, \ldots,z_{\ro - \dt + 1}$ and nowhere else.
 }
 \vspace{.3em}
 
\noi The formulation of problem and theorem here differs from that in \cite{KR22} by the inclusion of the case $\ro = 1$, $\dt = 2$, which refers to the gravitational case, whereas in \cite{KR22} exclusively the magnetic problem had been considered. The theorem, however, applies to this case as well as can easily be checked. 

\vspace{.3em}

\noi{\bf Counterexample (in the exterior space).} Let $\DD$ be a H\"older continuous, symmetric direction field satisfying condition (\ref{7.8}) with rotation number $\ro = 1$ that {\em violates} condition (\ref{7.4}), i.e.\ $\nn \cdot \DD$ is changing sign along $S_1$. According to the theorem above there is a unique solution $\HH_0 \in C^1 (A_\infty) \cap C(\overline{A_\infty})$ of Problem $P_\DD (A_\infty)$ with exact decay order $\dt = 2$, no zeroes in $A_\infty$, and boundary condition
$$
\HH_0 \big|_{S_1} = a\, \DD\; \mbox{ for some } \; a \in C (S_1, \real_+)\, .
$$
Obviously, $\HH_0$ is as well a solution of the intensity problem $P_I (A_\infty)$ with $I:= |\HH_0|\big|_{S_1}$ and rotation number $\ro =1$. The condition $\ro =1$ implies by relation (\ref{2.9}) the decay order $\dt =2$ and no zeroes in $A_\infty$ for {\em any} solution $\HH$ of $P_I (A_\infty)$. Hence, by theorem 2.4, $\HH_0$ is the only solution that fits to $I$ at $S_1$ and this solution violates by construction condition (\ref{7.4}). Solutions with $\ro >1$ violate condition (\ref{7.4}) in any case. Note that regularity conditions on $I$ (such as condition (\ref{2.33})) do not matter when only appealing to the uniqueness part of theorem 2.4.
\vspace{.2em}

\noi{\em Remark:} This example prescribes a direction field, not an intensity function. Therefore, the counterexample only demonstrates that there are intensity functions for which problems (\ref{7.2}) and (\ref{7.6}) differ. To read off a given intensity function if it allows ``gravitational solutions'', i.e., solutions satisfying condition (\ref{7.4}), is an interesting open problem. On the other hand theorem 2.4 guarantees always the existence of ``nongravitational solutions'' of the intensity problem for a given intensity function (just by prescribing zeroes in $A_\infty$).  
\setcounter{equation}{0}
\setcounter{section}{0}
\renewcommand{\theequation}{\thesection\arabic{equation}}
\renewcommand{\thesection}{\Alph{section}.}

\sect{Hardy-type inequalities}
Weighted Hardy-type inequalities in $L^2$ and, more generally, in $L^p$-spaces are the essential tool to get to grips the singularity at the symmetry axis. We collect here those inequalities from \cite{KR22}, which are needed in the present context, adapted to polar coordinates and free boundary conditions. We prove only those parts that have changed and rely otherwise on the proofs given in \cite{KR22}. 

Using polar coordinates $(r, \vi)$ on the rectangle $Q_R = (1, R)\ti (0, \pi)$ the following notation is convenient to describe weighted $L^p$-norms: 
$$
\| \cdot \|_{p, (\ga, \eta)} := \bigg( \int_{Q_R} | \cdot |^p\, \rd \mu_\ga^\eta \bigg)^{1/p} , \qquad \rd \mu_\ga^\eta  := \vih^{-\ga} \rd \vi\, r^{1-\eta} \rd r
$$
with $p \geq 1$, $\ga  \in \real$, $\eta  \in \real$, and $\vih$ denoting the ``tent function'' (see Fig.\ 5)
\begin{figure}
\begin{center}
\includegraphics[width=0.6\textwidth]{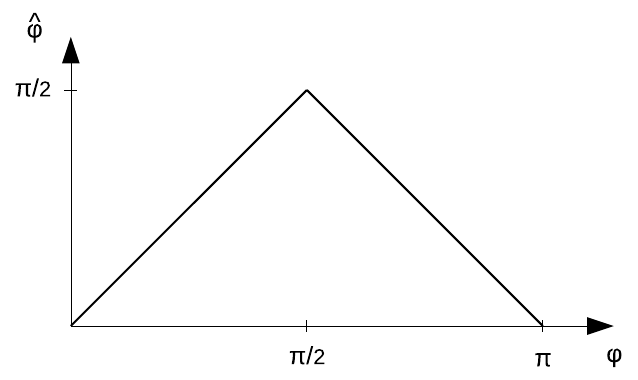}
%\resizebox{1\textwidth}{!}{\input{tent.eps}}
%\input{A4_loglike.tex}
\caption{Graph of the tent function $\vih : [0\, , \pi] \ra [0\, , \pi/2]$.}
\end{center}
\end{figure}
$$
\vih:= \left\{ 
\ba{cl}
\vi & \mbox{ for }\; 0< \vi < \pi/2\, , \\[1ex]
\pi - \vi & \mbox{ for }\; \pi/2 \leq \vi < \pi\, . 
\ea \right. 
$$
For $\ga = \eta = 0$ we use the simplified notation $\rd \mu_0^0 =: \rd \mu$, for $p=2$ we use: $\| \cdot \|_{2, (\ga, \eta)} =: \| \cdot \|_{(\ga, \eta)}$, and for $p= \infty$:  $\| \cdot \|_{\infty, (0, 0)} =: \| \cdot \|_\infty$. With the weighted gradient
$$
\na_{\!c} := \eee_r c\, \pa_r + \eee_\vi  \frac{1}{r} \pa_\vi\, ,\quad c > 0
$$
and the test function space 
$$
\CC^1 (\QRb) := \big\{\psi \in C^1(\QRb) : \psi (\cdot, 0) = \psi (\cdot, \pi) = 0 \big\}
$$
we define
$$
\HC_{p, \ga} := \HC_{p, \ga} (Q_R) := clos\big(\CC^1  (\QRb)\, ,\, \|\na_{\!c}\cdot \|_{p, (\ga, \eta)}\big)
$$
with the abbreviation $ \HC_{2,\ga}=: \HC_\ga$. Note that on the bounded domain $Q_R$ norms with different weights $\eta$ or $c> 0$ are equivalent. With respect to $p$ we have the well-known inclusion
$$
\HC_{p,\ga} (Q_R) \subset \HC_{q,\ga} (Q_R) \qquad \mbox{for } \; p \geq q\, ,
$$
and with respect to $\ga$ we have: 
$$
\HC_{p,\ga} (Q_R) \subset \HC_{p,\beta} (Q_R) \qquad \mbox{for } \; \ga \geq \beta\, .
$$
On $Q_\infty$ we define
$$
\HC^{loc}_{p, \ga} (Q_\infty) := \big\{ f: Q_\infty \ra \real\, : f|_{Q_R} \in \HC_{p,\ga} (Q_R) \mbox{ for any } R>1 \big\}\, .
$$
In $\HC_{p,\ga}$ the following Hardy-type inequalities hold:
\begin{lem}
 Let $f \in \HC_{p,\ga} (Q_R)$ or $\in \HC^{loc}_{p,\ga} (Q_\infty)$ with $\ga \neq 1-p$ and $p \geq 1$. Then, the following inequalities hold:
 \be \label{A.1}
 \Big\|\frac{f}{\vih}\Big\|_{p,(\ga, \eta)} \leq \frac{p}{|p + \ga -1|}\, \|\pa_\vi f\|_{p,(\ga, \eta)}
 \ee
 and, especially for $p=2$,
 \be \label{A.2}
 \Big\|\frac{f}{\vih}\Big\|_{(\ga, \eta)} \leq \frac{2}{|1 + \ga |}\, \|\pa_\vi f\|_{(\ga, \eta)}\, .
 \ee 
\end{lem}
\textsc{Proof:}
We split the domain $Q_R$ into subdomains $(1,R)\ti (0,\pi/2)$ and $(1,R)\ti (\pi/2, \pi)$ and map the latter subdomain by the substitution $\vi \mapsto \pi - \vi$ onto the former. On $(1,r)\ti (0,\pi/2)$ the proof proceeds as for proposition 3.1 in \cite{KR22}.
Note that no boundary condition is necessary at the right side wall $\{\vi = \pi/2\}$.
\qed

\vspace{1ex}

\noi By (\ref{A.1}) we obtain immediately the Poincar\' e-type inequalities:
\be \label{A.3}
\|f\|_{p,(\ga, \eta)} \leq \frac{\pi}{2}\, \frac{p}{|p + \ga -1|}\, \|\pa_\vi f\|_{p,(\ga, \eta)} \leq R\, \frac{\pi}{2}\, \frac{p}{|p + \ga -1|}\, \|\na_{\!c} f\|_{p,(\ga, \eta)}\, .
\ee
Let $H^1  (Q_R)$ be the usual Sobolev space $\mbox{clos}\big(C^1 (Q_R) , \|\cdot \|_{H^1} \big)$ with $H^1$-norm
$$
\|f \|^2_{H^1} = \int_{Q_R} \big(|f|^2 + |\na\, f|^2  \big)\, \rd \vi\, r \rd r \, ,
$$
then by (\ref{A.3}) we have the norm equivalence $\|\cdot \|_{H^1} \sim \|\na_{\!c} \cdot \|_{(0, \eta)}$ and by definition of $\HC_0$ the alternative characterization
$$
\HC_0 (Q_R) = \big\{ f \in H^1 (Q_R) : {\rm trace} f|_{\vi = 0} ={\rm trace} f|_{\vi = \pi} = 0 \big\}\, .
$$
By (\ref{A.2}) follows for $\ga \neq -1$ the further equivalence
\be \label{A.4}
\|\na_{\!c} \cdot \|_{(\ga,\eta)} \sim \| \na_{\!c} (\vih^{-\ga/2} \, \cdot\, )\|_{(0,\eta)}
\ee
and hence
\be \label{A.5}
\ba{rl}
\HC_\ga  (Q_R) & \disp = {\rm clos } \big( \CC^1 (\QRb)\, ,\, \|\na_{\!c}(\vih^{-\ga/2}\, \cdot \, )\|_{(0,\eta)} \big) \\[2ex]
&\disp  = \big\{ f : \vih^{-\ga/2} f \in \HC_0 (Q_R) \big\} \, .
\ea
\ee

Crucial for our solution procedure are two inequalities of type 
\be \label{A.6}
\Big\|\frac{f}{\vih^\de}\Big\|_{p,(0, \eta)} \leq C \, \|\na_{\!c} \, f \|_{(\beta, \eta)}
\ee
for ``large'' and ``small'' $p$ and as large as possible $\de$ for given $\beta$.

\vspace{.3em}
\noi{\bf Lemma A.2 (large-p-case).}
{\em
Let $c> 0$, $\eta > 0$, $p\geq 2$, $\beta > -2/p$, $\de < \beta/2 + 1/p$, $R \geq 2$, and  $f \in \HC_{\beta} (Q_R)$, then inequality (\ref{A.6}) holds with a constant $C = C_>$ depending on $c$, $\eta$, $p$, $\beta$, $\de$, and $R$.}
\vspace{.3em}

\noi\textsc{Proof:} Starting point are two one-dimensional inequalities, viz.,
\be \label{A.7}
\sup_{0 < \vi < \pi} \big\{|h(\vi)|\, \vih^{-\ga} \big\} \leq \int_0^\pi |\pa_\vi h |\, \vih^{-\ga} \, \rd \vi 
\ee
and
\be \label{A.8}
\sup_{1 < r < R} \big\{|g(r)|\big\} \leq \int_1^R |\pa_r g | \, \rd r + \frac{1}{R-1} \int_1^R |g|\, \rd r\, .
\ee
Inequality (\ref{A.7}) holds for functions $h \in H^1 ((0, \pi))$ with
$h(0) = h(\pi) = 0$ and corresponds to inequality (3.18) in \cite{KR22}. The mean-value version (\ref{A.8}) differs from the corresponding inequality in \cite{KR22} since no boundary conditions are available in the radial direction. It holds for $g \in H^1((1,R))$ and follows by 
$$
\ba{rl}
g(r) & \disp = \frac{1}{R -1} \int_1^R \big( g(r) - g (s) \big) \rd s + \frac{1}{R-1} \int_1^R g(s)\, \rd s \\[3ex]
&  \disp \qquad = \frac{1}{R -1} \int_1^R \int_s^r \pa_t g(t)\, \rd t \, \rd s + \frac{1}{R-1} \int_1^R g(s)\, \rd s \\[3ex]
& \disp \qquad \qquad \leq \int_1^R |\pa_t g(t) |\, \rd t  + \frac{1}{R-1} \int_1^R |g(s)|\,  \rd s\, .
\ea
$$
Let now $f \in \HC_{1, \ga} (Q_R)$. By (\ref{A.7}), (\ref{A.8}) one obtains 
$$
\ba{l}
\disp \int_{Q_R} f^2\, \vih^{-\ga}\, \rd \mu_\ga^\eta  = \int_1^R \int_{0}^\pi \frac{f}{r} \, \vih^{-\ga} \ti r f \, \vih^{-\ga} \rd \vi \, r^{1- \eta} \rd r \\[3ex]
\disp \qquad \leq \int_{1}^R \frac{1}{r}\, \sup_\vi \big\{|f (r,\vi)|\, \vih^{-\ga} \big\} r^{1-\eta} \rd r \times \int_0^\pi \sup_r \big\{|r f(r,\vi)|\big\}\, \vih^{-\ga} \rd \vi \\[3ex]
\disp \qquad \qquad \leq \int_{1}^R \frac{1}{r} \int_0^\pi |\pa_\vi f |\, \vih^{-\ga} \rd \vi\, r^{1- \eta} \rd r \times \bigg[ \int_0^\pi \int_{1}^R \big|\pa_r (r f)\big|\rd r\, \vih^{-\ga} \rd \vi \\[2ex]
\disp \hfill + \frac{1}{R - 1} \int_0 ^\pi \int_1^R r |f| \rd r\, \vih^{-\ga} \rd \vi \bigg]
\ea
$$
The first factor can immediately be estimated by
$$
\int_{Q_R} |\na_{\!c} f|\, \rd \mu_\ga^\eta\, ,
$$
whereas for the second we apply once more (\ref{A.7}) to obtain
$$
\ba{l}
\disp \int_1^R \int_0^\pi |\pa_r f|\, \vih^{-\ga} \rd \vi \, r \rd r + \int_1^R \int_0^\pi \Big(\frac{1}{R-1} + \frac{1}{r}\Big) |f|\, \vih^{-\ga} \rd \vi\, r \rd r \\[3ex]
\disp \qquad \leq \frac{R^\eta}{c} \int_{Q_R} |\na_{\!c} f|\, \rd \mu_\ga^\eta + \pi \int_1^R \int_0^\pi \Big(\frac{r}{R-1} + 1\Big) \frac{1}{r}\,  |\pa_\vi f|\, \vih^{-\ga} \rd \vi\, r \rd r \\[2ex]
\disp \qquad \qquad \leq (3 \pi + c^{-1})\, R^\eta \int_{Q_R} |\na_{\!c} f|\, \rd \mu_\ga^\eta\, .
\ea
$$
In summary we have 
$$
\int_{Q_R} f^2\, \vih^{-\ga}\, \rd \mu_\ga^\eta \leq (3 \pi + c^{-1})\, R^\eta \bigg(\int_{Q_R} |\na_{\!c} f|\, \rd \mu_\ga^\eta\bigg)^2 ,
$$
which corresponds to eq.\ (3.23) in \cite{KR22}. The rest of the proof proceeds as in \cite{KR22}.
\qed

%\vspace{1ex}

\vspace{.3em}
\noi{\bf Lemma A.3 (small-p-case).} {\em
Let $c> 0$, $\eta > 0$, $\beta > -1$, $2 \leq p \leq 4$, $\de = \beta/2 + (6-p)/(2 p)$, $R \geq 2$, and  $f \in \HC_{\beta} (Q_R)$, then inequality (\ref{A.6}) holds with a constant $C = C_<$ depending on $c$, $\eta$, $\beta$, $p$, and $R$.}
\vspace{.3em}

\noi\textsc{Proof:} Starting point are again one-dimensional inequalities, viz.,
\be \label{A.9}
\sup_{0 < \vi < \pi} \big\{|h(\vi)|^2 \, \vih^{-(1 + \beta)} \big\} \leq (1 + \beta)^{-1}\int_0^\pi |\pa_\vi h |^2\, \vih^{-\beta} \, \rd \vi\, , 
\ee
which corresponds to inequality (3.17) in \cite{KR22}, and
\be \label{A.10}
\sup_{1 < r < R} \big\{|r g(r)|^2\big\} \leq 2 \bigg(\int_1^R |\pa_r g | \, r \rd r \bigg)^2 + 18 \bigg( \int_1^R |g|\, \rd r\bigg)^2 ,
\ee
which follows from (\ref{A.8}). An analogous calculation as in the preceeding lemma yields by (\ref{A.9}), (\ref{A.10}):
$$
\ba{l}
\disp \int_{Q_R} f^4\, \vih^{-(1 + \beta)}\, \rd \mu_\beta^\eta  = \int_1^R \int_{0}^\pi \Big|\frac{f}{r}\Big|^2 \, \vih^{-(1 + \beta)} \ti |r f|^2 \, \vih^{-\beta} \rd \vi \, r^{1- \eta} \rd r \\[3ex]
\disp \qquad \leq \int_{1}^R \frac{1}{r^2}\, \sup_\vi \big\{|f (r,\vi)|^2 \, \vih^{-(1 + \beta)} \big\} r^{1-\eta} \rd r \times \int_0^\pi \sup_r \big\{|r f(r,\vi)|^2\big\}\, \vih^{-\beta} \rd \vi \\[3ex]
\disp \qquad \qquad \leq (1 + \beta)^{-1} \int_{1}^R \frac{1}{r^2} \int_0^\pi |\pa_\vi f |^2\, \vih^{-\beta} \rd \vi\, r^{1- \eta} \rd r 
\\[2ex]
\disp \qquad \qquad \quad \times \bigg[2 \int_0^\pi \bigg(\int_{1}^R |\pa_r f|\,\rd r\bigg)^2 \vih^{-\beta} \rd \vi + 18 \int_0 ^\pi \bigg(\int_1^R  |f|\, \rd r\bigg)^2 \vih^{-\beta} \rd \vi \bigg] .
\ea
$$
The first factor is estimated by 
$$
(1 + \beta)^{-1} \int_{Q_R} |\na_{\!c} f|^2\, \rd \mu_\beta^\eta\, ,
$$
the first term of the second factor using H\"olders's inequality is estimated by 
$$
(R^2 -1) \int_1^R \int_0^\pi | \pa_r f|^2\, \vih^{-\beta} \rd \vi \, r \rd r \leq \frac{R^{2 + \eta}}{c^2} \int_{Q_R} |\na_{\!c} f|^2\, \rd \mu_\beta^\eta\, ,
$$
and the second term using once more (\ref{A.7}) by 
$$
\ba{l}
\disp 18 \, \pi \bigg(\int_1^R \int_0^\pi |\pa_\phi f|\, \vih^{-\beta/2} \rd \vi\, \rd r\bigg)^2 \\[3ex]
\disp \quad\leq 9 \, \pi^2 (R^2 -1) \int_1^R \int_0^\pi \Big|\frac{1}{r} \pa_\phi f\Big|^2\, \, \vih^{-\beta} \rd \vi\, r \rd r \leq 9\, \pi^2 R^{2 + \eta} \int_{Q_R} |\na_{\!c} f|^2\, \rd \mu_\beta^\eta\, .
\ea
$$
In summary we have
$$
\int_{Q_R} f^4\, \vih^{-(1 + \beta)}\, \rd \mu_\beta^\eta \leq \big( 9\, \pi^2 + c^{- 2}\big) \frac{R^{2 + \eta}}{1 + \beta} \bigg(\int_{Q_R} |\na_{\!c} f|^2\, \rd \mu_\beta^\eta \bigg)^2 ,
$$
and interpolation with (\ref{A.2}) then yields as in \cite{KR22} the proposition.
\qed
\sect{Estimates of the data function $\Om$}
The data function $\Om= \Psi + q_l$ comprises the zero-positions angle $\Psi$ and the boundary value function $q_l$. The estimates of $\Psi$ have already been done in \cite{KR22} and we report here just the results in a slightly different notation, whereas the estimates of $q_l$ are new and are given with full proof.

Let us write the analytic function (\ref{2.11}) in the form
\be \label{B.1}
\ba{c}
\disp h(z) = z^{-\ro} \prod_{n =1}^N (z- z_n) (z - \zb_n) = z^{-\roh} \prod_{n =1}^N \Big(1- \frac{z_n}{z}\Big) \Big(1 - \frac{\zb_n}{z}\Big)\, ,\\[3ex]
\disp \ro - \roh = 2\, N\, ,\ N \in \nat_0
\ea
\ee
with points $\{z_1, \ldots , z_N\} \subset A_\infty^+$, which satisfy, in particular, $\Im z_n > 0$. Decomposing $\hb/h$ into real and imaginary parts and introducing polar coordinates $(r, \vi)$ one obtains
\be \label{B.2}
\frac{\overline{h(z)}}{h(z)} =: v(r, \vi) + i\, w(r , \vi)\, , \qquad (r, \vi) \in Q_\infty
\ee
with 
\be \label{B.3}
v= \frac{\hb^2 + h^2}{2 \,|h|^2}\, ,\qquad w = \frac{\hb^2 - h^2}{2 i\, |h|^2}\, , \qquad v^2 + w^2 = 1\, .
\ee
We then define 
\be \label{B.4}
\Psi (r, \vi) := \arctan \frac{w(r, \vi)}{v (r, \vi)}\, ,
\ee
where $\arctan$ means the principal branch of the inverse $\tan$ function with $\arctan 0 = 0$. This definition resolves the ambiguity in the choice of $\Psi$ and fixes the discontinuities of $\Psi$ at the zeroes of $v$. 
\begin{lem}
 Let $\Psi$ be as given by (\ref{B.1}) -- (\ref{B.4}). Then the following inequalities hold:
 \be \label{B.5}
 |\Psi ( \cdot , \vi) |\leq K\, \sin \vi \qquad \mbox{in } \ Q_\infty 
 \ee 
 and 
 \be \label{B.6}
 \int_{Q_\infty} \bigg(\frac{\Psi}{r \sin \vi}\bigg)^2\, \rd \mu_\beta^\eta \leq \Big(\frac{\pi}{2}\Big)^{1- \beta}\frac{2 K^2}{\eta (1 - \beta)}\, , \qquad 0 \leq \beta < 1\, , \quad \eta >0 
  \ee
 for some constant $K> 0$ that depends on the positions of the zeroes.
\end{lem}

The function $q_l$ is part of a conjugate pair $(p_l, q_l)$ satisfying the system (\ref{2.17}) in the annulus $A_R$ and with boundary conditions (\ref{2.18}) for the function $p_l$. The goal of the subsequent treatment of this well-known boundary value problem is to find suitable conditions on the boundary values of $p_l$ such that $q_l$ satisfies an estimate of type (\ref{B.5}).

The bounded and the unbounded cases have to be treated separately; the case 
$A_\infty$ turns out to be simpler and is treated first. Using complex notation $f:= p_l + i q_l$, $z = \ze + i \rho$ and switching from $A_\infty$ to the unit disk $D:= \{z : |z| < 1\}$, solutions $f$ of problem (\ref{2.17}), (\ref{2.18})$_a$ with continuous boundary data are well-known to be complex analytic in $D$ and continuous in $\overline{D}$. Cauchy-Schwarz' formula provides in this situation an explicit representation of $q_l$ in terms of the boundary values of $p_l$ (see e.g.\ \cite{B79}, p.\ 135):  
\be \label{B.7}
q_l (z) = \frac{1}{ 2 \pi} \int_0^{2\pi} \Im\, \bigg\{\frac{e^{i \theta} + z}{e^{i \theta} - z} \bigg\}\, p_l (e^{i\theta})\, \rd \theta + q_l (0)\, ,\quad z \in D\, .
\ee
The last term in (\ref{B.7}) vanishes by antisymmetry of $q_l$, so that a corresponding formula for the exterior domain $A_\infty$ may be obtained by the inversion $z \mapsto 1/\zb$. Switching to polar coordinates $(r, \vi)$ one obtains 
\be \label{B.7a}
\ba{l}
\disp \qt _l(r, \vi) = -\frac{r}{\pi} \int_0^{2\pi} \frac{\sin (\theta - \vi)}{1 - 2 r \cos (\theta - \vi) + r^2} \; \pt_l (\theta)\, \rd \theta \\[3ex]
\disp\qquad  \; = \frac{r}{\pi} \int_0^{2 \pi}\frac{\sin (\theta + \vi)}{1 - 2 r \cos (\theta + \vi) + r^2} \; \pt_l (-\theta)\, \rd \theta\, ,\quad (r, \vi) \in (1, \infty)\times (-\pi , \pi]\, ,
\ea
\ee
where we have set $z = r e^{i \vi}$, $q_l (z) = \qt_l (r, \vi)$, and $p_l (e^{i \theta} ) =\pt_l (\theta)$. The change of variable  
\be \label{B.7b}
s:= \frac{2 r}{1 + r^2}
\ee
maps $A_\infty$ back onto $D$ while simplifying the denominator under the integral: 
\be \label{B.8}
\qb _l(s, \vi) = \frac{s}{2 \pi} \int_0^{2\pi} \frac{\sin (\theta + \vi)}{1 - s \cos (\theta + \vi)} \; \pt_l (-\theta)\, \rd \theta \, ,\quad (s, \vi) \in (0, 1)\times (-\pi , \pi]
\ee
with $\qb_l (s, \cdot) = \qt_l (r, \cdot)$. Finally, exploiting the 
symmetries of the problem, viz.,
$$
\pt_l(-\theta) = \pt_l (\theta)\, \quad \mbox{and} \quad \qb_l (\cdot, -\vi) = - \qb_l (\cdot, \vi)\, ,
$$
the integral may be restricted onto the interval $[0,\pi]$:
\be \label{B.9}
\ba{l}
\disp \qb _l(s, \vi) = \frac{s}{2 \pi} \int_0^{\pi} \bigg(\frac{\sin (\theta + \vi)}{1 -  s \cos (\theta + \vi)} - \frac{\sin (\theta - \vi)}{1 -  s \cos (\theta - \vi)}\bigg)\, \pt_l (\theta)\, \rd \theta \\[3ex]
\disp \qquad \quad \; =\frac{s}{\pi} \sin \vi \int_0^\pi \frac{\cos \theta - s \cos \vi }{\big(1 - s \cos (\theta + \vi) \big) \big(1 - s \cos (\theta - \vi ) \big)}\; \pt_l (\theta)\, \rd \theta \\[3ex]
\disp \qquad \quad\; = \frac{s}{\pi} \sin \vi \int_0^\pi \frac{\cos \theta - s \cos \vi }{(\cos \theta  - s \cos \vi)^2 + (1 - s^2) \sin ^2 \theta }\; \pt_l (\theta)\, \rd \theta\, .
\ea
\ee
To obtain bounds on $q_l$ of the desired type we need more regularity of $\pt_l$ than just continuity. Let us assume H\"older continuity of $\pt_l$ in the form
\be \label{B.10}
|\pt_l (\theta) - \pt_l(\widetilde{\theta})| \leq C\, |\theta - \widetilde{\theta}|^\al\,  , \qquad \theta , \widetilde{\theta} \in [0\, , \pi]
\ee
for some $C>0$ and some $0< \al < 1$, and, moreover, at the symmetry axis $\{\theta =0\} \cup \{\theta = \pi\}$ the condition
\be \label{B.11}
\ba{rl}
\ba{c}
|\pt_l (\theta) - \pt_l(0)| \leq \Ct\, \theta^{1 + \at} \\[2ex]
|\pt_l (\theta) - \pt_l(\pi)| \leq \Ct\, (\pi -\theta)^{1 + \at} 
\ea & \; , \quad \theta \in [0\, , \pi]
\ea
\ee
for some $\Ct > 0$ and some $0< \at < 1$. To take advantage of the H\"older conditions we rewrite (\ref{B.9}) in the form
\be \label{B.12}
\qb _l(s, \vi) = \frac{s}{\pi} \sin \vi \int_0^\pi \frac{\cos \theta - s \cos \vi }{(\cos \theta  - s \cos \vi)^2 + (1 - s^2) \sin ^2 \theta }\; \big(\pt_l (\theta)- \pt_l (\vi)\big)\, \rd \theta\, ,
\ee
which follows by (\ref{B.8}), (\ref{B.9}) setting $\pt_l \equiv 1$:
$$
\ba{l}
\disp 
2 \sin \vi \int_0^\pi \frac{\cos \theta - s \cos \vi }{(\cos\theta  - s \cos \vi)^2 + (1 - s^2) \sin ^2 \theta }\; \rd \theta\, 
= \int_0^{2 \pi} \frac{\sin (\theta + \vi)}{ 1 - s \cos (\theta + \vi)}\; \rd \theta \\[3ex]
\disp \qquad \qquad \qquad = \frac{1}{s} \ln \big( 1 - s \cos (\theta + \vi) \big) \Big|_0^{2 \pi} = 0\, .
\ea
$$
We prove next the upper bound on the integrand in (\ref{B.12}):
\be \label{B.13}
\frac{|\cos \theta - s \cos \vi| }{(\cos\theta  - s \cos \vi)^2 + (1 - s^2) \sin ^2 \theta } \leq \frac{2}{|\cos \theta - \cos \vi|}\; ,\quad s \in [1/2\, , 1]\, .
\ee
Setting $\cos \theta =: x$, $\cos \vi =: y$, inequality (\ref{B.13}) is equivalent to
\be \label{B.14}
\big| x^2 + s\, y^2 - (1+s)\, x y\big| \leq 2 s \Big( s\, x^2 + s\, y^2 - 2 x y + \frac{1}{s} - s \Big)
\ee
on the square $[-1, 1]\times [-1, 1]$ with side condition $s \in [1/2\, , 1]$. We distinguish the cases (i) $xy \leq 0$ and (ii) $x y > 0$. In case (i) the left-hand side in (\ref{B.14}) satisfies 
$$
0\leq  x^2 + s\, y^2 - (1+s)\, x y \leq s\, x^2 + (1 - s) + s\, y^2 - 2 x y \, ,
$$
and hence (\ref{B.14}). In case (ii) the left-hand side satisfies the upper bound 
$$
\ba{rl}
x^2 + s\, y^2 - (1+s)\, x y & \leq s\, x^2 + (1 - s) + s\, y^2 - 2 x y + (1-s) \\[2ex]
& \disp \leq s\, x^2 + s\, y^2 - 2 x y + (1-s) + \Big( \frac{1}{s} -1 \Big) ,
\ea
$$
and hence one half of (\ref{B.14}). The other half requires a lower bound on the left-hand side. By
$$
(1 + s)\,  xy \leq x^2 + \Big(\frac{1+s}{2}\Big)^2 y^2
$$
one obtains
$$
\ba{l}
\disp -\big( x^2 + s\, y^2 - (1+s)\, x y\big) \leq \Big[\Big(\frac{1+ s}{2}\Big)^2 - s\Big] y^2 \leq 2 (1 -s)^2 y^2 \\[2ex]
\disp \qquad \qquad \leq 2 s \Big[ (1 -s)\Big( \frac{1 }{s} - 1\Big) y^2 + (x- y)^2 \Big] \\[2ex]
\disp \qquad \qquad \leq 2 s \Big[\Big( \frac{1 -s}{s} - (1-s)\Big) y^2 + x^2 + y^2 - 2 x y \Big] \\[2ex]
\disp \qquad \qquad \leq 2 s\Big(s\, x^2 + s\, y^2 - 2 x y + (1-s)x^2 + \frac{1- s}{s} y^2 \Big) \\[2ex]
\disp \qquad \qquad \leq 2 s\Big(s\, x^2 + s\, y^2 - 2 x y + (1-s) + \frac{1- s}{s} \Big) ,
\ea
$$
which completes the proof of (\ref{B.14}). Finally we prove
\be \label{B.15}
\frac{2}{\pi} \int_0^\pi \frac{| \pt_l (\theta) - \pt_l( \vi)|}{|\cos \theta - \cos \vi|} \, \rd \theta < K\, , \quad \vi \in [0\, , \pi]
\ee
with some constant $K$ depending on the parameters that characterize the H\"older conditions (\ref{B.10}) and (\ref{B.11}). We distinguish again several cases, viz., (i) $\vi = \pi/2$, (ii) $\vi \in [0\, , \pi/2)$, and (iii) $\vi \in (\pi/2\, , \pi]$. Case (i) is easily done by (\ref{B.10}):
$$
\ba{l}
\disp \int_0^\pi \frac{|\pt_l (\theta) - \pt_l( \pi/2)|}{|\cos \theta |} \, \rd \theta \leq C \int_0^\pi \frac{|\theta - \pi/2|^\al}{|\sin (\theta - \pi/2)|}\, \rd \theta \leq 2\, C\int_0^\pi |\theta -\pi/2|^{\al -1}\, \rd \theta \\[3ex]
\disp \qquad \qquad = \frac{4 C}{\al} \Big(\frac{\pi}{2}\Big)^\al ,
\ea
$$
where we made use of $\sin x \geq x/2$ for $x \in [0\, ,\pi/2]$. 

Case (ii) requires a balance between H\"older conditions (\ref{B.10}) and (\ref{B.11}) in the form
$$
\ba{l}
\disp \int_0^\pi \frac{| \pt_l (\theta) - \pt_l( \vi)|}{|\cos \theta - \cos \vi|} \, \rd \theta = \frac{1}{2} \int_0^\pi \frac{|\pt_l (\theta) - \pt_l (0) - \pt_l (\vi) + \pt_l (0)|^{1 - \eps}}{\sin \big((\vi + \theta)/2\big)} \ti \\[3ex]
\disp \qquad \qquad \qquad \qquad \qquad \qquad \quad\ti \frac{| \pt_l (\theta) - \pt_l( \vi)|^\eps}{\big|\sin\big((\vi - \theta)/2\big)\big|} \, \rd \theta \\[4ex]
\disp \qquad \qquad \leq \frac{1}{2}\, \Ct^{1 - \eps}\, C^\eps \int_0^\pi \frac{(\theta^{1 + \at} + \vi^{1 + \at})^{1 - \eps }}{(\theta + \vi)/8} \ti \frac{|\theta - \vi|^{\al \eps}}{|\theta - \vi|/4}\, \rd \theta\, ,
\ea
$$
where we made use of $\sin x \geq x/4$ for $x \in [0\, , 3\pi/4]$. The first factor under the integral is further estimated by 
$$
\ba{l}
\disp \frac{(\theta^{1 + \at} + \vi^{1 + \at})^{1 - \eps }}{\theta + \vi} \leq \frac{(\theta + \vi)^{(1 + \at)(1 - \eps) }}{\theta + \vi} 
\leq (\theta + \vi)^{(1 + \at)(1 - \eps) - 1} \leq \Big(\frac{3 \pi}{2} \Big)^{\at (1 - \eps) - \eps} ,
\ea
$$
provided that $\eps$ is chosen such that
$$
\at (1 - \eps) > \eps > 0 \, .
$$
The integral is then easily done with the result
$$
\ba{l}
\disp \int_0^\pi \frac{| \pt_l (\theta) - \pt_l( \vi)|}{|\cos \theta - \cos \vi|} \, \rd \theta \leq 16\, \Ct^{1 - \eps}\, C^\eps \Big(\frac{ 3 \pi}{2}\Big)^{\at (1 - \eps) - \eps} 
\int_0^\pi |\theta - \vi|^{\al \eps - 1}\, \rd \theta \\[3ex]
\disp \qquad \qquad \qquad \qquad \qquad \leq 32\, \Ct^{1 - \eps}\, C^\eps \Big(\frac{ 3 \pi}{2}\Big)^{\at (1 - \eps) - \eps}\; \frac{\pi^{\al \eps}}{ \al \eps}\, . 
\ea
$$
Case (iii) is treated along the lines of case (ii) with $\pt_l (0)$ replaced by $\pt_l(\pi)$ and using the estimate $\sin x \geq (\pi - x)/4$ for $x \in[\pi/4\, ,\pi]$:
$$
\ba{l}
\disp \int_0^\pi \frac{| \pt_l (\theta) - \pt_l( \vi)|}{|\cos \theta - \cos \vi|} \, \rd \theta \leq \frac{1}{2}\, \Ct^{1 - \eps}\, C^\eps \int_0^\pi \frac{\big((\pi -\theta)^{1 + \at} + (\pi -\vi)^{1 + \at}\big)^{1 - \eps }}{\sin \big((\theta + \vi)/2\big)} \ti \\[3ex]
\disp \qquad \qquad \qquad \qquad \qquad \qquad \qquad \qquad \quad \ti \frac{|\theta - \vi|^{\al \eps}}{|\sin (\theta - \vi)/2|}\, \rd \theta \\[3ex]
\disp \qquad \qquad \qquad \leq \frac{1}{2}\, \Ct^{1 - \eps}\, C^\eps \int_0^\pi \frac{(2 \pi -\theta - \vi)^{(1 + \at)(1 - \eps) }}{(2 \pi -\theta - \vi)/8} \ti \frac{|\theta - \vi|^{\al \eps}}{|\theta - \vi|/4}\, \rd \theta \\[3ex]
\disp \qquad \qquad \qquad \qquad \qquad \leq 32\, \Ct^{1 - \eps}\, C^\eps \Big(\frac{ 3 \pi}{2}\Big)^{\at (1 - \eps) - \eps}\; \frac{\pi^{\al \eps}}{ \al \eps}\, , 
\ea
$$
which is the same bound as in case (ii). This proves (\ref{B.15}).
Summarizing (\ref{B.12}) (\ref{B.13}), and (\ref{B.15}) we arrive at the bound
\be \label{B.16}
|\qb_l (s, \vi)| \leq K\; \sin \vi \, ,\quad (s, \vi) \in [1/2\, , 1) \ti [0 \, ,\pi]\, .
\ee
Note that 
$$
\frac{|\pt_l (\theta) - \pt_l(\vi)|}{|\cos \theta - \cos \vi|}
$$
represents by (\ref{B.13}), (\ref{B.15}) an integrable bound on the integrand in (\ref{B.12}). Thus, Lebesgue's theorem allows the passage $s \ra 1$ in (\ref{B.12}), which then holds for $s = 1$ as well:
\be \label{B.17}
\qb_l (1, \vi) = \frac{ \sin\vi}{\pi} \int_0^\pi \frac{|\pt_l (\theta) - \pt_l(\vi)|}{|\cos \theta - \cos \vi|}\, \rd \theta\, .
\ee
The other restriction, viz., $s \geq 1/2$, introduced in (\ref{B.13}), can be removed by the observation that the denominator in (\ref{B.12}) is strictly positive for $0\leq s \leq 1/2$:
\be \label{B.17a}
\ba{l}
(\cos \theta - s \cos \vi)^2 + (1 - s)^2 \sin^2 \theta \\[2ex]
\disp \qquad = 1 - s^2 + s^2 (\cos^2 \theta + \cos^2 \vi) - 2 s \cos \theta \cos \vi\\[1ex]
\disp \qquad \qquad = 1 - s^2 + s^2 (\cos \theta - \cos \vi)^2 + 2 s (s -1) \cos \theta \cos \vi \\[1ex]
\disp \qquad \qquad \qquad \geq \frac{3}{4} - \frac{1}{2} = \frac{1}{4} > 0\, .
\ea
\ee
Thus (\ref{B.16}) holds for all $s \in [0\, ,1]$ and undoing the substitution (\ref{B.7b}) yields the desired bound on $\qt_l (r, \vi)$ in $Q_\infty$.
\begin{lem}
 Let $(p_l, q_l)$ be a solution of the system (\ref{2.17}), (\ref{2.18})$_a$ in $A_\infty$ with boundary function $p(\theta)$ satisfying the H\"older conditions (\ref{B.10}), (\ref{B.11}). Then the following inequalities hold: \footnote{The tilde on $q_l$ in (\ref{B.7a}) is here omitted.}
 \be \label{B.18}
 |q_l ( \cdot , \vi) |\leq K\, \sin \vi \qquad \mbox{in } \ Q_\infty 
 \ee 
 and 
 \be \label{B.19}
 \int_{Q_\infty} \bigg(\frac{q_l}{r \sin \vi}\bigg)^2\, \rd \mu_\beta^\eta \leq \Big(\frac{\pi}{2}\Big)^{1- \beta}\frac{2 K^2}{\eta (1 - \beta)}\, , \qquad 0 \leq \beta < 1\, , \quad \eta >0 
  \ee
 for some constant $K> 0$.
\end{lem}
The bound (\ref{B.19}) is an immediate consequence of the basic estimate (\ref{B.18}).

\vspace{.5em}
\noi{\em Remark:} The mere H\"older condition (\ref{B.10}) still allows a bound of type   
$$
\sin\vi \int_0^\pi \frac{|\pt_l (\theta) - \pt_l(\vi)|}{|\cos \theta - \cos \vi|}\, \rd \theta < \widetilde{K}\, , \quad \vi \in [0, \pi]
$$
instead of (\ref{B.15}). This is enough to establish existence of $q_l$ up to and on the boundary but not enough for the bound (\ref{B.18}).
Generally it is well known that H\"older continuity of the boundary data, which are continuously assumed by one variable of the conjugate pair, is enough to guarantee existence (and H\"older continuity) of the other variable up to and on the boundary (see, e.g., Privaloff's theorem \cite{BJS64}).
\vspace{.5em}

In the bounded domain $A_R$
the boundary conditions at $S_1$ and $S_R$ are not completely independent. Writing the system (\ref{2.17}) in polar coordinates,
\be \label{B.22a}
\left.
\begin{array}{rl}
\disp \pa_r p_l - \frac{1}{r} \pa_\vi q_l & = 0\, ,  \\[1ex]
\disp \frac{1}{r} \pa_\vi p_l + \pa_r q_l & = 0\, ,  
\end{array} \right\} 
\end{equation}
we find by integration of (\ref{B.22a})$_1$ over $[-\pi, \pi]$ that
$\pa_r\int_{-\pi}^\pi p_l\, \rd \vi = 0$ and hence 
\be \label{B.22b}
\int_{-\pi}^\pi p_l (1, \vi)\,\rd \vi = \int_{-\pi}^\pi p_l (R, \vi)\, \rd \vi\, .
\ee
Inserting the boundary conditions (\ref{2.15}) and (\ref{2.16}), condition (\ref{B.22b}) can be met by a suitable choice of $p_0 =: p_{l,0}$.

A representation of $q_l$, analogous to (\ref{B.7}), is given by (see \cite{B79}, p.\ 398f.)
$$
\ba{l}
\disp q_l (z) = \frac{1}{ 2 \pi} \int_0^{2\pi} \Im\, \bigg\{\frac{R\,e^{i \theta} + z}{R\, e^{i \theta} - z} - H(R\, e^{i \theta} , z) \bigg\}\, p_l (R \, e^{i\theta})\, \rd \theta \\[3ex]
\disp \qquad \quad - \frac{1}{ 2 \pi} \int_0^{2\pi} \Im\, \bigg\{\frac{e^{i \theta} + z}{e^{i \theta} - z} - H(e^{i \theta} , z) \bigg\}\, p_l (e^{i\theta})\, \rd \theta\, , 
\ea
$$
where $H(w, z)$ has the series representation
$$
H(w,z) = 2 \sum_{k=1}^\infty \Big(\frac{w}{R^{2 k} z - w} - \frac{z}{R^{2 k} w - z}\Big) =  2 \sum_{k=1}^\infty \frac{1}{R^{2 k} - 1}\bigg(\Big(\frac{w}{z }\Big)^k  - \Big(\frac{z}{w}\Big)^k \bigg) .
$$
The series are absolutely converging for $1/R^2 < |z/w| < R^2$ and continuous on $\pa A_R \ti \overline{A_R}$. By $z = r e^{i \vi}$, $q_l (z) =: \qt_l (r, \vi)$, $p_l (e^{i \theta} ) =\pt_l (\theta)$, and $p_l (R\, e^{i \theta} ) =\ph_l (\theta)$ one obtains in polar coordinates $(r, \vi)$:
\be \label{B.20}
\ba{l}
\disp \qt _l(r, \vi) = \frac{1}{\pi} \int_0^{2\pi} \bigg[
\frac{r}{R}\, \frac{\sin (\vi - \theta)}{1 - 2 (r/R) \cos (\vi - \theta) + (r/R)^2} \\[3ex]
\disp \qquad \qquad \qquad \qquad + \sum_{k=1}^\infty \frac{1}{R^{2 k} - 1}\bigg(\Big(\frac{r}{R}\Big)^k  + \Big(\frac{R}{r}\Big)^k \bigg) \sin k (\vi - \theta) \bigg]\,\ph_l (\theta)\, \rd \theta \\[3ex]
\disp \qquad \qquad - \frac{1}{\pi} \int_0^{2\pi} \bigg[
r\, \frac{\sin (\vi - \theta)}{1 - 2 r \cos (\vi - \theta) + r^2} \\[3ex]
\disp \qquad \qquad \qquad \qquad \quad + \sum_{k=1}^\infty \frac{1}{R^{2 k} - 1} \big(r^k  + r^{-k} \big) \sin k (\vi - \theta) \bigg]\,\pt_l (\theta)\, \rd \theta \, .
\ea
\ee
The first terms under the integrals are analogous to (\ref{B.7a}) and corresponding estimates can by and large be taken over; the second terms represent harmonic functions, which much simpler can be estimated. Let us consider the first term in the second integral first. By the substitution (\ref{B.7b}) and by reducing the integral as in (\ref{B.9}) to the interval $(0\, ,\pi)$ we arrive precisely at (\ref{B.12}). Thus the bound (\ref{B.16}) applies and in the case that 
$$
s_R := \frac{2 \, R}{1 + R^2} \geq \frac{1}{2}
$$
applies, we are done; otherwise the argument based on (\ref{B.17a}) applies and we are done, too. Undoing the substitution (\ref{B.7b}) we arrive at (\ref{B.18}) with a constant $K$ depending on the constants that characterize the H\"older conditions (\ref{B.10}), (\ref{B.11}). 

When applying the substitution 
$$
\widehat{s} := \frac{2\, r/R}{1 + (r/R)^2}
$$
to the first term of the first integral in (\ref{B.20}) and performing the reduction to the interval $(0\, ,\pi)$ we arrive again at an integral of type (\ref{B.12}) with $\pt_l$, however, replaced by $\ph_l$. Under the conditions (\ref{B.10}), (\ref{B.11}) for $\ph_l$ one obtains as before the bound (\ref{B.16}) and then (\ref{B.18}).

The second terms under the integrals in (\ref{B.20}) allow sufficient estimates just by (\ref{B.10}). We demonstrate this for the term in the second integral. Reducing the integral to the interval $(0\, ,\pi)$ and taking advantage of the zero-mean property of the integrand one obtains:   
\be \label{B.21}
\ba{l}
\disp \frac{1}{\pi} \int_0^{2\pi} \bigg[
\sum_{k=1}^\infty \frac{r^k + r^{-k}}{R^{2 k} - 1}\,  \sin k (\vi - \theta) \bigg]\,\pt_l (\theta)\, \rd \theta \\[3ex]
\disp \qquad = \frac{1}{\pi} \sum_{k=1}^\infty \frac{r^k + r^{-k}}{R^{2 k} - 1} \int_0^{ \pi} \big( \sin k (\vi - \theta) + \sin k (\vi + \theta)\big)\,\big(\pt_l (\theta) - \pt_l (\vi)\big)\, \rd \theta \\[3ex]
\disp \qquad \qquad = \frac{2}{\pi} \sum_{k=1}^\infty \frac{r^k + r^{-k}}{(R^{ k} - 1)(R^k  + 1)} \int_0^{\pi}  \sin k \vi \, \cos k \theta\,\big(\pt_l (\theta) - \pt_l (\vi)\big)\, \rd \theta\, .
\ea
\ee
Thus by $|\sin k \vi| \leq k \sin \vi$ on $[0\, ,\pi]$ and by
$$
\sum_{k = 1}^\infty \frac{k}{R^k - 1} \leq \frac{1}{R -1} \sum_{k = 1}^\infty \frac{k}{R^{k -1} } = \frac{1}{R - 1}\, \frac{1}{(1 - 1/R)^2} = \frac{R^2}{(R-1)^3}\, ,
$$
(\ref{B.21}) allows the estimate:
$$
\ba{l}
\disp |(\ref{B.21})| \leq \frac{2}{\pi} \sin \vi\, \bigg(\sum_{k =1}^\infty \frac{k}{R^k - 1} \bigg) C \int_0^\pi |\theta - \vi|^\al \,\rd \theta \\[3ex]
\disp \qquad \qquad \leq 2 \, C\, \pi^\al \frac{R^2}{(R-1)^3}\, \sin \vi \leq 8 \, C\, \pi^\al \sin \vi \quad \mbox{ for }\; R \geq 2\, .
\ea
$$
The following lemma summarizes these estimates in the bounded domain $A_R$.
\begin{lem}
 Let $(p_l, q_l)$ be a solution of the system (\ref{2.17}), (\ref{2.18}) in $A_R$ with boundary function $p(\theta)$ satisfying condition (\ref{B.22b}) and 
 the H\"older conditions (\ref{B.10}), (\ref{B.11}) at $\pa A_R$. Then $q_l$ satisfies the estimates (\ref{B.18}) and (\ref{B.19}) (with $Q_\infty$ replaced by $Q_R$) for some constant $K$ that depends on the H\"older constants but does not depend on $R \geq 2$.
 \end{lem}
\sect{Estimate of the Fredholm constant $\widetilde{C}_F$}
The purpose of this appendix is to present a bound on $\Ct_F$
that is independent of the radius $R$ in $Q_R$ and independent of properties of the coefficient $\aaa$ other than the bounds (\ref{4.1}). To this end we estimate the spectrum of $\FCt$ from below by a positive constant. The idea of proof is to consider the parabolic problem associated to the (adjoint) eigenvalue problem and to construct a suitable supersolution, whose decay constant provides the desired spectral constant.

The spectrum $\Sigma_\KC$ of a compact operator $\KC$ is well known: it consists of a countable set of eigenvalues $\mu \in \cpl$ with finite multiplicity clustering at most at zero, which in our case is not an eigenvalue (cf.\ proposition 3.3). Let now $\KC_M$ be as defined by eq.\ (\ref{4.31aa}) on the Hilbert space $\LC_{\al + \ga}$, let $\KCt := \frac{1}{r}\, \KC_M \frac{1}{r}$ and $\FCt_M = \IC - M\, \KCt_M$. With
\be \label{C.1}
\KCt_M w = \mu\, w \qquad \Longleftrightarrow \qquad \FCt_M w = (1 -M \mu)\, w
\ee
we find for the constant $\Ct _F$ according to (\ref{4.32a})
\be \label{C.2}
\Ct_F^{-1} = \inf_{\mu \in \Sigma_{\KCt_M}} \big\{| 1 - \mu M|\big\}\, .
\ee
By inversion and with $\la =:\mu^{-1}$ (\ref{C.1})$_a$ takes the form 
$$
r \,\KC_M^{-1} r \, w = \la \, w\,,
$$
or, with $v=: r\, w$, 
$$
\KC^{-1}_M v = \la \, v/r^2\, ,
$$
which corresponds in a weak formulation to eq.\ (\ref{4.3}) with right-hand side replaced by 
$$
(\la - M) \int_{Q_R} v \, \psi\, \frac{1}{r^2}\, \rd \mu_\al^\eta\, .
$$
Finally, substituting $v=\vih^\al r^\eta u$ we arrive at eq.\ (\ref{4.2}) with right-hand side 
$$
(\la - M) \int_{Q_R} u \, \psi\, \frac{1}{r^2}\, \rd \mu\, .
$$
With $u \in \HC_0 (Q_R)$, $\psi \in C^1(\QRb)$ this is the weak formulation (cf.\ section 2) of the boundary value problem 
\be \label{C.3}
\left.
\ba{l}
\disp -\frac{1}{r} \pa_r (r \pa_r u) - \frac{1}{r^2} \pa_\vi^2 u - \frac{1}{r} \pa_r (b_r u) - \frac{1}{r^2} \pa_\vi (b_\vi u) = \frac{1}{r^2} (\la - M) u\, , \\[2ex]
\disp \quad \Big( \pa_r u + b_r \frac{u}{r} \Big)\Big|_{r = 1, R} = 0\; , \qquad u\big|_{\vi = 0, \pi} = 0\, ,
\ea
\right\}
\ee
where we have introduced  the abbreviations 
\be \label{C.4}
b_r := \frac{a_\zeta}{\tan \vi} + a_\rho\;, \qquad 
b_\vi := \frac{a_\rho}{\tan \vi} - a_\zeta \, .
\ee
It turns out that the adjoint problem is easier to handle:
let $L^2_c (Q_R)$ be the complex Hilbert space with scalar product $\langle u, w\rangle  := \int_{Q_R} u\, \overline{w}\, r \rd r \rd \vi$, then 
\be \label{C.5}
\left.
\ba{l}
\disp -\frac{1}{r} \pa_r (r \pa_r w) - \frac{1}{r^2} \pa_\vi^2 w + \frac{1}{r} b_r\, \pa_r w + \frac{1}{r^2} b_\vi\, \pa_\vi w = \frac{1}{r^2} (\lab - M) w \, ,\\[2ex]
\disp \quad \pa_r w \big|_{r = 1, R} = 0\; , \qquad w\big|_{\vi = 0, \pi} = 0\, ,
\ea
\right\}
\ee
denotes the adjoint boundary value problem with eigenvalue $\lab$.

Next we embed the elliptic problem (\ref{C.5}) trivially in a parabolic setting by $W (r, \vi, t) := w(r, \vi) \, e^{-(\lab - M)t}$, which is a solution of the problem 
\be \label{C.6}
\left.
\ba{l}
\disp \pa_t W - r \pa_r (r \pa_r W) - \pa_\vi^2 W + r b_r\, \pa_r W + b_\vi\, \pa_\vi W = 0 \, ,\\[2ex]
\disp \quad \pa_r W \big|_{r = 1, R} = 0\; , \quad W\big|_{\vi = 0, \pi} = 0\; ,\quad W\big|_{t =0} = w 
\ea
\right\}
\ee
in the cylinder $Q_R \ti [0, T)$, $T>0$. (\ref{C.6}) is a real equation, thus with $\la  - M =: \varrho + i \sigma$ the real part of $W$,
\be \label{C.7}
\ba{l}
\disp \RC:= \Re W = \Re \Big( \big(\Re w + i \, \Im w \big) e^{-(\varrho - i \si) t} \Big) \\[2ex]
\disp \qquad \qquad \;\; = e^{- \varrho t} \big(\Re w \cos  \si t - \Im w \sin \si t \big)\, ,
\ea
\ee
is a solution of problem (\ref{C.6}) as well, whose temporal decay is governed by the real part $\varrho$ of the eigenvalue $\la - M$. A positive lower bound $D$ on $\varrho$,
$$
\varrho \geq D \qquad \Longleftrightarrow \qquad \Re \la \geq M + D\, ,
$$
then yields the desired positive lower bound in (\ref{C.2}):
\be \label{C.8}
\Big|1 - \frac{M}{\la}\Big| \geq  1 - \frac{M}{|\la|} \geq 1 - \frac{M}{\Re \la} \geq 1 - \frac{M}{M + D} = \frac{D}{M + D} > 0\, .
\ee

The idea for a lower bound on $\varrho$ lies in the construction of a suitable supersolution\footnote{A supersolution of eq.\ (\ref{C.6})$_1$ is a solution of the inequality that arises from (\ref{C.6})$_1$ by replacing ``$=$'' by ``$\leq$''.} $\Fb$ with decay constant $D$, which, by means of a parabolic maximum principle, bounds $\RC$ for all time once it does so for $t= 0$ and if boundary values are under control. This clearly implies that $\Fb$ cannot shrink faster than $\RC$ and hence $\varrho \geq D$.

The maximum principle must be applicable to weak solutions of equations with merely bounded coefficients. We refer here to theorem 6.25 and corollary 6.26 in \cite{L96}, which apply to the following setting adapted to our situation: 

Let $Q_{\eps ,R}$ be the restricted rectangle $(1,R)\ti (\eps, \pi - \eps )$ for some $0 < \eps < \pi/2$, where the coefficients (\ref{C.4}) are then obviously bounded, and let $Q_{\eps, R}^T$ be the restricted cylinder $Q_{\eps, R}  \ti [0, T)$. Let
$$
C_0^1 (Q_{\eps, R}^T) := \big\{ \Psi \in C^1 (\overline{Q_{\eps, R}^T}) : \Psi = 0 \, \mbox{ on } \,
(1, R) \ti (\{\vi = \eps\} \cup \{\vi = \pi - \eps\})
\ti [0, T] \big\}
$$
be the set of space-time test functions and
$$
\VC := \big\{ \RC \in L^2 (Q_{\eps, R}^T ) : \| \RC \|_V  < \infty \big\}
$$
with
$$
\| \RC \|_V^2 := \int_{Q_{\eps, R}^T} |\na \RC |^2 \, \rd \mu  \rd t + \sup_{t \in (0,T)} \int_{Q_{\eps, R} \ti \{t\}} 
|\RC |^2 \, \rd \mu 
$$ 
%\end{comment}
be the set of admissible solutions of the following weak formulation of eq.\ (\ref{C.6}):
\be \label{C.9}
\ba{l}
\disp - \int_{Q_{\eps,R}^T} \RC \, \pa_t \Psi\, \rd \mu  \rd t + \int_{Q_{\eps,R}} \RC(\cdot , \cdot, T) \,  \Psi(\cdot, \cdot, T)\, \rd \mu \\[3ex]
\disp \qquad + \int_{Q_{\eps,R}^T} \Big(\pa_r \RC \, \pa_r \Psi + \frac{1}{r^2}\, \pa_\vi \RC\, \pa_\vi \Psi + \frac{1}{r}\, b_r\, \pa_r \RC\, \Psi + \frac{1}{r^2}\, b_\vi\, \pa_\vi \RC\, \Psi \Big) \rd \mu \rd t \\[2ex]
\disp \qquad \qquad = \int_{Q_{\eps,R}} \Re w \, \Psi(\cdot, \cdot, 0)\, \rd \mu\, .
\ea
\ee
Equation (\ref{C.9}) corresponds to the weak formulation in (\cite{L96}, p.\! 100 above) with coefficients $b^i = c^0 = 0$. In this situation we have (see \cite{L96}, p.\! 128):

\vspace{.3em}
\noi{\bf Proposition C.1 (weak parabolic maximum principle).}
{\em
Let $\RC \in \VC$ be a solution of (\ref{C.9}) as inequality, i.e.\! ``$\, =$'' is replaced by ``$\,\leq$'', for any nonnegative $\Psi \in C_0^1 (Q_{\eps,R}^T)$, then 
\be \label{C.10}
\sup_{Q_{\eps,R}^T} \RC \leq \sup_{\pa Q_{\eps, R}^T} \RC\, ,
\ee
where $\pa Q_{\eps, R}^T$ denotes the parabolic boundary $Q_{\eps, R} \ti \{t=0\} \cup \pa Q_{\eps, R} \ti [0,T)$.
 } 
 \vspace{.3em}
 
To apply proposition C.1 to our situation let us note that by assumption $w$ is a weak solution of eq. (\ref{C.5}) in $Q_R$. In the restricted domain $Q_{\eps, R}$, however,
standard elliptic regularity theory 
can be applied, which yields (at least) $H^2$-regularity  of $w$ and hence $\RC$. Moreover, applying theorem 15.1 in \cite{LU68} we obtain H\"older-continuity (up to the boundary) of $\RC$ and its derivarives for some $0< \al < 1$: 
\be \label{C.11}
\RC \in H^2 (Q_{\eps, R}) \cap C^{1,\al} (\overline{Q_{\eps,R}})\, .
\ee
This implies $\RC \in \VC$ and eq.\ (\ref{C.9}) can be obtained from eq.\ (\ref{C.6}) by multiplying with $\Psi/r^2$ and integration by parts. Note, moreover, that (\ref{C.11}) implies boundary conditions (including the Neumann condition) to hold in the classical sense.

The choice of the supersolution $\Fb$ has to be made with some care to guarantee the boundary control, i.e.\ $|\RC| \leq \Fb$ on $\pa Q_{\eps,R} \ti [0, T)$, which then by (\ref{C.10}) can be extended to all of $Q_{\eps, R}^T$. The (spatial) boundary components $(\{r = 1\} \cup \{r = R\}) \ti (\eps, \pi - \eps) =: \Gamma_R$ and $(1, R) \ti (\{\vi = \eps\} \cup \{\vi = \pi - \eps\}) =: \Gamma_\eps$ cause different problems, which affect the construction of the radial part $f(r)$ and of the angular part $g(\vi)$ in the ansatz for $\Fb$:
\be \label{C.12}
\Fb (r, \vi, t) = A\, f(r)\, g(\vi)\, e^{-D t} ,\qquad A\, ,D >0\, .
\ee

The solution $\RC$ of problem (\ref{C.6}) satisfies a (zero-) boundary condition on $\Gamma_0$ but not on $\Gamma_\eps$, $\eps >0$; however, it exhibits a certain asymptotic behaviour for $\vi \ra 0, \pi$ governed by the singular part in the coefficient $b_\vi$. The function $g(\vi)$ is constructed such that it exhibits a ``better'' asymptotics, which allows then control at $\Gamma_\eps$ for small $\eps >0$. Likewise at $\Gamma_R$ the Neumann boundary condition does not yield an obvious boundary control. Here we choose a slightly convex radial function $f(r)$ that has the property that any boundary point $\in  \Gamma_R$ with $\Fb = \RC$ contradicts the maximum principle.

To be useful, $- \Fb$ has to be a solution of eq.\ (\ref{C.9}) as inequality in the sense of proposition C.1. In terms of classical functions it is sufficient to construct a piecewise-$C^2$ supersolution of eq.\
(\ref{C.6})$_1$ that is globally $C^1$; integration by parts is then allowed and yields the weak formulation (\ref{C.9}).

Inserting the ansatz (\ref{C.12}) into (\ref{C.6})$_1$ leads to the inequality 
$$
D\, f g + r^2\Big(\frac{1}{r} \pa_r \big(r \pa_r f\big) - \frac{b_r}{r} \, \pa_r f \Big) g + f \big(\pa_\vi^2 g - b_\vi \pa_\vi g \big) \leq 0\, ,
$$
which suggests to consider the following two 1-dimensional subproblems:
\be \label{C.13}
f'' + \frac{1 - b_r}{r}\, f' - \frac{d_r}{r^2}\, f \leq 0 \qquad \mbox{ on }\; (1,R)
\ee
and 
\be \label{C.14}
g'' - b_\vi\, g' + d_\vi\, g \leq 0 \qquad \mbox{ on }\; (0,\pi)\, ,
\ee
where we have assumed that 
\be \label{C.15}
D = d_\vi - d_r\, ,\quad d_r,  d_\vi >0\, ,\quad f, g \geq 0\, .
\ee
Let us consider inequality (\ref{C.14}) first. Note that further estimate of $b_\vi$ has to observe the sign of $g'$ and, as described above, we would like to keep the singular part of $b_\vi$ on $(0, \pi)$. So, by (\ref{C.4}) and (\ref{4.1}) we make use of 
\be \label{C.16}
\left.
\ba{l}
\disp - b_\vi = - \frac{a_\rho}{\tan \vi} + a_\ze \leq \frac{0.22}{\vi} + 0.73 \leq \frac{1}{2\, \vi} + 1 \quad \mbox{ on }\; \Big(0, \frac{\pi}{2} \Big) ,\\[3ex]
\disp  - b_\vi \geq - \frac{0.22}{\vih} - 0.73 \geq - \frac{1}{2\, \vih} - 1 \quad \mbox{ on }\; \Big(\frac{\pi}{2} , \pi \Big) 
\ea \right\}
\ee
and look for solutions of the problem
$$
\left. \ba{c}
\disp g''_1 + \Big(1 + \frac{1}{2 \vi}\Big) g'_1 + d_\vi \, g_1 = 0\, ,\quad g'_1 \geq 0 \quad \mbox{ on }\; \Big(0, \frac{\pi}{2}\Big) , \\[2ex]
\disp g_1 (0) = 0\, .
\ea \right\}
$$
An analytic solution can be given in terms of Kummer's function $M(\cdot , \cdot , \cdot)$ (see \cite{AS72}, p.\ 504): 
\be \label{C.17} \disp
g_1 (\vi) = e^{-(1 + \ka )\vi/2}\, \vi^{1/2}\, M\Big(\frac{3}{4} + \frac{1}{4 \ka}\, , \frac{3}{2}\, , \ka\, \vi \Big)
\ee
with $\ka = (1 - 4 d_\vi)^{1/2}$. The condition $g'_1  (\pi/2) = 0$ fixes $d_\vi$ at 
\be \label{C.18}
d_\vi \approx 0.1945\, ,
\ee
and we have clearly $g'_1 \geq 0$ on $(0, \pi/2)$ (see Fig.\ 6). By reflection at $\vi = \pi/2$ we obtain a function $g_2 (\vi) := g_1 (\pi - \vi) = g_1 (\vih)$ on $(\pi/2, \pi)$ satisfying the problem
$$
\left. \ba{c}
\disp g''_2 - \Big(1 + \frac{1}{2 \vih}\Big) g'_2 + d_\vi \, g_2 = 0\, ,\quad g'_2 \leq 0 \quad \mbox{ on }\; \Big(\frac{\pi}{2} , \pi\Big) , \\[2ex]
\disp g_2 (\pi) = g'_2 \Big(\frac{\pi}{2}\Big) = 0\, .
\ea \right\}
$$
Thus by (\ref{C.16}), the function
$$
g:= \left\{
\ba{l}
g_1 \quad \mbox{ on }\; (0, \pi/2) \\[.5ex]
g_2 \quad \mbox{ on }\; (\pi/2, \pi)
\ea \right.
$$
is a $C^1$-solution of inequality (\ref{C.14}), vanishing  at the boundaries, and with decay constant (\ref{C.18}). 
\begin{figure}
\begin{center}
\includegraphics[width=0.5\textwidth]{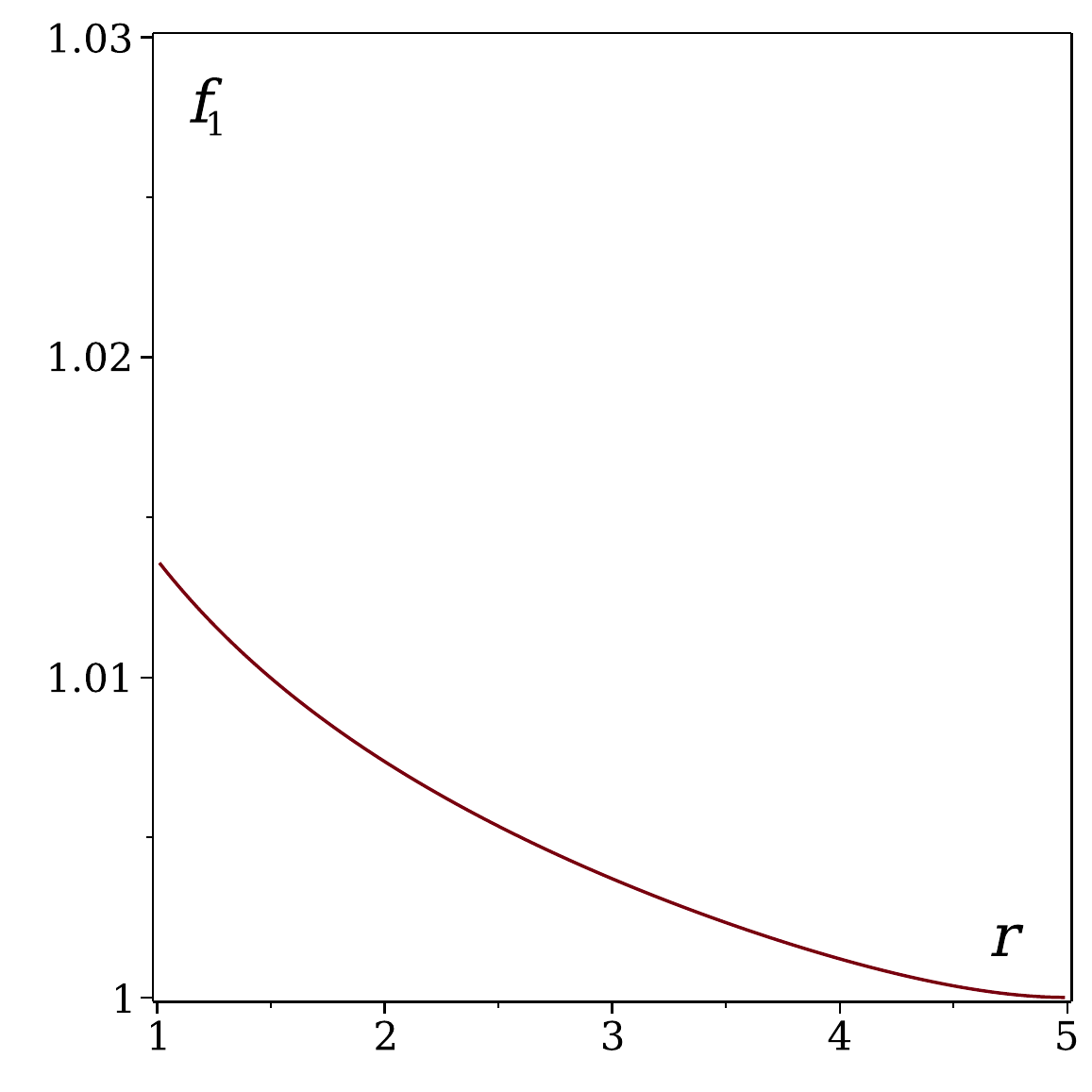}\includegraphics[width=0.5\textwidth]{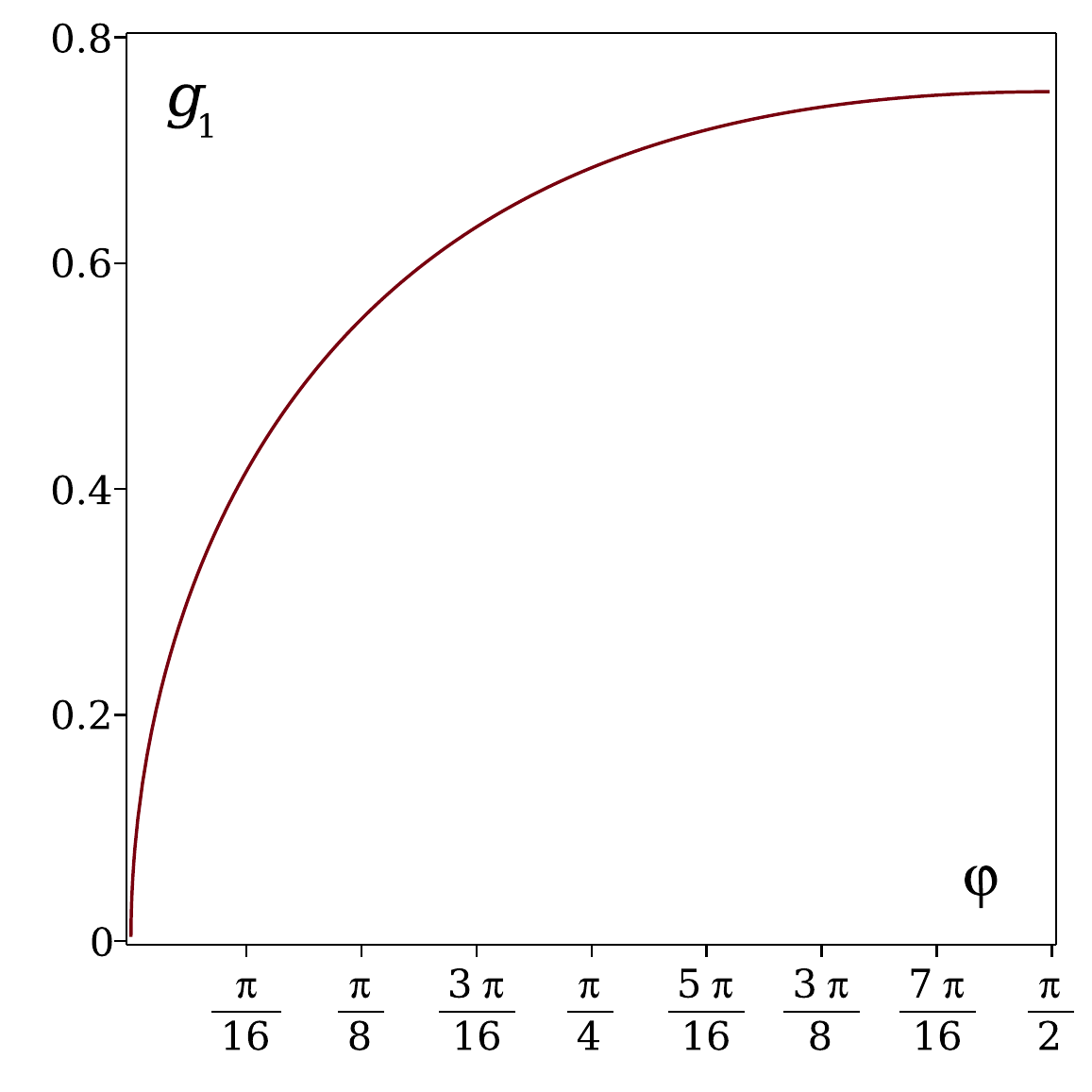}
%\mbox{\phantom{...}} % Das hier ist natürlich ein ziemlich dreckiger Trick, aber sonst kleben die Grafiken zu nahe aufeinander und die Achsengeschriftung von rechts kann nicht mehr eindeutig der rechten Grafik zugeordnet werden.
%\resizebox{0.7\textwidth}{!}{\input{Kummer.eps}}
\caption{Graphs of the functions $f_1: (1\, , R_0) \rightarrow\real_+$ with $R_0 = 5$, $c_\eps =10$, $d_r = 0.09$ and $g_1 : (0\, , \pi/2) \rightarrow \real_+$ with $d_\vi = 0.1945$.}
\end{center}
\end{figure}

Concerning the radial part $f(r)$ we make use of the estimate 
\be \label{C.19}
|b_r| = \Big| \frac{a_\ze }{ \tan \vi}+ a_\rho \Big| \leq \frac{1}{\eps} + 1 =: c_\eps  \quad \mbox{ on }\; Q_{\eps, R}
\ee
and look for solutions of the problem 
\be \label{C.20}
\left. \ba{c}
\disp f''_1 + \frac{1 - c_\eps }{r}\, f'_1 - \frac{d_r}{r^2} \, f_1 = 0\, ,\quad f'_1 \leq 0 \quad \mbox{ on }\; (1, R_0)\, , \\[2ex]
\disp f_1 (R_0) = 1\, ,\quad  f'_1 (R_0) = 0\, ,
\ea \right\}
\ee
where $R_0$ is some radius with $1 < R_0 < R$. The solution of the initial value problem (\ref{C.20}) is elementary:
$$
f_1 = \frac{\al_+ (r/R_0)^{\al_-} - \al_- (r/R_0)^{\al_+}}{\al_+ - \al_-} \, ,\qquad \al_\pm = \frac{1}{2} \big(c_\eps \pm ( c_\eps^2 + 4 d_r )^{1/2} \big)
$$
with
$$
f'_1 (r) = \frac{\al_+ \al_-}{\al_+ - \al_-} \, \frac{1}{r} \bigg(\Big(\frac{R_0}{r}\Big)^{-\al_-} - \Big( \frac{r}{R_0}\Big)^{\al_+}\bigg) \leq 0
$$
and, in particular, $f'_1 (1) < 0$ (see Fig.\ 6). Analogously we find $f_2 (r)$ as solution of the problem 
$$
\left. \ba{c}
\disp f''_2 + \frac{1 + c_\eps }{r}\, f'_2 - \frac{d_r}{r^2} \, f_2 = 0\, ,\quad f'_2 \geq 0 \quad \mbox{ on }\; (R_0, R)\, , \\[2ex]
\disp f_2 (R_0) = 1\, ,\quad  f'_2 (R_0) = 0\, .
\ea \right\}
$$
$f_1$ and $f_2$ constitute a convex $C^1$-function $f$ on $(1, R)$ satisfying inequality (\ref{C.13}) and 
\be \label{C.21}
f'(1) < 0\, ,\qquad f'(R) > 0\, .
\ee
Note that $d_r$ is here a free constant and by the choice $d_r := 0.09$ and by (\ref{C.18}) we find 
\be \label{C.22}
D = d_\vi - d_r > 0.1 > 0\, .
\ee
The following lemma proves $D$ to be a lower bound on the eigenvalues of problem (\ref{C.6}) or, equivalently, (\ref{C.5}).
\begin{lem}
Let $w$ be a weak solution of the eigenvalue problem (\ref{C.5}) 
on $Q_R$ with eigenvalue $\lab - M$, and with coefficients $\bb$ given by (\ref{C.4}). Then, the real part $\varrho$ of the eigenvalue has a positive lower bound $D$ given by (\ref{C.22}), which especially holds for any $R >1$ and for any vector field $\aaa$ satisfying the bounds (\ref{4.1}). 
\end{lem}
\noi\textsc{Proof:} Let us fix $Q_R$ with some $R>1 $ and assume that
\be \label{C.23}
D > \varrho\, ,
\ee
which means that the supersolution $\Fb$ shrinks faster than the (time-dependent) eigenfunction $\RC$ on $Q_R$. Let us further assume that $\RC$ has a positive maximum on $Q_R$ (otherwise $\RC$ has a negative minimum and the proof works for $-\RC$ instead of $\RC$). Let us fix an amplitude $A$ in (\ref{C.12}) such that 
\be \label{C.24}
\RC < \Fb \quad \mbox{ on }\; Q_R \ti \{t =0\}
\ee
and choose by (\ref{C.23}) $T > 0$ such that 
\be \label{C.25}
\RC > \Fb \quad \mbox{ at some point  }\; p \in Q_R^T\, .
\ee
At the boundary component $(1, R) \ti (\{\vi = 0\} \cup \{\vi = \pi\}) \ti [0,T)$ both functions $\RC$ and $\Fb$ vanish, however, with different asymptotics. 
With $M(\cdot, \cdot, 0) = 1$ the asymptotic bahaviour of $\Fb$ can directly be read off the representation (\ref{C.17}):
\be \label{C.26}
\Fb \sim \vi^{1/2} \quad \mbox{ for }\; \vi \ra  0\, .
\ee
On the other hand an asymptotic analysis of eq.\ (\ref{6.6})$_1$ yields the equation 
\be \label{C.27}
\RC''_0 - \frac{a_\rho}{\vi}\, \RC'_0 \approx 0 \quad \mbox{ for }\;\vi \ra 0
\ee
for the variable $\RC_0 := \RC (r_0, \vi, t_0)$ with $r = r_0$, $t= t_0$ fixed. Equation (\ref{C.27}) implies the asymptotic behaviour for $\vi \ra 0$: 
\be \label{C.28}
|\RC| \lesssim \vi^\de \quad \mbox{ with } \; \de \geq 1 - 0.22 > 3/4\, .
\ee
Comparing (\ref{C.26})  with (\ref{C.28}) one finds $\eps_0> 0$ such that for any $0 < \eps \leq \eps_0$ we have the desired control:
\be \label{C.29}
|\RC| \leq \Fb \quad \mbox{ on } \; (1, R) \ti \big(\{\vi = \eps\} \cup \{\vi = \pi - \eps \}\big) \ti [0,T)\, .
\ee
The argument works analogously in the case $\vi \ra \pi$, which  we have included into (\ref{C.29}).

At the other boundary component $(\{r = 1\} \cup \{r= R\}) \ti 
(\eps, \pi - \eps) \ti [0,T)$ we distinguish two cases. In the first case we assume 
$$
|\RC| < \Fb \quad \mbox{ on } \; \big(\{r = 1\} \cup \{r= R\}\big) \ti (\eps, \pi - \eps) \ti [0,T)\, ,
$$
which implies together with (\ref{C.24}) and (\ref{C.29}) at the parabolic boundary:
\be \label{C.30}
\RC - \Fb \leq 0 \quad \mbox{ on }\; \pa Q_{\eps, R}^T\, .
\ee
Applying proposition C.1 on $\RC - \Fb$ then yields $\RC \leq \Fb$ on $Q_{\eps,R}^T$ for any $\eps \leq \eps_0$ in contradiction to (\ref{C.25}). Otherwise we have by continuity (see (\ref{C.11})) boundary points $p$ with (say) $r = 1$ such that 
$$
\RC = \Fb \quad \mbox{ at } \; p \in \{r = 1\} \ti (\eps, \pi - \eps) \ti (0,T)\, .
$$
Let $p_0 = (1, \vi_0\, , T_0)$ be the earliest  (i.e., with minimal $T =: T_0$) such point, then inequality (\ref{C.30}) holds on the ``shorter'' boundary $\pa Q_{\eps, R}^{T_0}$ and by proposition C.1 on $Q_{\eps, R}^{T_0}$. This, however, is again a contradiction: the boundary condition (\ref{C.6})$_{2a}$ for $\RC$ and property (\ref{C.21}) of $\Fb$ imply obviously $\RC > \Fb$ in a neighbourhood of $p_0$ in $Q_{\eps, R}^{T_0}\, $.

In conclusion the initial assumption $D > \varrho$ has to be dropped, which proves the lemma.

\qed

By (\ref{C.2}) and (\ref{C.8}) lemma C.1 implies the desired bound on the Fredholm constant $\Ct_F$.

%
%
%      Danksagung
%
%
\begin{comment}
\section*{Conflict of interest}

The author has no competing interests to declare that are relevant to the content of this article.

\section*{Funding}

No funding was received for carrying out this work.

\section*{Data Availability}

Data sharing is not applicable to this article as no datasets were generated or analysed in this work.
\end{comment}

%\begin{acknowledgements}
%The authors would like to thank A.\ Tilgner for substantial help in %preparing the numerical code for the computations in appendix E.  
%\end{acknowledgements}
%
%
%
%
%
%    Literaturverzeichnis
%
%
%
%\section*{References}

%%%
\end{document}